\documentclass[12pt,a4paper]{article}

\usepackage{amsbsy,amsfonts,amsmath,amssymb,amsthm}
\usepackage{array}
\usepackage{bm}
\usepackage{color}
\usepackage{enumerate}
\usepackage{mathrsfs}
\usepackage{latexsym}
\usepackage[colorlinks=true]{hyperref}
\hypersetup{urlcolor=blue, citecolor=blue, linkcolor=blue}

\definecolor{wineRed}{rgb}{0.7,0,0.3}
\definecolor{grandBleu}{rgb}{0,0,0.8}
\definecolor{darkGreen}{rgb}{0,0.4,0}
\definecolor{blueViolet}{rgb}{0.4,0,1.0}
\definecolor{bloodOrange}{rgb}{0.85,0.05,0}
\definecolor{mycolor}{rgb}{0.8,0,0.2}
\definecolor{}{rgb}{0.8,0,0.2}

\usepackage[textsize=small]{todonotes}
\usepackage{cite}

\usepackage{comment}

\DeclareMathAlphabet{\mathpzc}{OT1}{pzc}{m}{it}

\usepackage[textsize=small]{todonotes}
\numberwithin{equation}{section}

\theoremstyle{plain}
\newtheorem{theorem}{Theorem}[section]

\newtheorem{main}{Main Theorem}
\newtheorem{lemma}[theorem]{Lemma}
\newtheorem{proposition}{Proposition}

\theoremstyle{definition}
\newtheorem{definition}[theorem]{Definition}
\newtheorem{remark}{Remark}

\def\Sgn{\mathop{\mathrm{Sgn}}\nolimits}

\begin{document}
\thispagestyle{plain}
\begin{center}
    \textbf{\Large Optimal Control Problems Governed by \\[0.5ex] 1-D {K}obayashi--{W}arren--{C}arter Type Systems}\footnotemark[1]
\end{center}
    \bigskip
\vspace{-0.5ex}
\begin{center}
    \textsc{Harbir Antil} 
    \\[1ex]
    {Department of Mathematical Sciences and the Center for Mathematics and \\ Artificial Intelligence (CMAI), George Mason University, \\ Fairfax, VA 22030, USA}
    \\[0ex]
    ({\ttfamily hantil@gmu.edu})
\end{center}
\begin{center}
    \textsc{Shodai Kubota}
    \\[1ex]
    {Department of Mathematics and Informatics, \\ Graduate School of Science and Engineering, Chiba University, \\ 1-33, Yayoi-cho, Inage-ku, 263-8522, Chiba, Japan}
    \\[0ex]
    ({\ttfamily skubota@chiba-u.jp})
\end{center}
\begin{center}
    \textsc{Ken Shirakawa}
    \\[1ex]
    {Department of Mathematics, Faculty of Education, Chiba University, \\ 1-33, Yayoi-cho, Inage-ku, 263-8522, Chiba, Japan}
    \\[0ex]
    ({\ttfamily sirakawa@faculty.chiba-u.jp})
\end{center}
\begin{center}
\vspace{-1ex}
    \textsc{and}
\vspace{-0.5ex}
\end{center}
\begin{center}
\textsc{Noriaki Yamazaki}
\\[1ex]
{Department of Mathematics, Faculty of Engineering, Kanagawa University, \\ 3-27-1, Rokkakubashi, Kanagawa-ku, Yokohama, 221-8686, Japan}
\\[0ex]
({\ttfamily noriaki@kanagawa-u.ac.jp})
\end{center}
\vspace{-1ex}

\footnotetext[1]{
The work of the third author supported by Grant-in-Aid for Scientific Research (C) No. 16K05224 and No. 20K03672, JSPS. The work of the forth author supported by Grant-in-Aid for Scientific Research (C) No. 20K03665 , JSPS.  In addition, the work of the first and the third authors is partially supported by the Air Force Office of Scientific Research (AFOSR) under Award NO: FA9550-19-1-0036 and NSF grants DMS-1818772 and DMS-1913004.  
\\[1ex]
\\[1ex]
}
\bigskip

\noindent
{\bf Abstract.} This paper is devoted to the study of a class of optimal control problems governed by 1--D Kobayashi--Warren--Carter type systems, which are based on a phase-field model of grain boundary motion, proposed by [Kobayashi et al, Physica D, 140, 141--150, 2000]. The class consists of an optimal control problem for a physically realistic state-system of Kobayashi--Warren--Carter type, and its regularized approximating problems. The results of this paper are stated in three Main Theorems 1--3. The first Main Theorem 1 is concerned with the solvability and continuous dependence for the state-systems. Meanwhile, the second Main Theorem 2 is concerned with the solvability of optimal control problems, and some semi-continuous association in the class of our optimal control problems. Finally, in the third Main Theorem 3, we derive the first order necessary optimality conditions for optimal controls of the regularized approximating problems. By taking the approximating limit, we also derive the optimality conditions for the optimal controls for the physically realistic problem.

\newpage

\section*{Introduction}
Let $ (0, T) $ be a time-interval with a constant $ 0 < T < \infty $, and let $\Omega := (0, 1) \subset \mathbb{R} $ be a one-dimensional spatial domain with a boundary $ \Gamma := \{0, 1\} $. Besides, we set $ Q := (0, T) \times \Omega $ and $ \Sigma := (0, T) \times \Gamma $, and  we define $ H := L^2(\Omega) $ and  $ \mathscr{H} := L^2(0, T; L^2(\Omega)) $  as the base spaces for our problems. 
\medskip

In this paper, we consider a class of optimal control problems governed by the following state-systems, which are denoted by (S)$_\varepsilon $, with $ \varepsilon \geq 0 $:
\medskip

(S)$ _\varepsilon $
\vspace{0ex}
\begin{equation}\label{1}
\left\{ \parbox{11cm}{
    $ \partial_{t} \eta -\partial_{x}^{2} \eta +g(\eta) +\alpha'(\eta) \sqrt{\varepsilon^2 +|\partial_{x} \theta|^2} = M_u u $ \quad in $ Q $,
\\[1ex]
$ \partial_{x} \eta(t, x) = 0 $, \quad $ (t, x) \in \Sigma $, 
\\[1ex]
$ \eta(0, x) = \eta_{0}(x) $, \quad $ x \in \Omega $;
}\right. 
\end{equation}
\begin{equation}\label{2}
\left\{\parbox{11cm}{
    $ \displaystyle \alpha_{0}(t, x) \partial_{t} \theta -\partial_{x} \left( \alpha(\eta) \frac{\partial_{x}\theta}{\sqrt{\varepsilon^2 +|\partial_{x} \theta|^2}} +\nu^{2}\partial_{x}\theta\right) = M_v v $ \quad in $Q$,
\\[1ex]
$ \theta(t, x) = 0 $, \quad $ (t, x) \in \Sigma $, 
\\[1ex]
$ \theta(0, x) = \theta_{0}(x)$, \quad $x \in \Omega $.
}\right. 
\end{equation}
For each $ \varepsilon \geq 0 $, we denote the optimal control problem by (OP)$ _\varepsilon $, and prescribe the problem as follows:
\begin{itemize}
    \item[\textmd{(OP)$_\varepsilon$}]Find a pair of functions $ [u^*, v^*] \in [\mathscr{H}]^2 $, called \emph{optimal control}, which minimizes a cost functional $ \mathcal{J}_\varepsilon = \mathcal{J}_\varepsilon(u, v) $, defined as:
        \begin{align}\label{J}
            \mathcal{J}_\varepsilon : [u, &\,  v] \in [\mathscr{H}]^2 \mapsto \mathcal{J}_\varepsilon(u, v) 
            \nonumber
            \\
            := &~ \frac{M_\eta }{2} \int_0^T |(\eta -\eta_\mathrm{ad})(t)|_H^2 \, dt  +\frac{M_\theta }{2} \int_0^T |(\theta -\theta_\mathrm{ad})(t)|_H^2 \, dt
            \\
            &~ +\frac{M_u }{2} \int_0^T |u(t)|_H^2 \, dt +\frac{M_v }{2} \int_0^T |v(t)|_H^2 \, dt \in [0, \infty),  
            \nonumber
        \end{align}
        where $ [\eta, \theta] \in [\mathscr{H}]^2 $ solves the state-system (S)$_\varepsilon$. 
\end{itemize}
The state-system (S)$_\varepsilon$ is a type of Kobayashi--Warren--Carter system, i.e. it is based on a phase-field model of grain boundary motion, proposed by Kobayashi et al \cite{MR1752970, MR1794359}. The order parameters, $ \eta \in \mathscr{H} $ and $ \theta \in \mathscr{H} $ indicate the \emph{orientation order} and \emph{orientation angle} of the polycrystal body, respectively. Moreover, $ [\eta_0, \theta_0] \in H^1(\Omega) \times H^1_0(\Omega)$ is an \emph{initial pair}, i.e. a pair of initial data of $ [\eta, \theta] $. The \emph{forcing pair} $ [u, v] \in [\mathscr{H}]^2 $ denotes the control variables that can control the profile of solution $ [\eta, \theta] \in [\mathscr{H}]^2 $ to (S)$_\varepsilon$. Additionally, $ 0 < \alpha_0 \in W^{1, \infty}(Q) $ and $ 0 < \alpha \in C^2(\mathbb{R}) $ are given functions to reproduce the mobilities of grain boundary motions. Finally, $ g \in W_\mathrm{loc}^{1, \infty}(\mathbb{R}) $ is a perturbation for the orientation order $ \eta $, and $ \nu > 0 $ is a fixed constant to relax the diffusion of the orientation angle $ \theta $. 
\medskip

In the state-system (S)$_\varepsilon$, the PDE part of the first initial-boundary value problem \eqref{1} is a type of Allen--Cahn equation, so that the forcing term $ u $ can be regarded as a \emph{temperature control} of the grain boundary formation. 
Also, the second problem \eqref{2} is to reproduce crystalline micro-structure of polycrystal, and the case of $ \varepsilon = 0 $ is the closest to the original setting adopted by Kobayashi et al \cite{MR1752970, MR1794359}. Indeed, when $ \varepsilon = 0 $, the quasi-linear diffusion as in \eqref{2} is described in a singular form $ -\partial_x \bigl( \alpha(\eta) \frac{\partial_x \theta}{|\partial_x \theta|} +\nu^2 \partial_x \theta \bigr) $, and it is known that this type of singularity is effective to reproduce the \emph{facet}, i.e. the locally uniform (constant) phase in each oriented grain (cf. \cite{MR1752970, MR1794359,MR2746654,MR2101878,MR2436794,MR2033382,MR2223383,MR3038131,MR3670006,MR3268865,MR2836557,MR1865089,MR1712447,MR3951294}). 
Hence, the systems (S)$_\varepsilon$, for positive $ \varepsilon $, can be said as regularized approximating systems, that are to approach to the physically realistic situation, reproduced by the limiting system (S)$_0$, as $ \varepsilon \downarrow 0 $. 
\medskip

On the other hand, the pair of functions $  [\eta_\mathrm{ad}, \theta_\mathrm{ad}] \in [\mathscr{H}]^2 $, in the optimal control problem (OP)$_\varepsilon$, is a given \emph{admissible target profile} of $ [\eta, \theta] \in [\mathscr{H}]^2 $. Moreover, $ M_\eta \geq 0 $, $ M_\theta \geq 0 $, $ M_u \geq 0 $, and $ M_v \geq 0 $ are fixed constants, that are to adjust the meaning of optimality in the problem (OP)$ _\varepsilon $. 
\medskip

This paper focuses on two issues:
\begin{itemize}
    \item[\textmd{~~$\sharp \, 1)$}]key-properties of the state-systems (S)$ _\varepsilon $, for $ \varepsilon \geq 0 $; 
        \vspace{0ex}
    \item[\textmd{~~$\sharp \, 2)$}]mathematical analysis of the optimal control problem (OP)$ _\varepsilon $, for $ \varepsilon \geq 0 $.
\end{itemize}
With regard to the first issue $\sharp \, 1)$, various singular systems, related to (S)$ _\varepsilon $, have been studied by several authors, e.g. \cite{MR2469586,MR2548486,MR2836555,MR2668289,MR3670006,MR3038131,MR3362773,MR3082861,MR3203495,MR3462536,MR3238848,MR3888636}. In particular, the mathematical theories developed in \cite[Theorems 2.1 and 2.2]{MR2469586} and \cite[Main Theorems 1 and 2]{MR3888636} are applicable for the well-posedness and $ \varepsilon $-dependence of the system (S)$ _\varepsilon $. However, since the previous works dealt with only homogeneous case, i.e., the case of $ [u, v] = [0, 0] $, some extension of the existing theories is needed for the application to our optimal control problem (OP)$_\varepsilon$. Meanwhile, for issue $ \sharp \, 2) $, the important point will be how to compute the G\^{a}teaux differential of the cost $ \mathcal{J}_\varepsilon $. This will be carried via a linearization of the state-system (S)$_\varepsilon$. When $ \varepsilon > 0 $, the problem (OP)$_\varepsilon$ admits sufficient regularity, and we can address the issue $ \sharp\,2) $ by using the standard linearization method. Although such linearization method does not work for the problem (OP)$_0$,  i.e. the case of $ \varepsilon = 0 $, it is possible to obtain some partial results by considering the limit as $\varepsilon \downarrow 0$ for (OP)$_\varepsilon$. 	
\bigskip

Now, based on these, the goal of this paper is to prove three Main Theorems, summarized as follows:
\begin{description}
    \item[{\boldmath Main Theorem 1}]mathematical results concerning the following items:
\item[{\boldmath~~(I-A)(Solvability of state-systems)}]Existence and uniqueness for the state-system (S)$_\varepsilon$, for any $ \varepsilon \geq 0 $.
\item[{\boldmath~~(I-B)(Continuous dependence among state-systems)}]Continuous dependence of solutions to the systems (S)$_\varepsilon$, with respect to $ \varepsilon \geq 0 $. Roughly summarized, the uniform convergence of the  solutions and governing convex energies, under the convergence of $\varepsilon$ to a value $\varepsilon_0 \geq 0$, weak $H^1$-convergence of initial values, and weak $L^2$-convergence of forces (controls).
\end{description}
\begin{description}
\item[{\boldmath Main Theorem 2}]mathematical results concerning the following items:
\item[{\boldmath~~(II-A)(Solvability of optimal control problems)}]Existence for the optimal control problem (OP)$_\varepsilon$, for any $ \varepsilon \geq 0 $.
\item[{\boldmath~~(II-B)($\varepsilon$-dependence of optimal controls)}]Some semi-continuous association between the optimal controls, with respect to $ \varepsilon $.    
\end{description}
\begin{description}
\item[{\boldmath Main Theorem 3}]mathematical results concerning the following items:
        \item[{\boldmath~~(III-A)(Necessary optimality conditions in cases of $ \varepsilon > 0 $)}] Derivation \linebreak of first order necessary optimality conditions for (OP)$_\varepsilon$ via adjoint method.
\item[{\boldmath~~(III-B)(Limiting optimality conditions as ~$\varepsilon \downarrow 0$)}]The limiting adjoint system as $ \varepsilon \downarrow 0 $.   
\end{description}

This paper is organized as follows. Preliminaries are given in Section 2, the auxiliary lemmas are given in Section 2 and the Main Theorems are proved in Sections 5-7, with an appendix in Section 8.

\section{Preliminaries} 

We begin by prescribing the notations used throughout this paper. 
\medskip

\noindent
\underline{\textbf{\textit{Abstract notations.}}}
For an abstract Banach space $ X $, we denote by $ |\cdot|_{X} $ the norm of $ X $, and denote by $ \langle \cdot, \cdot \rangle_X $ the duality pairing between $ X $ and its dual $ X^* $. In particular, when $ X $ is a Hilbert space, we denote by $ (\cdot,\cdot)_{X} $ the inner product of $ X $. For any subset $ A $ of a Banach space $ X $, let $ \chi_A : X \longrightarrow \{0, 1\} $ be the characteristic function of $ A $, i.e.:
    \begin{equation*}
        \chi_A: w \in X \mapsto \chi_A(w) := \begin{cases}
            1, \mbox{ if $ w \in A $,}
            \\[0.5ex]
            0, \mbox{ otherwise.}
        \end{cases}
    \end{equation*}

For two Banach spaces $ X $ and $ Y $,  we denote by $  \mathscr{L}(X; Y)$ the Banach space of bounded linear operators from $ X $ into $ Y $, and in particular, we let $ \mathscr{L}(X) := \mathscr{L}(X; X) $. 

For Banach spaces $ X_1, \dots, X_N $, with $ 1 < N \in \mathbb{N} $, let $ X_1 \times \dots \times X_N $ be the product Banach space endowed with the norm $ |\cdot|_{X_1 \times \cdots \times X_N} := |\cdot|_{X_1} + \cdots +|\cdot|_{X_N} $. However, when all $ X_1, \dots, X_N $ are Hilbert spaces, $ X_1 \times \dots \times X_N $ denotes the product Hilbert space endowed with the inner product $ (\cdot, \cdot)_{X_1 \times \cdots \times X_N} := (\cdot, \cdot)_{X_1} + \cdots +(\cdot, \cdot)_{X_N} $ and the norm $ |\cdot|_{X_1 \times \cdots \times X_N} := \bigl( |\cdot|_{X_1}^2 + \cdots +|\cdot|_{X_N}^2 \bigr)^{\frac{1}{2}} $. In particular, when all $ X_1, \dots,  X_N $ coincide with a Banach space $ Y $, we write:
\begin{equation*}
    [Y]^N := \overbrace{Y \times \cdots \times Y}^{\mbox{$N$ times}}.
\end{equation*}
Additionally, for any transform (operator) $ \mathcal{T} : X \longrightarrow Y $, we let:
\begin{equation*}
    \mathcal{T}[w_1, \dots, w_N] := \bigl[ \mathcal{T} w_1, \dots, \mathcal{T} w_N \bigl] \mbox{ in $ [Y]^N $, \quad for any $ [w_1, \dots, w_N] \in [X]^N $.}
\end{equation*}

\noindent
\underline{\textbf{\textit{Specific notations of this paper.}}} 
As is mentioned in the previous section, let $ (0, T) \subset \mathbb{R}$ be a bounded time-interval with a finite constant $ T > 0 $, and let  $ \Omega := (0, 1) \subset \mathbb{R} $ be a one-dimensional bounded spatial domain. We denote by $ \Gamma $ the boundary $ \partial \Omega = \{0, 1\} $ of $ \Omega $, and we let $ Q := (0, T) \times \Omega $ and $ \Sigma := (0, T) \times \Gamma $. Especially, we denote by $ \partial_t $ and $ \partial_x $ the distributional time-derivative and the distributional spatial-derivative, respectively. Also, the measure theoretical phrases, such as ``a.e.'', ``$dt$'', ``$dx$'', and so on, are all with respect to the Lebesgue measure in each corresponding dimension.
\medskip

On this basis, we define  
\begin{equation*} \begin{cases}
    H := L^2(\Omega) \mbox{ and } \mathscr{H} := L^2(0, T; H),
\\
    V := H^1(\Omega) \mbox{ and } \mathscr{V} := L^2(0, T; V),
\\
    V_0 := H_0^1(\Omega) \mbox{ and } \mathscr{V}_0 := L^2(0, T; V_0).
\end{cases}
\end{equation*}
Also, we identify the Hilbert spaces $ H $ and $ \mathscr{H} $ with their dual spaces. Based on the identifications, we have the following relationships of continuous embeddings:
\begin{equation*}
\begin{cases}
    V \subset H = H^* \subset V^* \mbox{ and } \mathscr{V} \subset \mathscr{H} = \mathscr{H}^* \subset \mathscr{V}^*,
\\
    V_0 \subset H = H^* \subset V_0^* \mbox{ and } \mathscr{V}_0 \subset \mathscr{H} = \mathscr{H}^* \subset \mathscr{V}_0^*,
\end{cases}
\end{equation*}
among the Hilbert spaces $ H $, $ V $, $ V_0 $, $ \mathscr{H} $, $ \mathscr{V} $, and $ \mathscr{V}_0 $, and the respective dual spaces $ H^* $, $ V^* $, $ V_0^* $, $ \mathscr{H}^* $, $ \mathscr{V}^* $, and $ \mathscr{V}_0^* $. Additionally, in this paper, we define the topology of the Hilbert space $ V_0 $ by using the following inner product:
\begin{equation*}
    (w, \tilde{w})_{V_0} := (\partial_x w, \partial_x \tilde{w})_H, \mbox{ for all $ w, \tilde{w} \in V_0 $.}
\end{equation*}
\begin{remark}\label{Rem.Prelim01}
    Due to the one-dimensional embeddings $ V \subset C(\overline{\Omega}) $ and $ V_0 \subset C(\overline{\Omega}) $, it is easily checked that: \begin{equation}\label{emb01}
        \begin{cases}
            \hspace{-2ex}
            \parbox{9cm}{
                \vspace{-2ex}
                \begin{itemize}
                    \item if $ \check{\mu} \in H $ and $ \check{p} \in V $, then $ \check{\mu} \check{p} \in H $, and $ |\check{\mu} \check{p}|_H \leq \sqrt{2} |\check{\mu}|_H |\check{p}|_V $,
                    \item if $ \hat{\mu} \in L^\infty(0, T; H) $ and $ \hat{p} \in \mathscr{V} $, then $ \hat{\mu} \hat{p} \in \mathscr{H} $, and $ |\hat{\mu} \hat{p}|_\mathscr{H} \leq \sqrt{2} |\hat{\mu}|_{L^\infty(0, T; H)} |\hat{p}|_\mathscr{V} $.
                \vspace{-2ex}
                \end{itemize}
            }
        \end{cases}
    \end{equation}
    Here, we note that the constant $ \sqrt{2} $ corresponds to the constant of embedding $ V \subset C(\overline{\Omega}) $. Moreover, under the setting $ \Omega := (0, 1) $, this $ \sqrt{2} $ can be used as a upper bound of the constants of embeddings $ V \subset L^q(\Omega) $ and $ V_0 \subset L^q(\Omega) $, for all $ 1 \leq q \leq \infty $.
\end{remark}

\noindent
\underline{\textbf{\textit{Notations in convex analysis. (cf. \cite[Chapter II]{MR0348562})}}} 
For a proper, lower semi-con- tinuous (l.s.c.), and convex function $ \Psi : X \to (-\infty, \infty] $ on a Hilbert space $ X $, we denote by $ D(\Psi) $ the effective domain of $ \Psi $. Also, we denote by $\partial \Psi$ the subdifferential of $\Psi$. The subdifferential $ \partial \Psi $ corresponds to a generalized derivative of $ \Psi $, and it is known as a maximal monotone graph in the product space $ X \times X $. The set $ D(\partial \Psi) := \bigl\{ z \in X \ |\ \partial \Psi(z) \neq \emptyset \bigr\} $ is called the domain of $ \partial \Psi $. We often use the notation ``$ [w_{0}, w_{0}^{*}] \in \partial \Psi $ in $ X \times X $\,'', to mean that ``$ w_{0}^{*} \in \partial \Psi(w_{0})$ in $ X $ for $ w_{0} \in D(\partial\Psi) $ '', by identifying the operator $ \partial \Psi $ with its graph in $ X \times X $.
\medskip

For Hilbert spaces $X_1, \cdots, X_N$, with $1<N \in \mathbb{N}$, let us consider a proper, l.s.c., and convex function on the product space $X_1 \times \dots \times X_N$:
\begin{equation*}
\tilde{\Psi}: w = [w_1,\cdots,w_N] \in X_1 \times\cdots\times X_N \mapsto \tilde{\Psi}(w)=\tilde{\Psi}(w_1,\cdots,w_N) \in (-\infty,\infty]. 
\end{equation*}
On this basis, for any $i \in \{1,\dots,N\}$, we denote by $\partial_{w_i} \tilde{\Psi}:X_1 \times \cdots \times X_N \to X_i$ a set-valued operator, which maps any $w=[w_1,\dots,w_i,\dots,w_N] \in X_1 \times \dots \times X_i \times \dots \times X_N$ to a subset $ \partial_{w_i} \tilde{\Psi}(w) \subset  X_i $, prescribed as follows:
\begin{equation*}
\begin{array}{rl}
\partial_{w_i}\tilde{\Psi}(w)&=\partial_{w_i}\tilde{\Psi}(w_1,\cdots,w_i,\cdots,w_N)
\\[2ex]
&:= \left\{\begin{array}{l|l}\tilde{w}^* \in X_i & \begin{array}{ll}\multicolumn{2}{l}{(\tilde{w}^*,\tilde{w}-w_i)_{X_i} \le \tilde{\Psi}(w_1,\cdots,\tilde{w},\cdots,w_N)}
\\[0.25ex] 
& \quad -\tilde{\Psi}(w_1,\cdots,w_i,\cdots,w_N), \mbox{ for any $\tilde{w} \in X_i$}\end{array}
\end{array}\right\}.
\end{array}
\end{equation*}
As is easily checked, 
\begin{equation}\label{prodSubDif}
    \partial \tilde{\Psi}(w) \subset \partial_{w_1} \tilde{\Psi}(w) \times \cdots \times \partial_{w_N} \tilde{\Psi}(w), \mbox{ for any $ w = [w_1, \dots, w_N] \in X_1 \times \cdots \times X_N $.}
\end{equation}
But, it should be noted that the converse inclusion of \eqref{prodSubDif} is not true, in general. 

\begin{remark}[Examples of the subdifferential]\label{exConvex}
    As one of the representatives of the subdifferentials, we exemplify the following set-valued function $ \Sgn^N: \mathbb{R}^N \rightarrow 2^{\mathbb{R}^N} $, with $ N \in \mathbb{N} $, which is defined as:
\begin{align*}
    \xi = [\xi_1, & \dots, \xi_N] \in \mathbb{R}^N \mapsto \Sgn^N(\xi) = \Sgn^N(\xi_1, \dots, \xi_N) 
    \\
    & := \left\{ \begin{array}{ll}
            \multicolumn{2}{l}{
                    \displaystyle \frac{\xi}{|\xi|} = \frac{[\xi_1, \dots, \xi_N]}{\sqrt{\xi_1^2 +\cdots +\xi_N^2}}, ~ } \mbox{if $ \xi \ne 0 $,}
                    \\[3ex]
            \mathbb{D}^N, & \mbox{otherwise,}
        \end{array} \right.
    \end{align*}
where $ \mathbb{D}^N $ denotes the closed unit ball in $ \mathbb{R}^N $ centered at the origin. Indeed, the set-valued function $ \Sgn^N $ coincides with the subdifferential of the Euclidean norm $ |{}\cdot{}| : \xi \in \mathbb{R}^N \mapsto |\xi| = \sqrt{\xi_1^2 + \cdots +\xi_N^2} \in [0, \infty) $, i.e.:
\begin{equation*}
\partial |{}\cdot{}|(\xi) = \Sgn^N(\xi), \mbox{ for any $ \xi \in D(\partial |{}\cdot{}|) = \mathbb{R}^N $,}
\end{equation*}
and furthermore, it is observed that:
\begin{equation*}
    \partial  |{}\cdot{}|(0) = \mathbb{D}^N \begin{array}{c} \subseteq_{\hspace{-1.25ex}\mbox{\tiny$_/$}}  
\end{array} [-1, 1]^N = \partial_{\xi_1}  |{}\cdot{}|(0) \times \cdots \times \partial_{\xi_N}  |{}\cdot{}|(0).
\end{equation*}
\end{remark}
\medskip

Finally, we mention about a notion of functional convergence, known as ``Mosco-convergence''. 
 
\begin{definition}[Mosco-convergence: cf. \cite{MR0298508}]\label{Def.Mosco}
    Let $ X $ be an abstract Hilbert space. Let $ \Psi : X \rightarrow (-\infty, \infty] $ be a proper, l.s.c., and convex function, and let $ \{ \Psi_n \}_{n = 1}^\infty $ be a sequence of proper, l.s.c., and convex functions $ \Psi_n : X \rightarrow (-\infty, \infty] $, $ n = 1, 2, 3, \dots $.  Then, it is said that $ \Psi_n \to \Psi $ on $ X $, in the sense of Mosco, as $ n \to \infty $, iff. the following two conditions are fulfilled:
\begin{description}
    \item[(\hypertarget{M_lb}{M1}) The condition of lower-bound:]$ \displaystyle \varliminf_{n \to \infty} \Psi_n(\check{w}_n) \geq \Psi(\check{w}) $, if $ \check{w} \in X $, $ \{ \check{w}_n  \}_{n = 1}^\infty \subset X $, and $ \check{w}_n \to \check{w} $ weakly in $ X $, as $ n \to \infty $. 
    \item[(\hypertarget{M_opt}{M2}) The condition of optimality:]for any $ \hat{w} \in D(\Psi) $, there exists a sequence \linebreak $ \{ \hat{w}_n \}_{n = 1}^\infty  \subset X $ such that $ \hat{w}_n \to \hat{w} $ in $ X $ and $ \Psi_n(\hat{w}_n) \to \Psi(\hat{w}) $, as $ n \to \infty $.
\end{description}
As well as, if the sequence of convex functions $ \{ \tilde{\Psi}_\varepsilon \}_{\varepsilon \in \Xi} $ is labeled by a continuous argument $\varepsilon \in \Xi$ with a infinite set $\Xi \subset \mathbb{R}$ , then for any $\varepsilon_{0} \in \Xi$, the Mosco-convergence of $\{ \tilde{\Psi}_\varepsilon \}_{\varepsilon \in \Xi}$, as $\varepsilon \to \varepsilon_{0}$, is defined by those of subsequences $ \{ \tilde{\Psi}_{\varepsilon_n} \}_{n = 1}^\infty $, for all sequences $\{ \varepsilon_n \}_{n=1}^{\infty} \subset \Xi$, satisfying $\varepsilon_{n} \to \varepsilon_{0}$ as $n \to \infty$.
\end{definition}

\begin{remark}\label{Rem.MG}
Let $ X $, $ \Psi $, and $ \{ \Psi_n \}_{n = 1}^\infty $ be as in Definition~\ref{Def.Mosco}. Then, the following hold:
\begin{description}
    \item[(\hypertarget{Fact1}{Fact\,1})](cf. \cite[Theorem 3.66]{MR0773850}, \cite[Chapter 2]{Kenmochi81}) Let us assume that
\begin{center}
$ \Psi_n \to \Psi $ on $ X $, in the sense of  Mosco, as $ n \to \infty $,
\vspace{0ex}
\end{center}
and
\begin{equation*}
\left\{ ~ \parbox{9cm}{
$ [w, w^*] \in X \times X $, ~ $ [w_n, w_n^*] \in \partial \Psi_n $ in $ X \times X $, $ n \in \mathbb{N} $,
\\[1ex]
$ w_n \to w $ in $ X $ and $ w_n^* \to w^* $ weakly in $ X $, as $ n \to \infty $.
} \right.
\end{equation*}
Then, it holds that:
\begin{equation*}
[w, w^*] \in \partial \Psi \mbox{ in $ X \times X $, and } \Psi_n(w_n) \to \Psi(w) \mbox{, as $ n \to \infty $.}
\end{equation*}
    \item[(\hypertarget{Fact2}{Fact\,2})](cf. \cite[Lemma 4.1]{MR3661429}, \cite[Appendix]{MR2096945}) Let $ N \in \mathbb{N} $ denote dimension constant, and let $  S \subset \mathbb{R}^N $ be a bounded open set. Then, a sequence $ \{ \widehat{\Psi}_n^S \}_{n = 1}^\infty $ of proper, l.s.c., and convex functions on $ L^2(S; X) $, defined as:
        \begin{equation*}
            w \in L^2(S; X) \mapsto \widehat{\Psi}_n^S(w) := \left\{ \begin{array}{ll}
                    \multicolumn{2}{l}{\displaystyle \int_S \Psi_n(w(t)) \, dt,}
                    \\[1ex]
                    & \mbox{ if $ \Psi_n(w) \in L^1(S) $,}
                    \\[2.5ex]
                    \infty, & \mbox{ otherwise,}
                \end{array} \right. \mbox{for $ n = 1, 2, 3, \dots $;}
        \end{equation*}
        converges to a proper, l.s.c., and convex function $ \widehat{\Psi}^S $ on $ L^2(S; X) $, defined as:
        \begin{equation*}
            z \in L^2(S; X) \mapsto \widehat{\Psi}^S(z) := \left\{ \begin{array}{ll}
                    \multicolumn{2}{l}{\displaystyle \int_S \Psi(z(t)) \, dt, \mbox{ if $ \Psi(z) \in L^1(S) $,}}
                    \\[2ex]
                    \infty, & \mbox{ otherwise;}
                \end{array} \right. 
        \end{equation*}
        on $ L^2(S; X) $, in the sense of Mosco, as $ n \to \infty $. 
\end{description}
\end{remark}
\begin{remark}[Example of Mosco-convergence]\label{Rem.ExMG}
    For any $ \varepsilon \geq 0 $, let $ f_\varepsilon : \mathbb{R} \longrightarrow [0, \infty) $ be a continuous and convex function, defined as:
    \begin{equation}\label{f_eps}
        f_\varepsilon : \xi \in \mathbb{R} \mapsto f_\varepsilon(\xi) := \sqrt{\varepsilon^2 +|\xi|^2} \in [0, \infty).
    \end{equation}
    Then, due to the uniform estimate:
    \begin{equation*}
        \bigl| f_\varepsilon(\xi) -|\xi| \bigr| \leq \varepsilon, \mbox{ for all $ \xi \in \mathbb{R} $,}
    \end{equation*}
    we easily see that:
    \begin{equation*}
        f_\varepsilon \to f_0 ~ (= |{}\cdot{}|) \mbox{ on $ \mathbb{R} $, in the sense of Mosco, as $ \varepsilon \downarrow 0 $.}
    \end{equation*}
    In addition, for any $ \varepsilon > 0 $, it can be said that the subdifferential $ \partial f_\varepsilon $ coincides with the usual differential:
    \begin{equation*}
        f_\varepsilon' : \xi \in \mathbb{R} \mapsto f_\varepsilon'(\xi) = \frac{\xi}{\sqrt{\varepsilon^2 +|\xi|^2}} \in \mathbb{R}.
    \end{equation*}
\end{remark}

\section{Auxiliary Lemmas}

In this section, we recall the previous work \cite{MR3888633}, and set up some auxiliary results. In what follows, we let $ \mathscr{Y} := \mathscr{V} \times \mathscr{V}_0 $, with the dual $ \mathscr{Y}^* := \mathscr{V}^* \times \mathscr{V}_0^* $.
Note that $\mathscr{Y}$ is a Hilbert space which is endowed with a uniform convex topology, based on the inner product for product space, as in the Preliminaries (see the paragraph of Abstruct notations).

Besides, we define:
    \begin{equation*}
        \mathscr{Z} := \bigl( W^{1, 2}(0, T; V^*) \cap \mathscr{V} \bigr) \times \bigl( W^{1, 2}(0, T; V_0^*) \cap \mathscr{V}_0 \bigr),
    \end{equation*}
    as a Banach space, endowed with the norm:
    \begin{align*}
        | [\tilde{p}, & \tilde{z}] |_{\mathscr{Z}} := |[\tilde{p}, \tilde{z}]|_{[C([0, T]; H)]^2} +\bigl( |[\tilde{p}, \tilde{z}]|_{\mathscr{Y}}^2 +|[\partial_t \tilde{p}, \partial_t \tilde{z}]|_{\mathscr{Y}^*}^2 \bigr)^{\frac{1}{2}}, \mbox{ for $ [\tilde{p}, \tilde{z}] \in \mathscr{Z} $.}
    \end{align*}

Based on this, let us consider the following linear system of parabolic initial-boundary value problem, denoted by (P):
\bigskip

(P)
\vspace{0ex}
\begin{equation*}
\left\{\parbox{11cm}{
$\partial_{t}p - \partial_{x}^{2} p + \mu(t, x)p + \lambda(t, x)p + \omega(t, x)\partial_{x} z = h(t, x)$, $(t, x) \in Q,$
\\[1ex]
$\partial_{x}p(t, x) = 0$, $(t, x) \in \Sigma$,
\\[1ex]
$p(0, x) = p_{0}(x)$, $x \in \Omega$;
}\right. 
\end{equation*}
\begin{equation*}
\left\{ \parbox{11cm}{
    $a(t, x)\partial_{t}z + b(t, x)z- \partial_{x} \bigl( A(t, x)\partial_{x}z + \nu^{2} \partial_{x}z + \omega(t, x)p \bigr) $ 
    \\
    $ ~~~~~~= k(t, x)$, $(t, x) \in Q $,
    \\[1ex]
    $ z(t, x) = 0 $, $ (t, x) \in \Sigma $,
    \\[1ex]
    $ z(0, x) = z_{0}(x) $, $ x \in \Omega $.
} \right. 
\end{equation*}
\noindent
This system is studied in \cite{MR3888633} as a key-problem for the G$\hat{\mbox{a}}$teaux differential of the cost $\mathcal{J}_{\varepsilon}$. 
In the context, $[a, b, \mu, \lambda, \omega, A] \in [\mathscr{H}]^{6}$ is a given sextuplet of functions which belongs to a subclass $\mathscr{S} \subset [\mathscr{H}]^{6}$, defined as:
\begin{eqnarray}\label{P01}
 \mathscr{S} := \left\{ [\tilde{a}, \tilde{b}, \tilde{\mu}, \tilde{\lambda}, \tilde{\omega}, \tilde{A}] \in [\mathscr{H}]^{6}\middle|
\begin{split}
    &\bullet  \tilde{a} \in W^{1, \infty}(Q)\ \mbox{and log }\tilde{a} \in L^{\infty}(Q),\\
    &\bullet  [\tilde{b}, \tilde{\lambda}, \tilde{\omega}]\in [L^{\infty}(Q)]^{3},\\
    &\bullet  \tilde{\mu} \in L^{\infty}(0, T; H)\ \mbox{with}\ \tilde{\mu} \geq 0\ \mbox{a.e. in}\ Q,\\ 
    &\bullet  \tilde{A} \in L^{\infty}(Q) \ \mbox{with log}\tilde{A} \in L^{\infty}(Q)
\end{split} 
\right\}.
\end{eqnarray}
Also, $[p_{0}, z_{0}] \in [H]^{2}$  and $[h, k] \in \mathscr{Y}^* $ are, respectively, an initial pair and forcing pair, in the system (P).
\medskip

Now, we refer to the previous work \cite{MR3888633}, to recall the key-properties of the system (P), in forms of Propositions. 

\begin{proposition}[cf. {\cite[Main Theorem 1 (I-A)]{MR3888633}}]\label{Prop(I-A)}
    For any sextuplet $[a, b, \mu, \lambda, \omega, A] \in \mathscr{S}$, any initial pair $[p_{0}, z_{0}] \in [H]^{2}$, and any forcing pair $[h, k] \in \mathscr{Y}^* $, the system (P) admits a unique solution, in the sense that: 
    \begin{equation}\label{ASY1-01}
        \begin{cases}
            p\in W^{1, 2}(0, T; V^{*})\cap L^{2}(0, T; V) \subset C([0, T]; H),
            \\
            z \in W^{1, 2}(0, T; V_{0}^{*})\cap L^{2}(0, T; V_{0}) \subset C([0, T]; H);
        \end{cases}
    \end{equation}
    \begin{align}\label{ASY1-02}
        \displaystyle \langle \partial_{t}p&(t),  \varphi \rangle_{V} + (\partial_{x}p(t), \partial_{x}\varphi)_{H} +( \mu(t)p(t), \varphi )_H  
        \nonumber
        \\
        & + (\lambda(t)p(t) + \omega(t)\partial_{x}z(t), \varphi)_{H} = \left<h(t), \varphi\right>_{V}, 
        \\
        & \mbox{ for any $ \varphi \in V $, a.e. $ t \in (0, T) $, subject to $ p(0) = p_0 $ in $ H $;}
        \nonumber
    \end{align}
    and
    \begin{align}\label{ASY1-03}
        \langle \partial_{t}z & (t), a(t)\psi \rangle_{V_{0}} + (b(t)z(t), \psi)_{H} 
        \nonumber
        \\
        & +\bigl( A(t)\partial_{x}z(t) + \nu^{2}\partial_{x}z(t) +p(t)\omega(t), \partial_{x}\psi \bigr)_H = \langle k(t), \psi \rangle_{V_{0}},
        \\
        & \mbox{for any $ \psi \in V_{0} $, a.e. $ t \in (0, T) $, subject to $ z(0) = z_0 $ in $ H $.}
        \nonumber
        \end{align}
\end{proposition}
\begin{proposition}[cf. {\cite[Main Theorem 1 (I-B)]{MR3888633}}]\label{Prop(I-B)}
    For every $ \ell = 1, 2 $, let us take arbitrary $ [a^\ell, b^\ell, \mu^\ell, \lambda^\ell, \omega^\ell, A^\ell] \in \mathscr{S} $, $ [p_0^\ell, z_0^\ell] \in [H]^2 $, and $ [h^\ell, k^\ell] \in \mathscr{Y}^* $, and let us denote by $ [p^\ell, z^\ell] \in [\mathscr{H}]^2 $ the solution to (P), corresponding to the sextuplet $  [a^\ell, b^\ell, \mu^\ell, \lambda^\ell, \omega^\ell, A^\ell] $, initial pair $ [p_0^\ell, z_0^\ell] $, and forcing pair $ [h^\ell, k^\ell] $. Besides, let $ C_0^* = C_0^*(a^1, b^1, \lambda^1, \omega^1) $ be a positive constant, depending on $a^1, b^1, \lambda^1,$ and $\omega^1$, which is defined as:
    \begin{align}\label{C_0*}
        C_0^* := & \frac{81(1 +\nu^2)}{\min \{ 1, \nu^2, \inf a^1(Q) \}} 
        \bigl( 1 +|a^1|_{W^{1, \infty}(Q)} +|b^1|_{L^\infty(Q)} +|\lambda^1|_{L^\infty(Q)} +|\omega^1|_{L^\infty(Q)}^2 \bigr).
    \end{align}
    Then, it is estimated that:
    \begin{align}\label{est_I-B} \frac{d}{dt} & \bigl( |(p^1 -p^2)(t)|_H^2 +|\sqrt{a^1(t)}(z^1 -z^2)(t)|_H^2 \bigr)
        \nonumber
        \\
        & \quad +\bigl( |(p^1 -p^2)(t)|_V^2 +\nu^2 |(z^1 -z^2)(t)|_{V_0}^2 \bigr)
        \nonumber
        \\
        \leq & 3 C_0^* \bigl( |(p^1 -p^2)(t)|_H^2 +|\sqrt{a^1(t)}(z^1 -z^2)(t)|_H^2 \bigr)
        \\
        & \quad +2C_0^*  \bigl( |(h^1 -h^2)(t)|_{V^*}^2 +|(k^1 -k^2)(t)|_{V_0^*}^2 +R_0^*(t) \bigr),
        \nonumber
        \\
        & \mbox{for a.e. $ t \in (0, T) $;}
        \nonumber
    \end{align}
    where
    \begin{align*}
        & R_0^* (t) := |\partial_t z^2(t)|_{V_0^*}^2 \bigl( |a^1 -a^2|_{C(\overline{Q})}^2 +|\partial_x (a^1 -a^2)(t)|_{L^4(\Omega)}^2 \bigr)
        \\
        & \quad +|p^2(t)|_{V}^2 \bigl( |(\mu^1 -\mu^2)(t)|_H^2 +|(\omega^1 -\omega^2)(t)|_{L^4(\Omega)}^2 \bigr)
        \\
        & \quad +|z^2(t)|_{V_0}^2 \bigl( |(b^1 -b^2)(t)|_{L^4(\Omega)}^2 +|p^2(t)(\lambda^1 -\lambda^2)(t)|_{H}^2 \bigr)
        \\
        & \quad +|\partial_x z^2(t) (\omega^1 -\omega^2)(t)|_{H}^2 +|(A^1 -A^2)(t) \partial_x z^2(t)|_{H}^2,
            \\
        & \mbox{for a.e. $ t \in (0, T) $.}
    \end{align*}
\end{proposition}
\begin{remark}\label{Rem.C_0*}
    In the previous work \cite{MR3888633}, the constant $ C_0^* $ for the estimate \eqref{est_I-B} is provided as:
    \begin{align}\label{C_0*org}
        C_0^* := & \frac{9(1 +\nu^2)}{\min \{ 1, \nu^2, \inf a^1(Q) \}} \cdot \bigl( 1 +(C_{V}^{L^4})^2 +(C_{V}^{L^4})^4 +(C_{V_0}^{L^4})^2 \bigr)
        \nonumber
        \\
        & \quad \cdot 
        \bigl( 1 +|a^1|_{W^{1, \infty}(Q)} +|b^1|_{L^\infty(Q)} +|\lambda^1|_{L^\infty(Q)} +|\omega^1|_{L^\infty(Q)}^2 \bigr),
    \end{align}
    with use of the constants $ C_{V}^{L^4} > 0 $ and $ C_{V_0}^{L^4} > 0  $ of the respective embeddings $ V \subset L^4(\Omega) $ and $ V_0 \subset L^4(\Omega) $. Note that the setting \eqref{C_0*} corresponds to the special case of the original one \eqref{C_0*org}, under the one-dimensional situation, as in Remark \ref{Rem.Prelim01}. 
\end{remark}
\begin{proposition}[cf. {\cite[Corollary 1]{MR3888633}}]\label{ASY_Cor.1}
    For any $ [a, b, \mu, \lambda, \omega, A] \in \mathscr{S} $, let us denote by $\mathcal{P} = \mathcal{P}(a, b, \mu, \lambda, \omega, A) : [H]^2 \times \mathscr{Y}^* \longrightarrow \mathscr{Z} $ a linear operator, which maps any $ \bigl[ [p_0, z_0], [h, k] \bigr] \in [H]^2 \times \mathscr{Y}^* $ to the solution $ [p, z] \in \mathscr{Z} $ to the linear system (P), for the sextuplet $ [a, b, \mu, \lambda, \omega, A] $, initial pair $ [p_0, z_0] $, and forcing pair $ [h, k] $. Then, for any sextuplet $ [a, b, \mu, \lambda, \omega, A] \in \mathscr{S} $, there exist positive constants $ M_0^* = M_0^*(a, b, \mu, \lambda, \omega, A) $ and $ M_1^* = M_1^*(a, b, \mu, \lambda, \omega, A) $, depending on $ a $, $ b $, $ \mu $, $ \lambda$, $\omega $, and $A$, such that:
    \begin{equation*}
        \begin{array}{c}
            M_0^* \bigl| \bigl[ [p_0, z_0], [h, k] \bigr] \bigr|_{[H]^2 \times \mathscr{Y}^*} \leq | [p, z] |_{\mathscr{Z}} \leq M_1^* \bigl| \bigl[ [p_0, z_0], [h, k] \bigr] \bigr|_{[H]^2 \times \mathscr{Y}^*},
            \\[1ex]
            \mbox{for all $ [p_0, z_0] \in [H]^2 $, $ [h, k] \in \mathscr{Y}^* $,}
            \\
            \mbox{and $ [p, z] = \mathcal{P}(a, b, \mu, \lambda, \omega, A)\bigl[ [p_0, z_0], [h, k] \bigr] \in \mathscr{Z} $,}
        \end{array}
    \end{equation*}
    i.e. the operator $ \mathcal{P} = \mathcal{P}(a, b, \mu, \lambda, \omega, A) $ is an isomorphism between the Hilbert space $ [H]^2 \times \mathscr{Y}^* $ and the Banach space $ \mathscr{Z} $. 
\end{proposition}
\begin{proposition}[cf. {\cite[Corollary 2]{MR3888633}}]\label{ASY_Cor.2}
    Let us assume: 
    \begin{equation*}
        [a, b, \mu, \lambda, \omega, A] \in \mathscr{S}, ~ \{ [a^n, b^n, \mu^n, \lambda^n, \omega^n, A^n] \}_{n = 1}^\infty \subset \mathscr{S}, 
    \end{equation*}
    \begin{align}\label{ASY01}
        [a^n, & \partial_t a^n, \partial_x a^n, b^n, \lambda^n, \omega^n, A^n] \to [a, \partial_t a, \partial_x a, b, \lambda, \omega, A] 
        \nonumber
        \\
        & \mbox{weakly-$*$ in $ [L^\infty(Q)]^7 $, and in the pointwise sense a.e. in  $ Q $,} 
        \\
        & \mbox{as}\ n \to \infty,
        \nonumber
    \end{align}
    and
    \begin{equation*}
        \begin{cases}
            \mu^n \to \mu \mbox{ weakly-$*$ in $ L^\infty(0, T; H) $,}
            \\
            \mu^n \to \mu \mbox{ in $ H $, for a.e. $ t \in (0, T) $,}
        \end{cases}
        \mbox{as $ n \to \infty $.}
    \end{equation*}
    Let us assume $ [p_0, z_0] \in [H]^2 $, $ [h, k] \in \mathscr{Y}^* $, and let us denote by $ [p, z] \in [\mathscr{H}]^2 $ the solution to (P), for the initial pair $ [p_0, z_0] $ and forcing pair $ [h, k] $. Also, let us assume $ \{ [p_0^n, z_0^n] \}_{n = 1}^\infty \subset [H]^2 $, $ \{ [h^n, k^n] \}_{n = 1}^\infty \subset \mathscr{Y}^* $, and for any $ n \in \mathbb{N} $, let us denote by $ [p^n, z^n] \in [\mathscr{H}]^2 $ the solution to (P), for the initial pair $ [p_0^n, z_0^n] $ and forcing pair $ [h^n, k^n] $. Then, the following two items hold.
    \begin{itemize}
        \item[\textmd{\em(A)}]The convergence:
            \begin{equation*}
                \begin{cases}
                    [p_0^n, z_0^n] \to [p_0, z_0] \mbox{ in $ [H]^2 $},
                    \\
                    [h^n, k^n] \to [h, k] \mbox{ in $ \mathscr{Y}^* $,}
                \end{cases}
                \mbox{as $ n \to \infty $,}
            \end{equation*}
            implies the convergence:
            \begin{align*}
                [p^n, z^n] & \to [p, z] \mbox{ in $ [C([0, T]; H)]^2 $, and in $ \mathscr{Y} $, as $ n \to \infty $.}
            \end{align*}
        \item[\textmd{\em(B)}]The following two convergences:
            \begin{equation*}
                \begin{cases}
                    [p_0^n, z_0^n] \to [p_0, z_0] \mbox{ weakly in $ [H]^2 $},
                    \\
                    [h^n, k^n] \to [h, k] \mbox{ weakly in $ \mathscr{Y}^* $,}
                \end{cases}
                \mbox{as $ n \to \infty $,}
            \end{equation*}
            and
            \begin{align*}
                [p^n, z^n] & \to [p, z] \mbox{ in $ [\mathscr{H}]^2 $, weakly in $ \mathscr{Y} $,}
                \\
                & \mbox{and weakly in $ W^{1, 2}(0, T; V^*) \times W^{1, 2}(0, T; V_0^*) $, as $ n \to \infty $,}
            \end{align*}
        are equivalent each other.
    \end{itemize}
\end{proposition}

\section{Main Theorems}

We begin by setting up some assumptions needed in our Main Theorems. 
\begin{itemize}
\item[\textmd{(A1)}]\label{given-f-c}
    $\nu>0$ is a fixed constant. Let $ [\eta_0, \theta_0] \in V \times V_0 $ be a fixed initial pair. Let $ [\eta_\mathrm{ad}, \theta_\mathrm{ad}] \in [\mathscr{H}]^2 $ be a fixed pair of functions, called the \emph{admissible target profile}. 
\item[\textmd{(A2)}]\label{g}
    $g : \mathbb{R} \longrightarrow  \mathbb{R}$ is a $C^{1}$-function, which is a Lipschitz continuous on $\mathbb{R}$. Also $g$ has a nonnegative primitive $ 0 \leq G \in C^{2}(\mathbb{R})$, i.e. the derivative $ G'= \frac{dG}{d\eta} $ coincides with $ g $ on $ \mathbb{R} $.
\item[\textmd{(A3)}]\label{alpha}
    $\alpha:\mathbb{R} \longrightarrow (0, \infty)$ and $\alpha_{0}: Q \longrightarrow (0, \infty)$ are Lipschitz continuous functions, such that:
        \begin{itemize}
            \item $\alpha \in C^2(\mathbb{R})$, with the first derivative $ \alpha' = \frac{d \alpha}{d \eta} $ and the second one $ \alpha'' = \frac{d^2 \alpha}{d\eta^2} $;
            \item $\alpha'(0) = 0$, $ \alpha'' \geq 0 $ on $ \mathbb{R} $, and $\alpha \alpha'$ is a Lipschitz continuous function on $\mathbb{R}$; 
            \item $ \alpha \geq \delta_* $ on $ \mathbb{R} $, and $ \alpha_0 \geq \delta_* $ on $ \overline{Q} $, for some constant $ \delta_* \in (0, 1) $.
        \end{itemize}
\end{itemize}
Additionally, for any $ \varepsilon \geq 0 $, let $ f_\varepsilon : \mathbb{R} \longrightarrow [0, \infty) $ be the convex function, defined in \eqref{f_eps}.
\bigskip

Now, the Main Theorems of this paper are stated as follows:
\begin{main}\label{mainTh01}
    Let us assume (A1)--(A3). Let us fix a constant $ \varepsilon \geq 0 $, an initial pair $ [\eta_0, \theta_0] \in V \times V_0 $, and a forcing pair $ [u, v] \in [\mathscr{H}]^2 $. Then, the following hold:
    \begin{itemize}
    \item[\textmd{(I-A)}]The state-system (S)$_{\varepsilon}$ admits a unique solution $[\eta, \theta] \in [\mathscr{H}]^{2}$, in the sense that:
        \begin{equation}\label{S1}
            \begin{cases}
            \eta \in W^{1, 2}(0, T; H) \cap L^{\infty}(0, T; V) \subset C(\overline{Q}),
            \\
            \theta \in W^{1, 2}(0, T; H) \cap L^{\infty}(0, T; V_{0}) \subset C(\overline{Q});
        \end{cases}
        \end{equation}
        \begin{equation}\label{S2}
                \begin{array}{c}
                    \displaystyle\bigl( \partial_{t}\eta(t), \varphi \bigr)_{H} + \bigl( \partial_{x}\eta(t), \partial_{x}\varphi \bigr)_{H} +\bigl( g(\eta(t)), \varphi \bigr)_{H} 
                    \\[1ex]
                    +\bigl( \alpha'(\eta(t))f_{\varepsilon}(\partial_{x}\theta(t)), \varphi \bigr)_{H}  = \bigl( M_u u(t), \varphi \bigr)_{H}, 
                    \\[1ex]
                    \mbox{for any $ \varphi \in V $, a.e. $ t \in (0, T) $, subject to } \eta(0) = \eta_0 \mbox{ in $ H $;}
                \end{array}
        \end{equation}
and
        \begin{equation}\label{S3}
                \begin{array}{c}
                    \displaystyle\bigl( \alpha_{0}(t)\partial_{t}\theta(t), \theta(t)-\psi \bigr)_{H} + \nu^2 \bigl( \partial_{x}\theta(t), \partial_{x}(\theta(t)-\psi) \bigr)_{H} 
                    \\[1ex]
                    \displaystyle +\int_{\Omega}\alpha(\eta(t))f_{\varepsilon}(\partial_{x}\theta(t))dx \leq \int_{\Omega}\alpha(\eta(t))f_{\varepsilon}(\partial_{x}\psi)dx
                    \\[2ex]
                     +\bigl( M_v v(t), \theta(t)-\psi \bigr)_{H}, \mbox{ for any } \psi \in V_{0},
                     \\[1ex]
                     \mbox{a.e. } t \in (0, T), \mbox{ subject to $ \theta(0) = \theta_0 $ in $ H $.} 
                \end{array}
        \end{equation}
        \item[\textmd{(I-B)}]Let $\{\varepsilon_{n} \}_{n=1}^{\infty} \subset [0, 1]$, $\{[\eta_{0, n}, \theta_{0, n}] \}_{n=1}^{\infty} \subset V \times V_{0}$, and $\{ [u_{n}, v_{n}]\}_{n=1}^{\infty} \subset [\mathscr{H}]^{2}$ be given sequences such that:
\begin{equation}\label{w.i}
    \begin{array}{c}
        \varepsilon_{n} \to \varepsilon,\ [\eta_{0, n}, \theta_{0, n}] \to [\eta_{0}, \theta_{0}] \mbox{ weakly\ in } V \times V_{0},
        \\[1ex]
        \mbox{and } [M_u u_{n}, M_v v_{n}] \to [M_u u, M_v v] \ \mbox{weakly\ in}\ [\mathscr{H}]^{2}, \mbox{ as } n \to \infty.
    \end{array}
\end{equation}
In addition, let $[\eta, \theta]$ be the unique solution to (S)$_{\varepsilon}$, for the forcing pair $[u, v]$, and for any $n \in \mathbb{N}$, let $[\eta_{n}, \theta_{n}]$ be the unique solution to (S)$_{\varepsilon_n}$, for the initial pair $ [\eta_{0, n}, \theta_{0, n}] $ and forcing pair $ [u_{n}, v_{n}] $. Then, it holds that:
\begin{align}\label{mThConv}
    [\eta_{n}, & \theta_{n}] \to [\eta, \theta] \mbox{ in } [C(\overline{Q})]^2, \mbox{ in } \mathscr{Y}, \mbox{ weakly in $ [W^{1, 2}(0, T; H)]^2 $,}
    \\
    & \mbox{and weakly-$*$ in $ L^\infty(0, T; V) \times L^\infty(0, T; V_0) $, as $ n \to \infty $,}
    \nonumber
\end{align}
and in particular,
        \begin{align}\label{mThConv00}
            \alpha''(\eta_{n}) & f_{\varepsilon_{n}}(\partial_x \theta_{n}) \to \alpha''(\eta)f_\varepsilon (\partial_x \theta) \mbox{ in $ \mathscr{H} $,}
            \nonumber
            \\
            & \mbox{and weakly-$*$ in $ L^\infty(0, T; H) $, as $ n \to \infty $.}
        \end{align}
    \end{itemize}
\end{main}
\begin{remark}\label{Rem.mTh01Conv}
    As a consequence of \eqref{mThConv} and \eqref{mThConv00}, we further find a subsequence $ \{ n_i \}_{i = 1}^\infty \subset \{n\} $, such that:
    \begin{equation}\label{mThConv01}
        \begin{array}{l}
            [\eta_{n_i}, \theta_{n_i}] \to [\eta, \theta], ~ [\partial_x \eta_{n_i}, \partial_x \theta_{n_i}] \to [\partial_x \eta, \partial_x \theta], 
            \\[0.5ex]
            \qquad \qquad \mbox{and } \alpha''(\eta_{n_i}) f_{\varepsilon_{n_i}}(\partial_x \theta_{n_i}) \to \alpha''(\eta) f_\varepsilon(\partial_x \theta), 
            \\[0.5ex]
            \qquad \qquad \qquad \qquad \mbox{in the pointwise sense a.e. in $ Q $, as $ i \to \infty $,}
        \end{array}
    \end{equation}
    and
    \begin{equation}\label{mThConv02}
        \begin{array}{l}
            [\eta_{n_i}(t), \theta_{n_i}(t)] \to [\eta(t), \theta(t)] \mbox{ in $ V \times V_0 $,}
            \\[0.5ex]
            \qquad \qquad \mbox{and } \alpha''(\eta_{n_i}(t)) f_{\varepsilon_{n_i}}(\partial_x \theta_{n_i}(t)) \to \alpha''(\eta(t)) f_\varepsilon(\partial_x \theta(t)) \mbox{ in $ H $,}
            \\[0.5ex]
            \qquad \qquad \qquad \qquad \mbox{in the pointwise sense for a.e. $ t \in (0, T) $, as $ i \to \infty $.}
        \end{array}
    \end{equation}
\end{remark}
\begin{main}\label{mainTh02}
    Let us assume (A1)--(A3), and fix any constant $ \varepsilon \geq 0 $. Then, the following two items hold.
    \begin{itemize}
        \item[\textmd{(II-A)}]The problem (OP)$_{\varepsilon}$ has at least one optimal control $[u^{*}, v^{*}] \in [\mathscr{H}]^{2}$, so that:
\begin{equation*}
\mathcal{J}_{\varepsilon}(u^{*}, v^{*}) = \min_{[u, v] \in [\mathscr{H}]^{2}}\mathcal{J}_{\varepsilon}(u, v).
\end{equation*}
        \item[\textmd{(II-B)}]Let $\{\varepsilon_{n} \}_{n=1}^{\infty} \subset [0, 1]$ and $\{[\eta_{0, n}, \theta_{0, n}] \}_{n=1}^{\infty} \subset V\times V_{0}$ be given sequences such that: 
\begin{equation}\label{ass.4}
\varepsilon_{n} \to \varepsilon, \mbox{and}\ [\eta_{0, n}, \theta_{0, n}] \to [\eta_{0}, \theta_{0}]\ {\rm weakly\ in}\ V\times V_{0},\ {\rm as}\ n \to \infty.
\end{equation}
In addition, for any $n \in \mathbb{N}$, let $[u_{n}^*, v_{n}^*] \in [\mathscr{H}]^{2}$ be the optimal control of (OP)$_{\varepsilon_{n}}$. Then, there exist a subsequence $ \{ n_{i} \}_{i = 1}^{\infty} \subset \{ n \} $ and a pair of functions $[u^{**}, v^{**}] \in [\mathscr{H}]^{2}$, such that: 
\begin{align*}
    \varepsilon_{n_{i}} \to \varepsilon,  \mbox{ and } & [M_u u_{n_{i}}^*, M_v v_{n_{i}}^*] \to [M_u u^{**}, M_v v^{**}] 
    \\
    & \mbox{weakly in } [\mathscr{H}]^{2}, \mbox{ as } i \to \infty,
\end{align*}
and 
\begin{equation*}
[u^{**}, v^{**}] \ \mbox{is an optimal control of (OP)}_{\varepsilon}.
\end{equation*}
\end{itemize}
\end{main}

\begin{main}
    \label{mainTh03}
    Under the assumptions (A1)--(A3), the following two items hold.
    \begin{itemize}
        \item[\textmd{\textit{(III-A)}}](Necessary condition for (OP)$_\varepsilon$ when $ \varepsilon > 0 $) 
            For any $ \varepsilon > 0 $, let $[u_{\varepsilon}^{*}, v_{\varepsilon}^{*}] \in [\mathscr{H}]^{2}$ be an optimal control of (OP)$_{\varepsilon}$, and let $ [\eta_{\varepsilon}^{*}, \theta_{\varepsilon}^{*}] $ be the solution to (S)$_{\varepsilon}$, for the initial pair $[\eta_0, \theta_0]$ and forcing pair $ [u_{\varepsilon}^{*}, v_{\varepsilon}^{*}] \in [\mathscr{H}]^{2} $. Then, it holds that:
    \begin{equation}\label{Thm.5-00}
        \displaystyle [M_u (u_\varepsilon^{*} + p_{\varepsilon}^{*}), M_v (v_\varepsilon^{*} + z_{\varepsilon}^{*})] = [0, 0] 
    \mbox{ in $ [\mathscr{H}]^2 $,}
\end{equation}
where $[p_{\varepsilon}^{*}, z_{\varepsilon}^{*}]$ is a unique solution to the following variational system:
\begin{align}\label{Thm.5-01}
    -\bigl< \partial_{t} p_{\varepsilon}^{*}(t), & \varphi \bigr>_{V} + \bigl( \partial_{x} p_{\varepsilon}^{*}(t), \partial_{x}\varphi \bigr)_{H} + \bigl( [\alpha''(\eta_{\varepsilon}^{*})f_{\varepsilon}(\partial_{x}\theta_{\varepsilon}^{*})](t) p_{\varepsilon}^{*}(t), \varphi \bigr)_{H}
    \nonumber
    \\
    & +\bigl( g'(\eta_{\varepsilon}^{*}(t)) p_{\varepsilon}^{*}(t), \varphi \bigr)_{H} +\bigl( [\alpha'(\eta_{\varepsilon}^{*}) f_{\varepsilon}'(\partial_{x}\theta_{\varepsilon}^{*})](t)\partial_{x}z_{\varepsilon}^{*}(t), \varphi \bigr)_{H}
    \\
    = & \bigl( M_\eta (\eta_{\varepsilon}^{*}-\eta_{\mbox{\scriptsize ad}})(t), \varphi \bigr)_{H}, \mbox{ for any $ \varphi \in V $,  and a.e. $ t \in (0, T) $;}
    \nonumber
\end{align}
and
\begin{align}\label{Thm.5-02}
    -\bigl< \partial_{t} & \bigl( \alpha_{0}z_{\varepsilon}^{*} \bigr)(t), \psi \bigr>_{V_{0}}  +\bigl( [\alpha(\eta_{\varepsilon}^{*})f_{\varepsilon}''(\partial_{x}\theta_{\varepsilon}^{*})](t)\partial_{x}z_{\varepsilon}^{*}(t) + \nu^2 \partial_{x}z_{\varepsilon}^{*}(t), \partial_x \psi \bigl)_H
    \nonumber
    \\
    & +\bigl( [\alpha'(\eta_{\varepsilon}^{*})f_{\varepsilon}'(\partial_{x}\theta_{\varepsilon}^{*})](t)p_{\varepsilon}^{*}(t), \partial_{x}\psi \bigr)_{H} = \bigl( M_\theta (\theta_{\varepsilon}^{*}-\theta_{\mbox{\scriptsize ad}})(t), \psi \bigr)_{H},
    \\
    & \mbox{for any}\ \psi \in V_{0},\ \mbox{and a.e.}\ t \in (0, T);
    \nonumber
\end{align}
subject to the terminal condition:
\begin{equation}\label{Thm.5-03}
[p_{\varepsilon}^{*}(T), z_{\varepsilon}^{*}(T)] = [0, 0]\ \mbox{in}\ [H]^{2}.
\end{equation}
    \item[\textmd{\textit{(III-B)}}]Let us define a Hilbert space $ \mathscr{W}_0 $ as:
        \begin{equation*}
            \mathscr{W}_0 := \left\{ \begin{array}{l|l}
                \psi \in W^{1, 2}(0, T; H) \cap \mathscr{V}_0 & \psi(0) = 0 \mbox{ in $ H $}
            \end{array} \right\}.
        \end{equation*}
        Then, there exists an optimal control $ [u^{\circ}, v^{\circ}] \in [\mathscr{H}]^2 $ of the problem (OP)$_0$, together with the solution $ [\eta^\circ, \theta^\circ] $ to the system (S)$_0$, for the initial pair $[\eta_0, \theta_0]$ and forcing pair $ [u^\circ, v^\circ] $, and there exist pairs of functions $ [p^{\circ}, z^{\circ}] \in \mathscr{Y} $, $ [\xi^\circ, \nu^\circ] \in \mathscr{H} \times L^\infty(Q) $, and a distribution $ \zeta^\circ \in \mathscr{W}_0^* $, such that:
    \begin{equation}\label{Thm.5-10}
        [M_u (u^\circ + p^\circ), M_v (v^\circ + z^\circ)] = [0, 0] \mbox{ in $ [\mathscr{H}]^2 $;}
    \end{equation}
    \begin{equation}\label{Thm.5-13}
        \begin{cases}
            p^\circ \in W^{1, 2}(0, T; V^*) ~ (\cap \mathscr{V}), 
            \mbox{ i.e. $ p^\circ \in C([0, T]; H) $,} 
            \\
            \nu^\circ \in \Sgn^1(\partial_x \theta^\circ), \mbox{ a.e. in $ Q $;}
        \end{cases}
    \end{equation}
    \begin{align}\label{Thm.5-11}
        \bigl< -\partial_t & p^\circ, \varphi \bigr>_{\mathscr{V}} +\bigl( \partial_x p^\circ, \partial_x \varphi \bigr)_{\mathscr{H}}  +\bigl( \alpha''(\eta^\circ)|\partial_x \theta^\circ| p^\circ, \varphi \bigl)_\mathscr{H}
        \nonumber
        \\
        & +\bigl( g'(\eta^\circ)p^\circ +\alpha'(\eta^\circ) 
        \xi^\circ, \varphi \bigr)_{\mathscr{H}} = \bigl( M_\eta  (\eta^\circ -\eta_\mathit{ad}), \varphi \bigr)_{\mathscr{H}}, 
        \\
        & 
        \mbox{for any $ \varphi \in \mathscr{V} $, subject to $ p^\circ(T) = 0 $ in $ H $;}
        \nonumber
    \end{align}
    and
    \begin{align}\label{Thm.5-12}
        \bigl( \alpha_0 z^\circ & , \partial_t \psi \bigr)_{\mathscr{H}} +\bigl< \zeta^\circ, \psi \bigr>_{\mathscr{W}_0} +\bigl( \nu^2 \partial_x z^\circ +\alpha'(\eta^\circ) \nu^\circ p^\circ, \partial_x \psi \bigr)_{\mathscr{H}}
        \nonumber
        \\
        & = \bigl( M_\theta (\theta^\circ -\theta_\mathit{ad}), \psi \bigr)_{\mathscr{H}}, \mbox{ for any $ \psi \in \mathscr{W}_0 $.}
    \end{align}
\end{itemize}
\end{main}

\begin{remark}\label{Rem.mTh03}
    Let $ \mathcal{R}_T \in \mathcal{L}(\mathscr{H}) $ be an isomorphism, defined as:
    \begin{equation*}
        \bigl( \mathcal{R}_T \varphi \bigr)(t) := \varphi(T -t) \mbox{ in $ H $, for a.e. $ t \in (0, T) $.}
    \end{equation*}
    Also, let us fix $ \varepsilon > 0 $, and define a bounded linear operator $ \mathcal{Q}_\varepsilon^* : [\mathscr{H}]^2 \longrightarrow \mathscr{Z} $ as the restriction $\mathcal{P}|_{\{[0, 0]\} \times \mathscr{Y}^{*}} $ of the linear isomorphism $ \mathcal{P} = \mathcal{P}(a, b, \mu, \lambda, \omega, A) : [H]^2 \times \mathscr{Y}^* \longrightarrow \mathscr{Z} $, as in Proposition \ref{ASY_Cor.1}, in the case when:
    \begin{equation}\label{setRem4}
        \begin{cases}
            [a, b] = \mathcal{R}_T [\alpha_0, -\partial_t \alpha_0] \mbox{ in $ W^{1, \infty}(Q) \times L^\infty(Q) $,}
            \\
            \mu = \mathcal{R}_T \bigl[ \alpha''(\eta_\varepsilon^*) f_\varepsilon(\partial_x \theta_\varepsilon^*) \bigr] \mbox{ in $ L^\infty(0, T; H) $,}
            \\
            [\lambda, \omega, A] = \mathcal{R}_T \bigl[ g'(\eta_\varepsilon^*), \alpha'(\eta_\varepsilon^*) f_\varepsilon'(\partial_x \theta_\varepsilon^*), \alpha(\eta_\varepsilon^*)f_\varepsilon''(\partial_x \theta_\varepsilon^*) \bigl] \mbox{ in $ [L^\infty(Q)]^3 $.}
        \end{cases}
    \end{equation}
    On this basis, let us define:
    \begin{equation*}
        \mathcal{P}_\varepsilon^* := \mathcal{R}_T \circ \mathcal{Q}_\varepsilon^* \circ \mathcal{R}_T \mbox{ in $ \mathscr{L}([\mathscr{H}]^2; \mathscr{Z}) $.}
    \end{equation*}
    Then, having in mind:
    \begin{equation}\label{productrule}
        \partial_t (\alpha_0 \tilde{z}) = \alpha_0 \partial_t \tilde{z} +\tilde{z} \partial_t \alpha_0 \mbox{ in $ \mathscr{V}_0^* $, for any $ \tilde{z} \in W^{1, 2}(0, T; V^*_0) $,}
    \end{equation}
    we can obtain the unique solution $[p_{\varepsilon}^{*}, z_{\varepsilon}^{*}] \in [\mathscr{H}]^2$ to the variational system \eqref{Thm.5-01}--\eqref{Thm.5-03} as follows:
\begin{equation*}
    [p_\varepsilon^*, z_\varepsilon^*] = \mathcal{P}_\varepsilon^* \bigl[ M_\eta  (\eta_\varepsilon^* -\eta_\mathrm{ad}), M_\theta  (\theta_\varepsilon^* -\theta_\mathrm{ad}) \bigr] \mbox{ in $ \mathscr{Z} $.}
\end{equation*}
\end{remark}

\section{Proof of Main Theorem \ref{mainTh01}} 

In this section, we give the proof of the first Main Theorem \ref{mainTh01}. Before the proof, we refer to the reformulation method as in \cite{MR3888636}, and consider to reduce the state-system (S)$_{\varepsilon}$ to an evolution equation in the Hilbert space $[H]^{2}$.

Let us fix any $\varepsilon \geq 0$. Besides, let us define  time-dependent operators $\mathcal{A}(t) \in \mathscr{L}([H]^{2})$, for $t \in [0, T]$, a nonlinear operator $\mathcal{G}:[H]^{2} \longrightarrow \mathbb [H]^{2}$, and a proper functional $\Phi_{\varepsilon}:[H]^{2} \longrightarrow [0, \infty]$, by setting:
\begin{equation}\label{AA_0}
    \mathcal{A}(t): {w} = [{\eta}, {\theta}] \in [H]^{2} \mapsto \mathcal{A}(t) {w} := [{\eta}, \, \alpha_0(t) {\theta}] \in [H]^{2}, \mbox{ for }t \in [0, T],
\end{equation}
\begin{equation}\label{GG}
    \mathcal{G}: w = [\eta, \theta] \in [H]^{2} \mapsto \mathcal{G}(w) := \bigl[ g(\eta) -\eta -\nu^{-2}\alpha(\eta)\alpha'(\eta), ~ 0  \bigr] \in [H]^{2},
\end{equation}
and
\begin{align}\label{Phi_eps}
    \Phi_{\varepsilon}: w & = [\eta, \theta] \in [H]^{2} \mapsto \Phi_{\varepsilon}(w) =  \Phi_{\varepsilon}(\eta, \theta)
    \nonumber
    \\
    &:= \left\{
        \begin{array}{ll}
            \multicolumn{2}{l}{\displaystyle\frac{1}{2}\int_{\Omega}|{\partial_x}\eta|^{2}dx +\frac{1}{2} \int_\Omega |\eta|^2 \, dx +\frac{1}{2}\int_{\Omega}\left(\nu f_{\varepsilon}({\partial_x}\theta) + \frac{1}{\nu}\alpha(\eta) \right)^{2}dx,}
            \\[2ex]
            & \mbox{if $[\eta, \theta] \in V\times V_{0}$,}
            \\[2ex]
            \infty, & \mbox{otherwise,}
        \end{array}
    \right. 
\end{align}
respectively. Note that the definition of $f_{\varepsilon}$, as in Remark \ref{Rem.ExMG}, and the assumption (A3) guarantee the lower semi-continuity and convexity of $\Phi_{\varepsilon}$ on $[H]^{2}$. 
\begin{remark}\label{Rem.eps>0}
    When $ \varepsilon > 0 $, we can easily check from Remark \ref{Rem.ExMG} and (A3) that the subdifferential $ \partial \Phi_\varepsilon \subset [H]^2 \times [H]^2 $ is single-valued, and
\begin{equation*}
    [w, w^*] \in \partial \Phi_\varepsilon \mbox{ in } [H]^2 \times [H]^2 \mbox{ for } w = [\eta, \theta] \in [H]^2 \mbox{ and } w^* = [\eta^*, \theta^*] \in [H]^2, 
\end{equation*}
iff.
\begin{equation*}
    \begin{cases}
        \hspace{-2ex}
        \parbox{14cm}{
            \vspace{-2ex}
            \begin{itemize}
                \item $ w = [\eta, \theta] \in H^2(\Omega) \times V_0 $ with $ {\partial_x} \eta(\ell) = 0 $, for $ \ell \in \Gamma = \{0, 1\} $, and $ \alpha(\eta) f_{\varepsilon}'(\partial_x \theta) +\nu^2 {\partial_x} \theta \in V_0 $,
                \item $ w^* = \rule{0pt}{18pt}^\mathrm{t} \hspace{-0.5ex} \begin{bmatrix}
                        \,\eta^* \\ \,\theta^*
                \end{bmatrix} = \rule{0pt}{18pt}^\mathrm{t} \hspace{-0.5ex} \begin{bmatrix}
                    -{\partial_x^2}\eta +\eta + \alpha'(\eta)f_{\varepsilon}({\partial_x}\theta) + \nu^{-2}\alpha(\eta)\alpha'(\eta)  
                    \\
                    -\partial_x \bigl( \alpha(\eta) f_{\varepsilon}'({\partial_x}\theta) + \nu^{2}{\partial_x}\theta \bigr)
                \end{bmatrix} $ in $ [H]^2 $.
                \vspace{-2ex}
            \end{itemize}
        }
    \end{cases}
\end{equation*}
Therefore, in the case of $ \varepsilon > 0 $, the state-system (S)$_\varepsilon$ will be equivalent to the following Cauchy problem (E)$_\varepsilon$ of an evolution equation: 
\begin{equation*}
    \mathrm{(E)}_\varepsilon~~\begin{cases}
        \mathcal{A}(t) w'(t) +\partial \Phi_\varepsilon(w(t)) +\mathcal{G}(w(t)) \ni \mathfrak{f}(t) \mbox{ in $ [H]^2 $, $ t \in (0, T) $,}
        \\
        w(0) = w_0 \mbox{ in $[H]^2$.}
    \end{cases}
\end{equation*}
In the context, ``\,$'$\,'' denotes the time-derivative, $ w_0 := [\eta_0, \theta_0] \in V \times V_0 $ and $ \mathfrak{f} := [M_{u} u, M_{v} v] \in [\mathscr{H}]^2 $ are the initial pair and forcing pair, as in (S)$_\varepsilon$, respectively. 
\end{remark}
\begin{remark}\label{Rem.Sec03-01}
    In the case of $ \varepsilon = 0 $, the equivalence between the corresponding state-system (S)$_0$ and Cauchy problem (E)$_0$ is not so obvious. However, we can show a partial relation, such that:
    \begin{enumerate}
        \item[$(\star \, 0)$]if $ w = [\eta, \theta] $ is a solution to (E)$_0$, then it is also a solution to (S)$_0$. 
    \end{enumerate}
    In fact, as is easily seen, the operator $ \partial_\eta \Phi_0 : [H]^2 \longrightarrow H $ is single-valued. Besides, for any $ \tilde{\theta} \in V_0 $, it follows that $ [\eta, \eta^*] \in \partial_\eta \Phi_0({}\cdot, \tilde{\theta}) $ in $ H \times H $, iff.:
    \begin{align*}\label{var_eta}
        \bigl( \eta^*, \varphi \bigr)_H &= \bigl( {\partial_x} \eta, {\partial_x} \varphi \bigr)_H +\bigl(\eta, \varphi \bigr)_H \nonumber
        \\
         &\qquad + \bigl( \alpha'(\eta)|{\partial_x} \tilde{\theta}| +\nu^{-2} \alpha(\eta) \alpha'(\eta), \varphi \bigr)_H, \mbox{ for any $ \varphi \in V $.}
    \end{align*}
    Similarly, for any $ \tilde{\eta} \in V $, one can see that $ [\theta, \theta^*] \in \partial_\theta \Phi_0(\tilde{\eta}, \cdot{}) $ in $ H \times H $, iff.:
    \begin{equation}\label{var_th}
        \begin{array}{c}
            \displaystyle \bigl( -\theta^*, \theta -\psi \bigr)_H +\nu^2 \bigl( \partial_x \theta, \partial_x(\theta -\psi) \bigr)_H + \int_\Omega \alpha(\tilde{\eta})|\partial_x \theta| \, dx
            \\[2ex]
            \displaystyle \leq \int_\Omega \alpha(\tilde{\eta})|\partial_x \psi| \,dx, \mbox{ for any $ \psi \in V_{0} $.}
        \end{array}
    \end{equation}
    Taking into account \eqref{AA_0}--\eqref{var_th}, we deduce that the variational problem as in \eqref{S1}--\eqref{S3} is equivalently reformulated to the following Cauchy problem:
\begin{equation*}
    \mbox{($\widetilde{\mathrm{E}}$)}~~
    \begin{cases}
        \mathcal{A}(t) w'(t) +\bigl[ \partial_\eta \Phi_0(w(t)) \times \partial_\theta \Phi_0(w(t)) \bigr] +\mathcal{G}(w(t)) \ni \mathfrak{f}(t) \mbox{ in $ [H]^2 $, $ t \in (0, T) $,}
        \\
        w(0) = w_0 \mbox{ in $[H]^2$.}
    \end{cases}
\end{equation*}
The item $ (\star\,0) $ is a straightforward consequence of this reformulation and the inclusion $ \partial \Phi_0 \subset \bigl[ \partial_\eta \Phi_0 \times \partial_\theta \Phi_0 \bigr] $ in $ [H]^2 \times [H]^2 $, mentioned in \eqref{prodSubDif}. 
\end{remark}
\medskip

Now, we are ready to prove the Main Theorem \ref{mainTh01}.
\paragraph{\textbf{Proof of Main Theorem \ref{mainTh01} (I-A)}}{
    First, we verify the existence part. Under the setting \eqref{AA_0}--\eqref{Phi_eps}, we immediately check that:
    \begin{itemize}
        \item[\textmd{(ev.0)}]for any $ t \in [0, T] $, $ \mathcal{A}(t) \in \mathscr{L}([H]^2) $ is positive and selfadjoint, and  
            \begin{equation*}
                (\mathcal{A}(t) w, w)_{[H]^2} \geq \delta_* |w|_{[H]^2}^2, \mbox{ for any $ w \in [H]^2 $,}
            \end{equation*}
            with the constant $ \delta_* \in (0, 1) $ as in (A3);
        \item[\textmd{(ev.1)}]$ \mathcal{A} \in W^{1, \infty}(0, T; \mathscr{L}([H]^2)) $, and  
            \begin{equation*}
                A^* := \mathrm{ess} \sup_{\hspace{-3ex}t \in (0, T)} \left\{ \max \{ |\mathcal{A}(t)|_{\mathscr{L}([H]^2)}, |\mathcal{A}'(t)|_{\mathscr{L}([H]^2)} \} \right\} \leq 1 +|\alpha_0|_{W^{1, \infty}(Q)} < \infty;
            \end{equation*}
        \item[\textmd{(ev.2)}]$ \mathcal{G} : [H]^2 \longrightarrow [H]^2 $ is a Lipschitz continuous operator with a Lipschitz constant:
            \begin{equation*}
                L_* := 1 +|g'|_{L^\infty(\mathbb{R})} +\nu^{-2} {\textstyle{\bigl| \frac{d}{d\eta}(\alpha \alpha') \bigr|_{L^\infty(\mathbb{R})}},}
            \end{equation*}
            and $ \mathcal{G} $ has a $ C^1 $-potential functional
            \begin{align*}
                \widehat{\mathcal{G}}: w =  [\eta, \theta] \in [H]^2 & \mapsto \widehat{\mathcal{G}}(w) := \int_\Omega \left( G(\eta) -\frac{\eta^2}{2} -\frac{\alpha(\eta)^2}{2 \nu^2} \right) \, dx \in \mathbb{R};
            \end{align*}
        \item[\textmd{(ev.3)}]$ \Phi_\varepsilon \geq 0 $ on $ [H]^2 $, and the sublevel set $ \bigl\{  \tilde{w} \in [H]^2 \, \bigl| \, \Phi_\varepsilon(\tilde{w}) \leq r \bigr\} $ is contained in a compact set $ K_\nu(r) $ in $ [H]^2 $, defined as
            \begin{equation*}
                K_\nu(r) := \left\{ \begin{array}{l|l}
                    \tilde{w} = [\tilde{\eta}, \tilde{\theta}] \in V \times V_0 &  |\tilde{\eta}|_V^2 \leq 2r \mbox{ and } |\tilde{\theta}|_{V_0}^2 \leq 2 \nu^{-2} r
            \end{array} \right\}, 
        \end{equation*}
        for any $ r \geq 0 $.
    \end{itemize}
    On account of \eqref{AA_0}--\eqref{Phi_eps} and (ev.0)--(ev.3), we can apply Lemma \ref{Lem.CP} in Appendix, as the case when:
    \begin{equation*}
        X = [H]^2, ~ \mathcal{A}_0 = \mathcal{A} \mbox{ in $ W^{1, \infty}(0, T; \mathscr{L}([H]^2)) $, } \mathcal{G}_0 = \mathcal{G} \mbox{ on $ [H]^2 $, and } \Psi_0 = \Phi_\varepsilon \mbox{ on $ [H]^2 $,}
    \end{equation*}
    and we can find a solution $ w = [\eta, \theta] \in [\mathscr{H}]^2 $ to the Cauchy problem (E)$_\varepsilon$. In the light of Remarks \ref{Rem.eps>0} and \ref{Rem.Sec03-01}, finding this $ w = [\eta, \theta] $ directly leads to the existence of solution to the state-system (S)$_\varepsilon$. 
    \medskip

    Next, for the verification of the uniqueness part, we suppose that the both pairs of functions $ [\eta^\ell, \theta^\ell] \in [\mathscr{H}]^2 $,  $\ell = 1, 2$, solve the state-system (S)$_\varepsilon$ for the common initial pair $[\eta_0, \theta_0]$ and forcing pair $ [u, v] \in [\mathscr{H}]^2 $. Besides, let us take the difference between two variational forms \eqref{S2} for $ \eta^\ell $, $ \ell = 1, 2 $, and put $ \varphi = \eta^1 -\eta^2 $. Then, by using the assumptions (A1)--(A3), and H\"{o}lder's and Young's inequalities, we have:
    \begin{subequations}\label{IA-00}
    \begin{equation}\label{IA-01}
        \frac{1}{2} \frac{d}{dt} |(\eta^1 -\eta^2)(t)|_H^2 +|\partial_x (\eta^1 -\eta^2)(t)|_{H}^2  =  I_A^1 +I_A^2,
    \end{equation}
    with
    \begin{align}
        I_A^1 & := -\bigl(g(\eta^1(t)) -g(\eta^2(t)), (\eta^1 -\eta^2)(t) \bigr)_H \leq L_* |(\eta^1 -\eta^2)(t)|_H^2,
        \label{IA02}
    \end{align}
    and
    \begin{align}
        I_A^2  := & -\bigl( \alpha'(\eta^1(t)) f_\varepsilon(\partial_x \theta^1(t)) -\alpha'(\eta^2(t)) f_\varepsilon(\partial_x \theta^2(t)), (\eta^1 -\eta^2)(t) \bigr)_H
        \nonumber
        \\
        = & \displaystyle \int_\Omega f_\varepsilon(\partial_x \theta^1(t)) \bigl( \alpha'(\eta^1(t))(\eta^2 -\eta^1)(t) \bigr) \, dx
        \nonumber
        \\
        & \displaystyle \qquad +\int_\Omega f_\varepsilon(\partial_x \theta^2(t)) \bigl( \alpha'(\eta^2(t))(\eta^1 -\eta^2)(t) \bigr) \, dx
        \nonumber
        \\
        \leq & \displaystyle -\int_\Omega \bigl( f_\varepsilon(\partial_x \theta^1(t)) -f_\varepsilon(\partial_x \theta^2(t)) \bigr) \bigl( \alpha(\eta^1(t)) -\alpha(\eta^2(t)) \bigr) \, dx
        \nonumber
        \\
        \leq & |\alpha'|_{L^\infty(\mathbb{R})}|\partial_x(\theta^1 -\theta^2)(t)|_H |(\eta^1 -\eta^2)(t)|_H
        \nonumber
        \\
        \leq & \displaystyle \frac{\nu^2}{4} |\partial_x (\theta^1 -\theta^2)(t)|_H^2 +\frac{|\alpha'|_{L^\infty(\mathbb{R})}^2}{\nu^2} |(\eta^1 -\eta^2)(t)|_H^2, \mbox{for a.e. $ t \in (0, T) $.}
        \label{IA-03}
    \end{align}
    \end{subequations}
    Meanwhile, for any $ \ell \in \{1, 2\} $, let us take $ \ell^\perp \in \{1, 2\} \setminus \{\ell\}  $, and put $ \psi =\theta^{\ell^\perp} $ in the variational inequality \eqref{S3} for $ \theta^\ell $. Then, adding those two variational inequalities, and using H\"{o}lder's and Young's inequalities, one can observe that:
    \begin{subequations}\label{IA-10}
    \begin{equation}\label{IA-11}
        \frac{1}{2} \frac{d}{dt} |\sqrt{\alpha_0(t)}(\theta^1 -\theta^2)(t)|_H^2 +\nu^2 |\partial_x(\theta^1 -\theta^2)(t)|_H^2 \leq I_A^3 +I_A^4, 
    \end{equation}
    with
    \begin{align}
        I_A^3 & := \frac{1}{2} \int_\Omega \partial_t \alpha_0 (t)|(\theta^1 -\theta^2)(t)|^2 \, dx \leq \frac{|\partial_t \alpha_0|_{L^\infty(Q)}}{2} |(\theta^1 -\theta^2)(t)|_H^2,
        \label{IA12}
    \end{align}
    and
    \begin{align}
        I_A^4  := & \displaystyle -\int_\Omega \alpha(\eta^2(t)) f_\varepsilon(\partial_x \theta^2(t)) \, dx +\int_\Omega \alpha(\eta^2(t)) f_\varepsilon(\partial_x \theta^1(t))
        \nonumber
        \\
        & \displaystyle \qquad \displaystyle -\int_\Omega \alpha(\eta^1(t)) f_\varepsilon(\partial_x \theta^1(t)) \, dx +\int_\Omega \alpha(\eta^1(t)) f_\varepsilon(\partial_x \theta^2(t))
        \nonumber
        \\
        = & \displaystyle -\int_\Omega \bigl( f_\varepsilon(\partial_x \theta^1(t)) -f_\varepsilon(\partial_x \theta^2(t)) \bigr) \bigl( \alpha(\eta^1(t)) -\alpha(\eta^2(t)) \bigr) \, dx
        \nonumber
        \\
        \leq & \displaystyle \frac{\nu^2}{4} |\partial_x (\theta^1 -\theta^2)(t)|_H^2 +\frac{|\alpha'|_{L^\infty(\mathbb{R})}^2}{\nu^2} |(\eta^1 -\eta^2)(t)|_H^2, \mbox{for a.e. $ t \in (0, T) $.}
        \label{IA-13}
    \end{align}
    \end{subequations}
    As the summation of \eqref{IA-00} and \eqref{IA-10}, we obtain that:
    \begin{align}\label{IA-20}
        \displaystyle \frac{1}{2} \frac{d}{dt} \bigl( |(\eta^1 & -\eta^2)(t)|_H^2 +|\sqrt{\alpha_0(t)}(\theta^1 -\theta^2)(t)|_H^2 \bigr) 
        \nonumber
        \\
        \leq ~ & C_A^1 \bigl( |(\eta^1 -\eta^2)(t)|_H^2 +|\sqrt{\alpha_0(t)}(\theta^1 -\theta^2)(t)|_H^2 \bigr),
        \\
        & \mbox{for a.e. $ t \in (0, T) $,  with } C_A^1 := L_* + \frac{2|\alpha'|_{L^\infty(\mathbb{R})}^2}{\nu^2} +\frac{|\partial_t \alpha_0|_{L^\infty(Q)}}{2 \delta_*}.
        \nonumber
    \end{align}
    
    Now, with (A3) in mind, we can verify the uniqueness part of (I-A), just by applying Gronwall's lemma to the estimate \eqref{IA-20}. 
    \qed
}
\medskip
\begin{remark}\label{Rem.sols}
    As a consequence of the uniqueness result in (I-A), we can say that the converse of $ (\star\,0) $ in Remark \ref{Rem.Sec03-01} is also true, i.e. the three problems (S)$_0$, (E)$_0$, and ($\widetilde{\mathrm{E}}$) are equivalent each other.
\end{remark}
\paragraph{\textbf{Proof of Main Theorem \ref{mainTh01} (I-B)}}{
    By Remarks \ref{Rem.eps>0}--\ref{Rem.sols}, the solution $ w := [\eta, \theta] \in [\mathscr{H}]^2 $ to the state-system (S)$_\varepsilon$ coincides with that to the Cauchy problem (E)$_\varepsilon$ for the initial data $ w_0 := [\eta_0, \theta_0] \in V \times V_0 $ and forcing term $ \mathfrak{f} := [M_u u, M_v v] \in [\mathscr{H}]^2 $. Also, putting:
    \begin{equation*}
        \begin{array}{c}
            w_n := [\eta_n, \theta_n] \mbox{ in $ [\mathscr{H}]^2 $, } w_{0, n} := [\eta_{0, n}, \theta_{0, n}] \mbox{ in $ [H]^2 $,} 
            \\
            \mbox{and } \mathfrak{f}_n := [M_u u_n, M_v v_n] \mbox{ in $ [\mathscr{H}]^2 $, for }n = 1, 2, 3, \dots,
        \end{array}
    \end{equation*}
    we can suppose that the sequence $ \{ w_n \}_{n = 1}^\infty = \{ [\eta_n, \theta_n] \}_{n = 1}^\infty $ of solutions to systems (S)$_{\varepsilon_n}$, $ n = 1, 2, 3, \dots $, coincides with that of solutions to the problems (E)$_{\varepsilon_n} $, for the initial data $ w_{0, n} $ and forcing terms $ \mathfrak{f}_n $, $ n = 1, 2, 3, \dots $. In addition:
    
\begin{itemize}
    \item[\textmd{(ev.4)}]$ \Phi_{\varepsilon_n} \geq 0 $ on $ [H]^2 $, for $ n = 1, 2, 3, \dots $, and the union $ \bigcup_{n = 1}^\infty \bigl\{  \tilde{w} \in [H]^2 \, \bigl| \, \Phi_{\varepsilon_n}(\tilde{w}) \leq r \bigr\} $ of sublevel sets is contained in the compact set $ K_\nu(r) \subset [H]^2 $, as in (ev.3), for any $ r > 0 $;
    \item[\textmd{(ev.5)}]$ \Phi_{\varepsilon_n} \to \Phi_\varepsilon $ on $ [H]^2 $, in the sense of Mosco, as $ n \to \infty $, more precisely, the following estimate 
        \begin{align}\label{ev.M01}
            |\Phi_{\varepsilon_n}({w}) & -\Phi_{\varepsilon}({w})| 
            \nonumber
            \\
            = & \frac{1}{2} \left| \int_\Omega \left( \bigl( \nu f_{\varepsilon_n}(\partial_x \theta) +\nu^{-1} \alpha(\eta) \bigr)^2 -\bigl( \nu f_{\varepsilon}(\partial_x \theta) +\nu^{-1} \alpha(\eta) \bigr)^2 \right) \, dx  \right|
            \nonumber
            \\
            \leq & \frac{\nu}{2} \int_\Omega \bigl| \nu \bigl( f_{\varepsilon_n}(\partial_x \theta) +f_{\varepsilon}(\partial_x \theta) \bigr) +2 \nu^{-1} \alpha(\eta) \bigr| \bigl| f_{\varepsilon_n}(\partial_x \theta) -f_\varepsilon(\partial_x \theta) \bigr| \, dx
            \nonumber
            \\
            \leq & \frac{\nu^2}{2} |\varepsilon_n -\varepsilon| \int_\Omega \left( (\varepsilon_n +\varepsilon) +2|\partial_x {\theta}| +\frac{2}{\nu^2} \alpha({\eta}) \right) \, dx
            \nonumber
            \\
            \leq & \nu^2 |\varepsilon_n -\varepsilon| \int_\Omega \left( 1 +|\partial_x {\theta}| +\frac{1}{\nu^2} \alpha({\eta}) \right) \, dx,
            \nonumber
            \\
            & \mbox{for any $ {w} = [{\eta}, {\theta}] \in V \times V_0 $, $ n = 1, 2, 3, \dots $,}
        \end{align}
        where we use the following inequality:
\begin{align*}
      |f_\varepsilon(\omega) - f_\delta(\omega)| & = \left|\frac{\varepsilon^2 - \delta^2}{\sqrt{\varepsilon^2 + |\omega|^2} + \sqrt{\delta^2 + |\omega|^2}} \right|
      \\
      & = \frac{|\varepsilon + \delta|}{\sqrt{\varepsilon^2 + |\omega|^2} + \sqrt{\delta^2 + |\omega|^2}}|\varepsilon - \delta|
      \\
      & \leq |\varepsilon - \delta|, \ \mbox{for any}\ \varepsilon, \delta \in [0, 1],\ \mbox{and}\ \omega \in \mathbb{R}. 
\end{align*}
        Immediately leads to the corresponding lower bound condition and optimality condition, in the Mosco-convergence of $ \{ \Phi_{\varepsilon_n} \}_{n=1}^{\infty} $;
    \item[\textmd{(ev.6)}]$ \sup_{n \in \mathbb{N}} \Phi_{\varepsilon_n}(w_{0, n}) < \infty $, and 
        \begin{align*}
        w_{0, n} \to w_{0} \ \mbox{in}\ [H]^{2},\ \mbox{as}\ n \to \infty,
        \end{align*}        
        more precisely, it follows from \eqref{w.i} and (A3) that
        \begin{align*}
            \sup_{n \in \mathbb{N}} \Phi_{\varepsilon_n}(w_{0, n}) & \leq \sup_{n \in \mathbb{N}} \left( \frac{1}{2}|\eta_{0, n}|_V^2 +\nu^2 ( 1 +|\theta_{0, n}|_{V_0}^2) +\frac{1}{\nu^2} |\alpha(\eta_{0, n})|_H^2 \right) < \infty,
        \end{align*}
        and moreover, the weak convergence of $\{w_{0, n} \}_{n=1}^{\infty}$ in $V \times V_{0}$ and the compactness of embedding $V \times V_{0} \subset [H]^{2}$ imply the strong convergence of $\{w_{0, n} \}_{n=1}^{\infty}$ in $[H]^{2}$.
\end{itemize}

On account of \eqref{w.i} and (ev.0)--(ev.6), we can apply Lemma \ref{Lem.CP02}, to show that:
\begin{subequations}\label{convKS}
\begin{equation}\label{convKS01}
    \begin{cases}
        w_n \to w \mbox{ in $ C([0, T]; [H]^2) $ (i.e. in $ [C([0, T]; H)]^2 $),}
        \\
        \quad \mbox{weakly in $ W^{1, 2}(0, T; [H]^2) $ (i.e. weakly in $ [W^{1, 2}(0, T; H)]^2 $),}
        \\[1ex]
        \displaystyle \int_0^T \Phi_{\varepsilon_n}(w_n(t)) \, dt \to \int_0^T \Phi_{\varepsilon}(w(t)) \, dt,
    \end{cases}
    \mbox{as $ n \to \infty $,}
\end{equation}
\begin{align*}
    \displaystyle \sup_{n \in \mathbb{N}} |w_n| & _{L^\infty(0, T; V) \times L^\infty(0, T; V_0)}^2 \leq 4 \sup_{n \in \mathbb{N}} |w_n|_{L^\infty(0, T; V \times V_0)}^2 
    \\
    & \displaystyle \leq \frac{8}{\min \, \{1, \nu^2 \}} \sup_{n \in \mathbb{N}} |\Phi_{\varepsilon_n}(w_n)|_{L^\infty(0, T)} < \infty, 
\end{align*}
and hence,
\begin{equation}\label{convKS02}
    w_n \to w \mbox{ weakly-$*$ in $ L^\infty(0, T; V) \times L^\infty(0, T; V_0) $, as $ n \to \infty $.}
\end{equation}
\end{subequations}
Also, as a consequence of the one-dimensional compact embeddings $ V \subset C(\overline{\Omega}) $ and $ V_0 \subset C(\overline{\Omega}) $, the uniqueness of solution $ w $ to (E)$_\varepsilon$, and Ascoli's theorem (cf. \cite[Corollary 4]{MR0916688}), we can derive from \eqref{convKS01} that
\begin{equation}\label{convKS07s}
    w_n \to w \mbox{ in $ [C(\overline{Q})]^2 $,\ as $ n \to \infty $.}
\end{equation}
Furthermore, from \eqref{convKS}, \eqref{convKS07s}, and the assumptions (A1) and (A3), one can observe that:
\begin{subequations}\label{convKSs}
\begin{equation}\label{convKS03}
    \begin{cases}
        \displaystyle \varliminf_{n \to \infty} \frac{1}{2} | \eta_n |_{\mathscr{V}}^2 \geq \frac{1}{2} |\eta|_{\mathscr{V}}^2, \quad \varliminf_{n \to \infty} \frac{\nu^2}{2} |\theta_n |_{\mathscr{V}_0}^2 \geq \frac{\nu^2}{2} |\theta|_{\mathscr{V}_0}^2,
        \\[2ex]
        \displaystyle \lim_{n \to \infty} \frac{1}{2 \nu^2} |\alpha(\eta_n)|_{\mathscr{H}}^2 = \frac{1}{2 \nu^2} |\alpha(\eta)|_{\mathscr{H}}^2,
    \end{cases}
\end{equation}
and
\begin{align}\label{convKS04}
    \displaystyle \varliminf_{n \to \infty} & \displaystyle \bigl| \alpha(\eta_n) f_{\varepsilon_n}(\partial_x \theta_n) \big|_{L^1(Q)} = \varliminf_{n \to \infty} \int_0^T \int_\Omega \alpha(\eta_n(t)) f_{\varepsilon_n}(\partial_x \theta_n(t)) \, dx dt
    \nonumber
    \\
    & \geq \varliminf_{n \to \infty} \int_0^T \int_\Omega \alpha(\eta(t)) f_{\varepsilon_n}(\partial_x \theta_n(t)) \, dx dt
    \nonumber
    \\
    & \qquad -\lim_{n \to \infty} |\alpha(\eta_n) -\alpha(\eta)|_{C(\overline{Q})} \cdot \sup_{n \in \mathbb{N}} \bigl( T \varepsilon_n +|\partial_x \theta_n|_{L^1(0, T; L^1(\Omega))} \bigr)
    \nonumber
    \\
    & \geq \varliminf_{n \to \infty} \int_0^T \int_\Omega \alpha(\eta(t)) f_{\varepsilon}(\partial_x \theta_n(t)) \, dx dt -|\alpha(\eta)|_{C(\overline{Q})} \cdot \lim_{n \to \infty} \bigl( T|\varepsilon_n -\varepsilon| \bigr)
    \nonumber
    \\
    & \geq \int_0^T \int_\Omega \alpha(\eta(t)) f_{\varepsilon}(\partial_x \theta(t)) \, dx dt = \bigl| \alpha(\eta) f_{\varepsilon}(\partial_x \theta) \big|_{L^1(Q)}.
\end{align}
\end{subequations}
Here, from \eqref{Phi_eps}, it is seen that:
\begin{align}\label{convKS'}
    & \displaystyle \int_0^T {\Phi}_{\tilde{\varepsilon}}(\tilde{w}(t)) \, dt =  \int_0^T {\Phi}_{\tilde{\varepsilon}}(\tilde{\eta}(t), \tilde{\theta}(t)) \, dt
    \nonumber
    \\
    & \qquad = \frac{1}{2} |\tilde{\eta}|_{\mathscr{V}}^2 +\frac{\nu^2}{2} |\tilde{\theta}|_{\mathscr{V}_0}^2 +\bigl| \alpha(\tilde{\eta}) f_{\tilde{\varepsilon}}(\partial_x \tilde{\theta}) \bigr|_{L^1(Q)} +\frac{1}{2 \nu^2} |\alpha(\tilde{\eta})|_{\mathscr{H}}^2 +\frac{\nu^2 \tilde{\varepsilon}^2}{2} T
    \nonumber
    \\[1ex]
    & \qquad \qquad \mbox{for all $ \tilde{\varepsilon} > 0 $ and $ \tilde{w} = [\tilde{\eta}, \tilde{\theta}] \in D(\Phi_{\tilde{\varepsilon}}) = \mathscr{Y} $.}
\end{align}
Taking into account \eqref{convKS01}, \eqref{convKSs}, and \eqref{convKS'}, we deduce that:
\begin{equation}\label{convKS05}
        |\eta_n|_\mathscr{V}^2 +\nu^2 |\theta_n|_{\mathscr{V}_0}^2 \to |\eta|_\mathscr{V}^2 +\nu^2 |\theta|_{\mathscr{V}_0}^2,
        \mbox{and hence, } |w_n|_{\mathscr{Y}} \to |w|_{\mathscr{Y}}, \mbox{ as $ n \to \infty $.}
\end{equation}
Since the norm of Hilbert space $ \mathscr{Y} $ is uniformly convex, the convergences \eqref{convKS02} and \eqref{convKS05} imply the strong convergences:
\begin{subequations}\label{convKS06-07}
\begin{equation}\label{convKS06}
    w_n \to w \mbox{ in $ \mathscr{Y} $, as $ n \to \infty $,}
\end{equation}
and 
\begin{align}\label{convKS07}
    |f_{\varepsilon_n} (\partial_x \theta_n) & -f_\varepsilon(\partial_x \theta)|_\mathscr{H} \leq  |f_{\varepsilon_n}(\partial_x \theta_n) -f_{\varepsilon_n}(\partial_x \theta)|_\mathscr{H} +|f_{\varepsilon_n}(\partial_x \theta) -f_{\varepsilon}(\partial_x \theta)|_\mathscr{H} 
    \nonumber
    \\
    \leq &  |\theta_n -\theta|_{\mathscr{V}_0} +\sqrt{T} |\varepsilon_n -\varepsilon| \to 0, \mbox{ as $ n \to \infty $.}
\end{align}
\end{subequations}
The convergences \eqref{convKS} and \eqref{convKS06-07} are sufficient to verify the conclusions \eqref{mThConv} and \eqref{mThConv00} of Main Theorem \ref{mainTh01} (I-B). 
\qed
}

\section{Proof of Main Theorem \ref{mainTh02}}

In this section, we prove the second Main Theorem \ref{mainTh02}. Let $[\eta_0, \theta_0] \in V \times V_0$ be the initial pair, fixed in (A1). Also, let us fix arbitrary forcing pair  $ [\bar{u}, \bar{v}] \in [\mathscr{H}]^2 $, and let us invoke the definition of the cost function \eqref{J}, to estimate that:
\begin{align}\label{mTh02-00}
    0 &~ \leq \underline{J}_{\varepsilon}:= \inf \mathcal{J}_{\varepsilon}([\mathscr{H}]^{2}) \leq \overline{J}_{\varepsilon}:= \mathcal{J}_\varepsilon(\bar{u}, \bar{v}) < \infty, \mbox{ for all } \varepsilon \geq 0.
\end{align}
Also, for any $ \varepsilon \geq 0 $, we denote by $[\bar{\eta}_{\varepsilon}, \bar{\theta}_{\varepsilon}]$ the solution to (S)$_{\varepsilon}$, for the initial pair $[\eta_0, \theta_0]$ and forcing pair $ [\bar{u}, \bar{v}]$.
\bigskip

Based on these, the proof of Main Theorem \ref{mainTh02} is demonstrated as follows.

\paragraph{\textbf{Proof of Main Theorem \ref{mainTh02} (II-A)}}
Let us fix any $\varepsilon \geq 0$. 
Then, from the estimate \eqref{mTh02-00}, we immediately find a sequence of forcing pairs $\{[u_{n}, v_{n}] \}_{n=1}^{\infty} \subset [\mathscr{H}]^{2}$, such that:
\begin{subequations}\label{mTh02-0102}
\begin{equation}\label{mTh02-01}
\mathcal{J}_{\varepsilon}(u_{n}, v_{n}) \downarrow \underline{J}_{\varepsilon},\ \mbox{as}\ n \to \infty,
\end{equation}
and
\begin{align}\label{mTh02-02}
    \sup_{n \in \mathbb{N}}& \bigl| [\sqrt{M_u} u_{n}, \sqrt{M_v} v_{n}] \bigr|_{[\mathscr{H}]^{2}}^{2} \leq \mathcal{J}_\varepsilon(\bar{u}, \bar{v}) < \infty.
\end{align}
\end{subequations}
Also, the estimate \eqref{mTh02-02} enables us to take a subsequence of $\{[u_{n}, v_{n}] \}_{n=1}^{\infty} \subset [\mathscr{H}]^{2}$ (not relabeled), and to find a pair of functions $[u^{*}, v^{*}] \in [\mathscr{H}]^{2}$, such that:
\begin{equation*}
[\sqrt{M_u} u_{n}, \sqrt{M_v} v_{n}] \to [\sqrt{M_u} u^{*}, \sqrt{M_v} v^{*}]\ \mbox{weakly in}\ [\mathscr{H}]^{2}, \mbox{as}\ n \to \infty, 
\end{equation*}
and as well as,
\begin{equation}\label{mTh02-03}
[M_u u_{n}, M_v v_{n}] \to [M_u u^{*}, M_v v^{*}]\ \mbox{weakly in}\ [\mathscr{H}]^{2}, \ \mbox{as}\ n \to \infty.
\end{equation}

Let $[\eta^{*}, \theta^{*}] \in [\mathscr{H}]^{2}$ be the solution to (S)$_{\varepsilon}$, for the initial pair $[\eta_0, \theta_0]$ and forcing pair $[u^{*}, v^{*}]$. 
As well as, for any $n \in \mathbb{N}$, let $[\eta_{n}, \theta_{n}] \in [\mathscr{H}]^{2}$ be the solution to (S)$_{\varepsilon}$, for the forcing pair $[u_{n}, v_{n}]$. 
Then, having in mind \eqref{mTh02-03} and the initial condition:
\begin{equation*}
    [\eta_{n}(0), \theta_{n}(0)] = [\eta^{*}(0), \theta^{*}(0)] = [\eta_{0}, \theta_{0}] \mbox{ in $[H]^{2}$, for $n = 1, 2, 3, \dots $,}
\end{equation*}
we can apply Main Theorem \ref{mainTh01} (I-B), to see that:
\begin{equation}\label{mTh02-04}
[\eta_{n}, \theta_{n}] \to [\eta^{*}, \theta^{*}]\ \mbox{in}\ [C(\overline{Q})]^{2},\ \mbox{as}\ n \to \infty.
\end{equation}
On account of \eqref{mTh02-01}, \eqref{mTh02-03}, and \eqref{mTh02-04}, it is computed that:
\begin{align*}
    \mathcal{J}_{\varepsilon}(u^{*}, v^{*}) & = \frac{1}{2} \bigl| [\sqrt{M_\eta} (\eta^{*}-\eta_\mathrm{ad}), \sqrt{M_\theta} (\theta^{*}-\theta_\mathrm{ad}) ] \bigr|_{[\mathscr{H}]^{2}}^{2}\\
                                            &\qquad + \frac{1}{2} \bigl| [\sqrt{M_u} u^{*}, \sqrt{M_v} v^{*}] \bigr|_{[\mathscr{H}]^{2}}^{2}
    \\
                                            & \leq \frac{1}{2}\lim_{n \to \infty} \bigl| [\sqrt{M_\eta} (\eta_{n}-\eta_\mathrm{ad}), \sqrt{M_\theta} (\theta_{n}-\theta_\mathrm{ad}) ] \bigr|_{[\mathscr{H}]^{2}}^{2}\\
                                            &\qquad + \frac{1}{2} \varliminf_{n \to \infty} \bigl| [\sqrt{M_u} u_{n}, \sqrt{M_v} v_{n}] \bigr|_{[\mathscr{H}]^{2}}^{2}
    \\
    & = \lim_{n \to \infty}\mathcal{J}_{\varepsilon}(u_{n}, v_{n}) = \underline{J}_{\varepsilon} ~ (\leq \mathcal{J}_{\varepsilon}(u^{*}, v^{*})),
\end{align*}
and it implies that
\begin{equation*}
    \mathcal{J}_\varepsilon(u^{*}, v^{*}) = \min_{[u, v] \in [\mathscr{H}]^{2}} \mathcal{J}_{\varepsilon}(u, v).
\end{equation*}

Thus, we conclude the item (II-A). 
\qed

\paragraph{\textbf{Proof of Main Theorem \ref{mainTh02} (II-B)}}
Let $ \varepsilon  \in [0, 1] $ and $\{\varepsilon_{n} \}_{n=1}^{\infty} \subset [0, 1]$ be as in \eqref{ass.4}. Let $[\bar{\eta}_{\varepsilon}, \bar{\theta}_{\varepsilon}] \in [\mathscr{H}]^{2}$ be the solution to the system (S)$_\varepsilon$, for the initial pair $ [\eta_0, \theta_0] $ and forcing pair $ [\bar{u}, \bar{v}] $, and let $[\bar{\eta}_{\varepsilon_{n}}, \bar{\theta}_{\varepsilon_{n}}] \in [\mathscr{H}]^{2}$, $n=1, 2, 3,\ldots$, be the solutions to (S)$_{\varepsilon_n}$, for the respective initial pairs $ [\eta_{0, n}, \theta_{0, n}] $, $ n = 1, 2, 3, \dots $, and the fixed forcing pair $ [\bar{u}, \bar{v}] $. On this basis, let us first apply Main Theorem \ref{mainTh01} (I-B) to the solutions $[\bar{\eta}_{\varepsilon}, \bar{\theta}_{\varepsilon}] \in [\mathscr{H}]^{2}$ and $[\bar{\eta}_{\varepsilon_{n}}, \bar{\theta}_{\varepsilon_{n}}] \in [\mathscr{H}]^{2}$, $n=1, 2, 3,\ldots$. Then, we have 
\begin{equation}\label{mTh02-09}
    \begin{cases}
        [\bar{\eta}_{\varepsilon_{n}}, \bar{\theta}_{\varepsilon_{n}}] \to [\bar{\eta}_{\varepsilon}, \bar{\theta}_{\varepsilon}]\ \mbox{in}\ [C(\overline{Q})]^{2},
        \\[1ex]
        [\bar{\eta}_{n}(0), \bar{\theta}_{n}(0)] = [\eta_{0, n}, \theta_{0, n}] 
        \\
        \hspace{3.15ex}\to [\eta_0, \theta_0] = [\bar{\eta}_{\varepsilon}(0), \bar{\theta}_{\varepsilon}(0)] \mbox{ in $ [C(\overline{\Omega})]^2 $,}
    \end{cases}
    \mbox{as}\ n \to \infty,
\end{equation}
and hence, 
\begin{equation}\label{mTh02-10}
    \overline{J}_\mathrm{sup} := \sup_{n \in \mathbb{N}} \mathcal{J}_{\varepsilon_n}(\bar{u}, \bar{v}) < \infty.
\end{equation}

Next, for any $n \in \mathbb{N}$, let us denote by $[\eta^{*}_{n}, \theta^{*}_{n}] \in [\mathscr{H}]^{2}$ the solution to (S)$_{\varepsilon_{n}}$, for the initial pair $[\eta_{0, n}, \theta_{0, n}]$ and forcing pair $[u^{*}_{n}, v^{*}_{n}]$. Then, in the light of \eqref{mTh02-00} and \eqref{mTh02-10}, we can see that:
\begin{align*}
    0 & \leq \frac{1}{2}|[\sqrt{M_u} u^{*}_{n}, \sqrt{M_v} v^{*}_{n}] |_{[\mathscr{H}]^{2}}^{2} \leq \underline{J}_{\varepsilon_{n}} \leq \overline{J}_\mathrm{sup} < \infty,  \mbox{ for } n = 1, 2, 3, \dots.
\end{align*}
Therefore, we can find a subsequence $\{n_{i} \}_{i=1}^{\infty} \subset \{n\}$, together with a pair of functions $[u^{**}, v^{**}] \in [\mathscr{H}]^{2}$, such that:
\begin{equation*}
[\sqrt{M_u} u^{*}_{n_i}, \sqrt{M_v} v^{*}_{n_i}] \to [\sqrt{M_u} u^{**}, \sqrt{M_v} v^{**}]\ \mbox{weakly in}\ [\mathscr{H}]^{2},\ \mbox{as}\ i \to \infty, 
\end{equation*}
and as well as,
\begin{equation}\label{mTh02-11}
    [M_u u^{*}_{n_i}, M_v v^{*}_{n_i}] \to [M_u u^{**}, M_v v^{**}]\ \mbox{weakly in}\ [\mathscr{H}]^{2}, \ \mbox{as}\ i \to \infty.
\end{equation}
Here, let us denote by $[\eta^{**}, \theta^{**}] \in  [\mathscr{H}]^2 $ the solution to (S)$_{\varepsilon}$, for the initial pair $ [\eta_0, \theta_0] $ and forcing pair $[u^{**}, v^{**}]$.
Then, applying Main Theorem \ref{mainTh01} (I-B), again, to the solutions $[\eta^{**}, \theta^{**}]$ and $[\eta^{*}_{n_i}, \theta^{*}_{n_i}]$, $i = 1, 2, 3, \dots$, we can observe that:
\begin{equation}\label{mTh02-12}
[\eta^{*}_{n_i}, \theta^{*}_{n_i}] \to [\eta^{**}, \theta^{**}] \mbox{ in } [C(\overline{Q})]^{2}, \mbox{ as } i \to \infty.
\end{equation}

Now, as a consequence of \eqref{mTh02-09}, \eqref{mTh02-11}, and \eqref{mTh02-12}, it is verified that:
\begin{align*}
    \mathcal{J}_{\varepsilon} (u^{**}, v^{**}) &= \frac{1}{2} \bigl| [\sqrt{M_\eta} (\eta^{**}-\eta_\mathrm{ad}), \sqrt{M_\theta} (\theta^{**}-\theta_\mathrm{ad})] \bigr|_{[\mathscr{H}]^{2}}^{2}
    \\
    &\qquad + \frac{1}{2} \bigl| [\sqrt{M_u} u^{**}, \sqrt{M_v} v^{**}] \bigr|_{[\mathscr{H}]^{2}}^{2}
    \\
    & \leq \frac{1}{2} \lim_{i \to \infty} \bigl| [\sqrt{M_\eta} (\eta^{*}_{n_i} -\eta_\mathrm{ad}), \sqrt{M_\theta} (\theta^{*}_{n_i} -\theta_\mathrm{ad})] \bigr|_{[\mathscr{H}]^{2}}^{2}
    \\
    &\qquad  + \frac{1}{2}\varliminf_{i \to \infty} \bigl| [\sqrt{M_u} u^{*}_{n_i}, \sqrt{M_v} v^{*}_{n_i}] \bigr|_{[\mathscr{H}]^{2}}^{2}
    \\
    & \leq \varliminf_{i \to \infty}\mathcal{J}_{\varepsilon_{n_i}}(u^{*}_{n_i}, v^{*}_{n_i}) \leq \lim_{i \to \infty}\mathcal{J}_{\varepsilon_{n_i}}(\bar{u}, \bar{v})
    \\
    & = \frac{1}{2}\lim_{i \to \infty} \bigl| [\sqrt{M_\eta} (\bar{\eta}_{\varepsilon_{n_i}}-\eta_\mathrm{ad}), \sqrt{M_\theta} (\bar{\theta}_{\varepsilon_{n_i}}-\theta_\mathrm{ad})] \bigr|_{[\mathscr{H}]^{2}}^{2}
    \\
    &\qquad  + \frac{1}{2}\bigl| [\sqrt{M_u} \bar{u}, \sqrt{M_v} \bar{v}] \bigr|_{[\mathscr{H}]^{2}}^{2}
    \\
    & = \mathcal{J}_{\varepsilon}(\bar{u}, \bar{v}).
\end{align*}
Since the choice of $[\bar{u}, \bar{v}] \in [\mathscr{H}]^{2}$ is arbitrary, we conclude that:
\begin{equation*}
    \mathcal{J}_\varepsilon(u^{**}, v^{**}) = \min_{[u, v] \in [\mathscr{H}]^{2}} \mathcal{J}_{\varepsilon}(u, v),
\end{equation*}
and complete the proof of the item (II-B).
\qed

\section{Proof of Main Theorem \ref{mainTh03}}

This section is devoted to the proof of Main Theorem \ref{mainTh03}. To this end, we need to start with the case of $ \varepsilon > 0 $, and prepare some Lemmas, associated with the G\^{a}teaux differential of the regular cost function $ \mathcal{J}_\varepsilon $. 
\medskip

Let $ \varepsilon > 0 $ be a fixed constant, and let $ [\eta_0, \theta_0] \in V \times V_0 $ be the initial pair, fixed in (A1). Let us take any forcing pair $ [u, v] \in [\mathscr{H}]^2 $, and consider the unique solution $ [\eta, \theta] \in [\mathscr{H}]^2 $ to the state-system (S)$_\varepsilon$. Also, let us take any constant $ \delta \in (-1, 1) \setminus \{0\} $ and any pair of functions $ [h, k] \in [\mathscr{H}]^2 $, and consider another solution $ [\eta^\delta, \theta^\delta] \in [\mathscr{H}]^2 $ to the system (S)$_\varepsilon$, for the initial pair $ [\eta_0, \theta_0] $ and a perturbed forcing pair $ [u +\delta h, v +\delta k] $. On this basis, we consider a sequence of pairs of functions $ \{ [\chi^\delta, \gamma^\delta] \}_{\delta \in (-1, 1) \setminus \{0\}} \subset [\mathscr{H}]^2 $, defined as:
\begin{equation}\label{set00}
    [\chi^\delta, \gamma^\delta] := \left[ \frac{\eta^\delta -\eta}{\delta}, ~ \frac{\theta^\delta -\theta}{\delta} \right] \in [\mathscr{H}]^2, \mbox{ for $ \delta \in (-1, 1) \setminus \{0\} $.}
\end{equation}
This sequence acts a key-role in the computation of G\^{a}teaux differential of the cost function $ \mathcal{J}_\varepsilon $, for $ \varepsilon > 0 $. 
\begin{remark}\label{Rem.GD01}
    Note that for any  $ \delta \in (-1, 1) \setminus \{0\} $, the pair of functions $ [\chi^\delta, \gamma^\delta] \in [\mathscr{H}]^2 $ fulfills the following variational forms:
\begin{align*}
    (\partial_{t} & \chi^\delta(t), \varphi)_{H} + (\partial_{x} \chi^\delta(t), \partial_{x}\varphi)_{H} 
    \\
    & +\int_{\Omega} \left(\int_{0}^{1} g'(\eta(t)+\varsigma \delta \chi^\delta(t)) \, d \varsigma \right) \chi^\delta(t) \varphi \, dx
    \\
    &  + \int_{\Omega} \left(f_{\varepsilon}(\partial_{x} \theta(t)) \int_{0}^{1} \alpha''(\eta(t) +\varsigma \delta \chi^\delta(t)) \, d\varsigma \right) \chi^\delta(t) \varphi \, dx
    \\
    & +\int_{\Omega}\left(\alpha'(\eta^\delta(t)) \int_{0}^{1} f_{\varepsilon}'(\partial_{x}\theta(t) +\varsigma \delta \partial_{x} \gamma^\delta(t)) \, d\varsigma \right) \partial_{x}\gamma^\delta(t) \varphi \, dx
    \\
    = & (M_u h(t), \varphi)_{H}, \mbox{ for any $ \varphi \in V $, a.e. $ t \in (0, T) $, subject to $ \chi^\delta(0) = 0 $ in $ H $,}
\end{align*}
and
\begin{align*}
    (\alpha_{0}(t) & \partial_{t} \gamma^\delta(t), \psi)_{H} + \nu^{2}(\partial_{x} \gamma^\delta(t), \partial_{x}\psi)_{H} 
    \nonumber 
    \\
 & + \int_{\Omega} \left(\alpha(\eta^\delta(t)) \int_{0}^{1} f_{\varepsilon}''(\partial_{x}\theta(t) + \varsigma \delta \partial_{x} \gamma^\delta(t)) \, d \varsigma \right) \partial_{x} \gamma^\delta(t) \partial_{x}\psi \, dx 
    \nonumber 
    \\
& + \int_{\Omega}\left(f'_{\varepsilon}(\partial_{x}\theta(t)) \int_{0}^{1}\alpha'(\eta(t) +\varsigma \delta \chi^\delta(t)) \, d \varsigma \right) \chi^\delta(t) \partial_{x} \psi \, dx 
    \nonumber
    \\
    = & (M_v k(t), \psi)_{H}, \mbox{ for any $ \psi \in V_{0} $,  a.e. $ t \in (0, T) $, subject to $ \gamma^\delta(0) = 0 $ in $ H $.}
\end{align*}
    In fact, these variational forms are obtained by taking the difference between respective two variational forms for $ [\eta^\delta, \theta^\delta] $ and $ [\eta, \theta] $, as in Main Theorem \ref{mainTh01} (I-A), and by using the following linearization formulas: 

\begin{align*}
    & \frac{1}{\delta} \bigl( g(\eta^{\delta})-g(\eta) \bigr) = \left( \int_{0}^{1}g'(\eta + \varsigma \delta \chi^\delta) \, d\varsigma \right) \chi^{\delta} \mbox{ in $ \mathscr{H} $,}
\end{align*}
\begin{align*}
    & \frac{1}{\delta} \bigl( \alpha' (\eta^{\delta})f_{\varepsilon}(\partial_{x}\theta^{\delta}) - \alpha'(\eta)f_{\varepsilon}(\partial_{x}\theta) \bigr)
    \\
    & \qquad = \frac{1}{\delta} \bigl( (\alpha' (\eta^{\delta})-\alpha'(\eta))f_{\varepsilon}(\partial_{x} \theta) \bigr) + \frac{1}{\delta} \bigl( \alpha'(\eta^\delta)(f_{\varepsilon}(\partial_{x}\theta^{\delta})-f_{\varepsilon}(\partial_{x}\theta)) \bigr)
    \\
    & \qquad = \left( f_{\varepsilon}(\partial_{x}\theta) \int_{0}^{1}\alpha''(\eta+\varsigma \delta \chi^\delta) \, d\varsigma \right) \chi^{\delta}
    \\
    &  \qquad \qquad + \left( \alpha'(\eta^{\delta}) \int_{0}^{1}f_{\varepsilon}'(\partial_{x}\theta +\varsigma \delta \partial_{x} \gamma^\delta) \, d\varsigma \right) \partial_{x} \gamma^{\delta} \mbox{ in $ \mathscr{H} $,}
\end{align*}
and
\begin{align*}
    & \frac{1}{\delta} \bigl( \alpha(\eta^{\delta}) f'_{\varepsilon}(\partial_{x}\theta^{\delta}) - \alpha(\eta)f'_{\varepsilon}(\partial_{x}\theta) \bigr)
    \\
    & \qquad = \frac{1}{\delta} \bigl( \alpha(\eta^{\delta})(f'_{\varepsilon}(\partial_{x}\theta^{\delta})-f'_{\varepsilon}(\partial_{x}\theta)) \bigr) + \frac{1}{\delta} \bigl( (\alpha (\eta^{\delta})-\alpha(\eta))f'_{\varepsilon}(\partial_{x}\theta) \bigr)
    \\
    & \qquad = \left( \alpha(\eta^{\delta}) \int_{0}^{1}f_{\varepsilon}''(\partial_{x}\theta +\varsigma \delta \partial_{x} \gamma^\delta) \, d\varsigma \right) \partial_{x}\gamma^{\delta}
    \\
    & \qquad \qquad +\left( f'_{\varepsilon}(\partial_{x}\theta) \int_{0}^{1}\alpha'(\eta +\varsigma \delta \chi^\delta) \, d \varsigma \right) \chi^{\delta} \mbox{ in $ \mathscr{H} $.}
\end{align*}
Incidentally, the above linearization formulas can be verified as consequences of the assumptions (A1)--(A3) and the mean-value theorem (cf. \cite[Theorem 5 in p. 313]{lang1968analysisI}).
\end{remark}

Now, we verify the following two Lemmas. 
\begin{lemma}\label{Lem.GD01}
    Let us fix $ \varepsilon > 0 $, and assume (A1)--(A3). Then, for any $ [u, v] \in [\mathscr{H}]^2 $, the cost function $ \mathcal{J}_\varepsilon $ admits the G\^{a}teaux derivative $ \mathcal{J}_\varepsilon'(u, v) \in [\mathscr{H}]^2  $ $ (= ([\mathscr{H}]^2 )^*) $, such that:
    \begin{align}\label{GD01}
        \bigl( \mathcal{J}_\varepsilon'(u, v), & [h, k] \bigr)_{[\mathscr{H}]^2} = \bigl( [M_\eta (\eta -\eta_\mathrm{ad}), M_\theta (\theta -\theta_\mathrm{ad})], \bar{\mathcal{P}}_\varepsilon [M_u h, M_v k] \bigr)_{[\mathscr{H}]^2}
        \nonumber
        \\
        & +\bigl( [M_u u, M_v v], [h, k] \bigr)_{[\mathscr{H}]^2}, \mbox{ for any $ [h, k] \in [\mathscr{H}]^2 $.}
    \end{align}
    In the context, $ [\eta, \theta] $ is the solution to the state-system (S)$_\varepsilon$, for the initial pair $[\eta_0, \theta_0]$ and forcing pair $ [u, v] $, and $ \bar{\mathcal{P}}_\varepsilon : [\mathscr{H}]^2 \longrightarrow \mathscr{Z} $ is a bounded linear operator, which is given as a restriction $ \mathcal{P}|_{\{[0, 0]\} \times [\mathscr{H}]^2} $ of the (linear) isomorphism $ \mathcal{P} = \mathcal{P}(a, b, \mu, \lambda, \omega, A) : [H]^2 \times \mathscr{Y}^* \longrightarrow \mathscr{Z} $, as in Proposition \ref{ASY_Cor.1}, in the case when:
    \begin{equation}\label{set01}
        \begin{cases}
            [a, b]= [\alpha_0, 0] \mbox{ in $ W^{1, \infty}(Q) \times L^\infty(Q) $,}
            \\
            \mu = \bar{\mu}_\varepsilon := \alpha''(\eta) f_\varepsilon(\partial_x \theta) \mbox{ in $ L^\infty(0, T; H) $,} 
            \\
            [\lambda, \omega, A] = [\bar{\lambda}_\varepsilon, \bar{\omega}_\varepsilon, \bar{A}_\varepsilon] := \bigl[ g'(\eta), \alpha'(\eta) f_\varepsilon'(\partial_x \theta), \alpha(\eta) f_\varepsilon''(\partial_x \theta) \bigr] \mbox{in $ [L^\infty(Q)]^3 $.}
            \end{cases}
    \end{equation}
\end{lemma}
\paragraph{\textbf{Proof.}}{Let us fix any $ [u, v] \in [\mathscr{H}]^2 $, and take any $ \delta \in (-1, 1) \setminus \{0\} $ and any $ [h, k] \in [\mathscr{H}]^2 $. Then, it is easily seen that:
\begin{align}\label{GD02}
    \frac{1}{\delta} \bigl(  \mathcal{J}_\varepsilon &( u +\delta h, v +\delta k) -\mathcal{J}_\varepsilon(u, v) \bigr)
    \nonumber
    \\
    = & \left( \frac{M_\eta}{2} (\eta^\delta +\eta -2 \eta_\mathrm{ad}), \chi^\delta \right)_{\hspace{-0.5ex}\mathscr{H}} +\left( \frac{M_\theta}{2} (\theta^\delta +\theta -2 \theta_\mathrm{ad}), \gamma^\delta \right)_{\hspace{-0.5ex}\mathscr{H}}
    \\
    & \quad +\left( \frac{M_u}{2}(2u +\delta h), h \right)_{\hspace{-0.5ex}\mathscr{H}} +\left( \frac{M_v}{2} (2v +\delta k), k \right)_{\hspace{-0.5ex}\mathscr{H}}.
    \nonumber
\end{align}
Here, let us set:
\begin{subequations}\label{sets02}
\begin{equation}\label{set02}
    \begin{cases}
        \displaystyle \bar{\mu}_\varepsilon^\delta := f_\varepsilon(\partial_x \theta) \int_0^1 \alpha''(\eta +\varsigma \delta \chi^\delta) \, d \varsigma \mbox{ in $ L^\infty(0, T; H) $,} 
        \\
        \displaystyle \bar{\lambda}_\varepsilon^\delta := \int_0^1 g'(\eta +\varsigma \delta \chi^\delta) \, d\varsigma \mbox{ in $ L^\infty(Q) $,}
        \\
        \displaystyle \bar{\omega}_\varepsilon^\delta := \alpha'(\eta^\delta) \int_0^1 f_\varepsilon'(\partial_x \theta +\varsigma \delta \partial_x \gamma^\delta) \, d \varsigma \mbox{ in $ L^\infty(Q) $,}
        \\
        \displaystyle \bar{A}_\varepsilon^\delta := \alpha(\eta^\delta) \int_0^1 f_\varepsilon''(\partial_x \theta +\varsigma \delta \partial_x \gamma^\delta) \, d \varsigma \mbox{ in $ L^\infty(Q) $,}
    \end{cases} 
\end{equation}
and
\begin{align}\label{set03}
    \displaystyle \bar{k}_\varepsilon^\delta := M_v k + \partial_x & \left[ \rule{-1pt}{16pt} \right. 
    \chi^\delta f_\varepsilon'(\partial_x \theta) \int_0^1 \alpha'(\eta +\varsigma \delta \chi^\delta) \, d \varsigma 
    \nonumber
    \\
    & \qquad \displaystyle  -\chi^\delta \alpha'(\eta^\delta) \int_0^1 f_\varepsilon'(\partial_x \theta +\varsigma \delta \partial_x \gamma^\delta) \, d \varsigma 
    \left. \rule{-1pt}{16pt} \right] \mbox{ in $ \mathscr{V}_0^* $,}
    \\
    & \mbox{for all $ \delta \in (-1, 1) \setminus \{0\} $.}
    \nonumber
\end{align}
\end{subequations}
Then, in the light of Remark \ref{Rem.GD01}, one can say that:
\begin{equation*}
    [\chi^\delta, \gamma^\delta] = \bar{\mathcal{P}}_\varepsilon^\delta [M_u h, \bar{k}_\varepsilon^\delta] \mbox{ in $ \mathscr{Z} $, for $ \delta \in (-1, 1) \setminus \{0\} $,}
\end{equation*}
by using the restriction $ \bar{\mathcal{P}}_\varepsilon^\delta := \mathcal{P}|_{\{[0, 0]\} \times \mathscr{Y}^*} : \mathscr{Y}^* \longrightarrow \mathscr{Z} $ of the (linear) isomorphism $ \mathcal{P} = \mathcal{P}(a, b, \mu, \lambda, \omega, A) : [H]^2 \times \mathscr{Y}^* \longrightarrow \mathscr{Z} $, as in Proposition \ref{ASY_Cor.1}, in the case when:
\begin{equation*}
    \begin{cases}
        [a, b, \lambda, \omega, A]= [\alpha_0, 0, \bar{\lambda}_\varepsilon^\delta, \bar{\omega}_\varepsilon^\delta, \bar{A}_\varepsilon^\delta] \mbox{ in $ W^{1, \infty}(Q) \times [L^\infty(Q)]^4 $,}
        \\
        \mu = \bar{\mu}_\varepsilon^\delta \mbox{ in $ L^\infty(0, T; H) $, for $ \delta \in (-1, 1) \setminus \{0\} $.} 
    \end{cases}
\end{equation*}
Besides, taking into account \eqref{f_eps}, \eqref{sets02}, (A2), (A3), and Remarks \ref{Rem.Prelim01} and \ref{Rem.C_0*}, we have:
\begin{subequations}\label{est3-01}
    \begin{align}\label{est3-01a}
    \bar{C}_0^* := & \frac{81 (1 +\nu^2)}{\min\{ 1, \nu^2, \inf \alpha_0(Q) \}}  \left( 1 +|\alpha_0|_{W^{1, \infty}(Q)} +|g'|_{L^\infty(\mathbb{R})} +|\alpha'|_{L^\infty(\mathbb{R})} \right)
        \\
        \geq & \frac{81 (1 +\nu^2)}{\min\{ 1, \nu^2, \inf \alpha_0(Q) \}} \sup_{0 < |\delta| < 1} \left\{ 1 +|\alpha_0|_{W^{1, \infty}(Q)} +|\bar{\lambda}_\varepsilon^\delta|_{L^\infty(Q)} +|\bar{\omega}_\varepsilon^\delta|_{L^\infty(Q)} \right\},
        \nonumber
    \end{align}
and
    \begin{align}\label{est3-01b}
    \bigl| \bigl< [M_u h(t), \, & \bar{k}_\varepsilon^\delta(t)], [\varphi, \psi] \bigr>_{V \times V_0} \bigr| \leq |\langle M_u h(t), \varphi \rangle_{V}| +|\langle \bar{k}_\varepsilon^\delta(t), \psi \rangle_{V_0}| 
    \nonumber
        \\
    \leq \, & |M_u h(t)|_H |\varphi|_H +|M_v k(t)|_H |\psi|_H +2|\alpha'|_{L^\infty(\mathbb{R})}|\chi^\delta(t)|_H|\partial_x \psi|_H
    \nonumber
        \\
        \leq \, & M_u |h(t)|_H |\varphi|_V +\bigl( \sqrt{2} M_v|k(t)|_H +2|\alpha'|_{L^\infty(\mathbb{R})}|\chi^\delta(t)|_H \bigr) |\psi|_{V_0},
        \\
    & \mbox{for a.e. $ t \in (0, T) $, any $ [\varphi, \psi] \in V \times V_0 $, and any $  \delta \in (-1, 1) \setminus \{0\} $,}
    \nonumber
\end{align} 
so that 
    \begin{align}\label{est3-01c}
    \bigl| [M_u h(t), \bar{k}_\varepsilon^\delta(t)] & \bigr|_{V^* \times V_0^*}^2 \leq \bar{C}_1^* \bigl( \bigl| [h(t), k(t)] \bigr|_{[H]^2}^2 +|\chi^\delta(t)|_H^2 \bigr),
    \nonumber
        \\
    & \mbox{for a.e. $ t \in (0, T) $, and any $ \delta \in (-1, 1) \setminus \{ 0 \} $,}
\end{align}
\end{subequations}
with a positive constant $ \bar{C}_1^* := 4 \bigl(  M_u^2 + M_v^2 + |\alpha'|_{L^\infty(\mathbb{R})}^2  \bigr) $.
\medskip

Now, having in mind \eqref{est3-01}, let us apply Proposition \ref{Prop(I-B)} to the case when:
\begin{equation*}
    \begin{array}{c}
    \begin{cases}
        [a^1, b^1, \mu^1, \lambda^1, \omega^1, A^1] = [a^2, b^2, \mu^2, \lambda^2, \omega^2, A^2] = [\alpha_0, 0, \overline{\mu}_\varepsilon^\delta, \bar{\lambda}_\varepsilon^\delta, \bar{\omega}_\varepsilon^\delta, \bar{A}_\varepsilon^\delta],
        \\[0.5ex]
        [p_0^1, z_0^1] = [p_0^2, z_0^2] = [0, 0], ~ [h^1, k^1] = [M_u h, \bar{k}_\varepsilon^\delta], ~ [h^2, k^2] = [0, 0],
        \\[0.5ex]
        [p^1, z^1] = [\chi^\delta, \gamma^\delta] = \bar{\mathcal{P}}_\varepsilon^\delta [M_u h, \bar{k}_\varepsilon^\delta], ~ [p^2, z^2] = [0, 0] = \bar{\mathcal{P}}_\varepsilon^\delta [0, 0],
    \end{cases}
        \ \\
        \ \\[-1.5ex]
    \quad\mbox{for $ \delta \in (-1, 1) \setminus \{0\} $.}
    \end{array}
\end{equation*}
Then, we estimate that:
\begin{align*}
    \frac{d}{dt} & \bigl( |\chi^\delta(t)|_H^2 +|\sqrt{\alpha_0(t)} \gamma^\delta(t)|_H^2 \bigr) +\bigl( |\chi^\delta(t)|_V^2 +\nu^2 |\gamma^\delta(t)|_{V_0}^2 \bigr)
    \\
    & \leq 3 \bar{C}_0^* \bigl( |\chi^\delta(t)|_H^2 +|\sqrt{\alpha_0(t)} \gamma^\delta(t)|_H^2 \bigr) +2\bar{C}_0^* \big( |M_u h(t)|_{V^*}^2 +|\bar{k}_\varepsilon^\delta(t)|_{V_0^*}^2 \bigr)
    \\
    & \leq 3 \bar{C}_0^* (1 +\bar{C}_1^*) \bigl( |\chi^\delta(t)|_H^2 +|\sqrt{\alpha_0(t)} \gamma^\delta(t)|_H^2 \bigr) +2\bar{C}_0^* \bar{C}_1^* \big( |h(t)|_{H}^2 +|k(t)|_{H}^2 \bigr), 
    \\
    & \qquad\mbox{for a.e. $ t \in (0, T) $,}
\end{align*}
and subsequently, by using (A3) and Gronwall's lemma, we observe that:
\begin{itemize}
    \item[\textmd{$(\star\,1)$}]the sequence $ \{ [\chi^\delta, \gamma^\delta] \}_{\delta \in (-1, 1) \setminus \{0\}} $ is bounded in $ [C([0, T]; H)]^2 \cap \mathscr{Y} $. 
\end{itemize}

Meanwhile, as consequences of \eqref{set00}, \eqref{set01}--\eqref{est3-01}, $(\star\,1)$, (A1)--(A3), Main Theorem \ref{mainTh01}, Remark \ref{Rem.mTh01Conv}, and Lebesgue's dominated convergence theorem, one can find a sequence $ \{ \delta_n \}_{n = 1}^\infty \subset \mathbb{R} $, such that:
\begin{subequations}\label{convs3-01}
    \begin{align}\label{convs3-00a}
        0 < |\delta_n| < 1, \mbox{ and } \delta_n \to 0, \mbox{ as $ n \to \infty $,}
    \end{align}
    \begin{align}\label{convs3-01a}
        &
        \begin{cases}
            [\delta_n \chi^{\delta_n}, \delta_n \gamma^{\delta_n}] = [\eta^{\delta_n} -\eta, \theta^{\delta_n} -\theta] \to [0, 0]
            \\
            \qquad \mbox{in $ [C(\overline{Q})]^2 $, and in $ \mathscr{Y} $,}
            \\[1ex]
            [\delta_n \partial_x \chi^{\delta_n}, \delta_n \partial_x \gamma^{\delta_n}] = [\partial_x (\eta^{\delta_n} -\eta), \partial_x (\theta^{\delta_n} -\theta)] \to [0, 0]
            \\
            \qquad \mbox{in $ [\mathscr{H}]^2 $, and in the pointwise sense a.e. in $ Q $,}
        \end{cases}
        \mbox{as $ n \to \infty $,}
    \end{align}
    \begin{align}\label{convs3-01b}
        [\bar{\lambda}_\varepsilon^{\delta_n}, & \bar{\omega}_\varepsilon^{\delta_n}, \bar{A}_\varepsilon^{\delta_n}] \to [\bar{\lambda}_\varepsilon, \bar{\omega}_\varepsilon, \bar{A}_\varepsilon]  \mbox{ weakly-$ * $ in $ [L^\infty(Q)]^3 $,}
        \nonumber
        \\
        & \mbox{and in the pointwise sense a.e. in $ Q $, ~ as $ n \to \infty $,}
    \end{align}
    \begin{equation}\label{convs3-01c}
        \begin{cases}
            \bar{\mu}_\varepsilon^{\delta_n} \to \bar{\mu}_\varepsilon \mbox{ weakly-$ * $ in $ L^\infty(0, T; H) $,}
            \\
            \bar{\mu}_\varepsilon^{\delta_n}(t) \to \bar{\mu}_\varepsilon(t) \mbox{ in $ H $,}
            \mbox{ for a.e. $ t \in (0, T) $,}
        \end{cases}
        \mbox{as $ n \to \infty $,}
    \end{equation}
    and
    \begin{align}\label{convs3-01d}
        \langle \bar{k}_\varepsilon^{\delta_n} -M_v k, \psi \rangle_{\mathscr{V}_0} =  &~ -\left( \chi^{\delta_n}, ~ f_\varepsilon'(\partial_x \theta) \left( \rule{-1pt}{14pt} \right. \int_0^1 \alpha'(\eta +\varsigma \delta_n \chi^{\delta_n}) \, d \varsigma  \left. \rule{-1pt}{14pt} \right) \partial_x \psi\right)_{\hspace{-0.5ex}\mathscr{H}}
        \nonumber
        \\
        &~ +\left( \chi^{\delta_n}, ~ 
        \alpha'(\eta^{\delta_n}) \left( \rule{-1pt}{14pt} \right. \int_0^1 f_\varepsilon'(\partial_x \theta +\varsigma \delta_n \partial_x \gamma^{\delta_n}) \, d \varsigma \left. \rule{-1pt}{14pt} \right) \partial_x \psi \right)_{\hspace{-0.5ex}\mathscr{H}}
        \\
        \to &~ 0, \mbox{ as $ n  \to \infty $.}
        \nonumber
    \end{align}
\end{subequations}

On account of \eqref{set00} and \eqref{set01}--\eqref{convs3-01}, we can apply Proposition \ref{ASY_Cor.2} (B), and can see that:
\begin{align}\label{conv3-02}
    [\chi^{\delta_n}, & \gamma^{\delta_n}] = \bar{\mathcal{P}}_\varepsilon^{\delta_n}[M_u h, \bar{k}_\varepsilon^{\delta_n}] \to [\chi, \gamma] := \bar{\mathcal{P}}_\varepsilon[M_u h, M_v k] \mbox{ in $ [\mathscr{H}]^2 $, weakly in $ \mathscr{Y} $,}
    \nonumber
    \\
    & \mbox{and weakly in $ W^{1, 2}(0, T; V^*) \times W^{1, 2}(0, T; V_0^*) $, as $ n \to \infty $.}
\end{align}
Since the uniqueness of the solution $ [\chi, \gamma] = \bar{\mathcal{P}}_\varepsilon [M_u h, M_v k] $ is guaranteed by Proposition \ref{Prop(I-A)}, the observations \eqref{GD02}, \eqref{convs3-01}, and \eqref{conv3-02} enable us to compute the directional derivative $ D_{[h, k]}\mathcal{J}_\varepsilon(u, v) \in \mathbb{R} $, as follows:
\begin{align*}
    D_{[h, k]} & \mathcal{J}_\varepsilon(u, v) := \lim_{\delta \to 0} \frac{1}{\delta} \bigl( \mathcal{J}_\varepsilon(u +\delta h, v +\delta k) -\mathcal{J}_\varepsilon(u, v) \bigr) 
    \\
    = & \bigl( [M_\eta (\eta -\eta_\mathrm{ad}), M_\theta (\theta -\theta_\mathrm{ad})], \bar{\mathcal{P}}_\varepsilon [M_u h, M_v k] \bigr)_{[\mathscr{H}]^2}
    +\bigl( [M_u u, M_v v], [h, k] \bigr)_{[\mathscr{H}]^2}, 
    \\
    & \quad \mbox{for any $ [u, v] \in [\mathscr{H}]^2 $, and any direction $ [h, k] \in [\mathscr{H}]^2 $.}
\end{align*}
Moreover, with Proposition \ref{ASY_Cor.1} and Riesz's theorem in mind, we deduce the existence of the G\^{a}teaux derivative $ \mathcal{J}_\varepsilon'(u, v) \in ([\mathscr{H}]^2)^* $ $ (= [\mathscr{H}]^2) $ at $ [u, v] \in [\mathscr{H}]^2 $, i.e.:
\begin{equation*}
    \bigl( \mathcal{J}_\varepsilon'(u, v), [h, k] \bigr)_{[\mathscr{H}]^2} =  D_{[h, k]} \mathcal{J}_\varepsilon(u, v), \mbox{ for every $ [u, v], [h, k] \in [\mathscr{H}]^2 $.}
\end{equation*}

Thus, we conclude this lemma with the required property \eqref{GD01}. 
\hfill \qed
}
\begin{lemma}\label{Lem.GD02}
    Under the assumptions (A1)--(A3), let $ [u_\varepsilon^*, v_\varepsilon^*] \in [\mathscr{H}]^2 $ be an optimal control of the problem (OP)$_\varepsilon$, and let $ [\eta_\varepsilon^*, \theta_\varepsilon^*] $ be the solution to the system (S)$_\varepsilon$, for the initial pair $[\eta_0, \theta_0] $ and forcing pair $ [u_\varepsilon^*, v_\varepsilon^*] $. Also, let $ \mathcal{P}_\varepsilon^* : [\mathscr{H}]^2 \longrightarrow \mathscr{Z} $ be the bounded linear operator, defined in Remark \ref{Rem.mTh03}, with the use of the solution $ [\eta_\varepsilon^*, \theta_\varepsilon^*] $. Let $ \mathcal{P}_\varepsilon : [\mathscr{H}]^{2} \longrightarrow \mathscr{Z} $ be a bounded linear operator, which is defined as a restriction $ \mathcal{P}|_{\{[0, 0]\} \times [\mathscr{H}]^2} $ of the linear isomorphism $ \mathcal{P} = \mathcal{P}(a, b, \mu, \lambda, \omega, A) : [H]^2 \times \mathscr{Y}^* \longrightarrow \mathscr{Z} $, as in Proposition \ref{ASY_Cor.1}, in the case when:
    \begin{equation}\label{set_GD02-01}
        \begin{cases}
            [a, b]= [\alpha_0, 0] \mbox{ in $ W^{1, \infty}(Q) \times L^\infty(Q) $,}
            \\
            \mu = \alpha''(\eta_\varepsilon^*) f_\varepsilon(\partial_x \theta_\varepsilon^*) \mbox{ in $ L^\infty(0, T; H) $,} 
            \\
            [\lambda, \omega, A] = \bigl[ g'(\eta_\varepsilon^*), \alpha'(\eta_\varepsilon^*) f_\varepsilon'(\partial_x \theta_\varepsilon^*), \alpha(\eta_\varepsilon^*) f_\varepsilon''(\partial_x \theta_\varepsilon^*) \bigr] \mbox{ in $ [L^\infty(Q)]^3 $.}
            \end{cases}
    \end{equation}
    Then, the operators $ \mathcal{P}_\varepsilon^* $ and $ \mathcal{P}_\varepsilon $  have a conjugate relationship, in the following sense:
    \begin{align*}
        \bigl( \mathcal{P}_\varepsilon^*[u, v] & , [h, k] \bigr)_{[\mathscr{H}]^2} = \bigl( [u, v], \mathcal{P}_\varepsilon [h, k] \bigr)_{[\mathscr{H}]^2}, 
        \\
        & \mbox{for all $ [h, k], [u, v] \in [\mathscr{H}]^2 $.}
    \end{align*}
\end{lemma}
\paragraph{\textbf{Proof.}}{
    Let us fix arbitrary pairs of functions $ [h, k], [u, v] \in [\mathscr{H}]^2 $, and let us put:
    \begin{equation*}
    [\chi_\varepsilon, \gamma_\varepsilon] := \mathcal{P}_\varepsilon [h, k] \quad \mbox{ and } \quad [p_\varepsilon, z_\varepsilon] := \mathcal{P}_\varepsilon^* [u, v], \mbox{ in $ [\mathscr{H}]^2 $.} 
\end{equation*}
Then, invoking Proposition \ref{Prop(I-A)}, and the settings as in \eqref{setRem4} and \eqref{set_GD02-01}, we compute that:
\begin{align*}
    \bigl( \mathcal{P}_\varepsilon^* & [u, v], [h, k] \bigr)_{[\mathscr{H}]^2} = \int_0^T \bigl( p_\varepsilon(t), h(t) \bigr)_{H} \, dt +\int_0^T \bigl( z_\varepsilon(t), k(t) \bigr)_{H} \, dt
    ~~~~~~~~~~~~~~~~~~~~
    \\
    = & \int_0^T \langle h(t), p_\varepsilon(t) \rangle_{V} \, dt +\int_0^T \langle k(t), z_\varepsilon(t) \rangle_{V_0} \, dt
\end{align*}
\begin{align*}
    = & \int_0^T \left[ \rule{-1pt}{16pt} \right. \bigl< \partial_t \chi_\varepsilon(t), p_\varepsilon(t) \bigr>_{V} +\bigl( \partial_x \chi_\varepsilon(t), \partial_x p_\varepsilon(t) \bigr)_H 
    \\
    & \qquad +\bigl( \alpha''(\eta_\varepsilon^*(t)) f_\varepsilon(\partial_x \theta_\varepsilon^*(t)) \chi_\varepsilon(t),  p_\varepsilon(t) \bigr)_H 
    \\
    & \qquad +\bigl( g'(\eta_\varepsilon^*(t)) \chi_\varepsilon(t), p_\varepsilon(t) \bigr)_H +\bigl( \alpha'(\eta_\varepsilon^*(t))f_\varepsilon'(\partial_x \theta_\varepsilon^*(t)) \partial_x \gamma_\varepsilon(t), p_\varepsilon(t) \bigr)_{H} \left. \rule{-1pt}{16pt} \right] \, dt
    \\
    & +\int_0^T \left[ \rule{-1pt}{16pt} \right. \bigl< \alpha_0(t) \partial_t \gamma_\varepsilon(t), z_\varepsilon(t) \bigr>_{V_0} +\bigl( \alpha'(\eta_\varepsilon^*(t)) f_\varepsilon'(\partial_x \theta_\varepsilon^*(t)) \chi_\varepsilon(t), \partial_x z_\varepsilon(t) \bigr)_{H} 
   \\
    & \qquad +\bigl( \alpha(\eta_\varepsilon^*(t)) f_\varepsilon''(\partial_x \theta_\varepsilon^*(t)) \partial_x \gamma_\varepsilon(t), \partial_x z_\varepsilon(t) \bigr)_H +\nu^2 \bigl( \partial_x \gamma_\varepsilon(t), \partial_x z_\varepsilon(t) \bigr)_H \left. \rule{-1pt}{16pt} \right] \, dt
\end{align*}
\begin{align*}
    = & \bigl( p_\varepsilon(T), \chi_\varepsilon(T) \bigr)_H -\bigl( p_\varepsilon(0), \chi_\varepsilon(0) \bigr)_H +\int_0^T \left[ \rule{-1pt}{16pt} \right. \bigl< -\partial_t p_\varepsilon(t), \chi_\varepsilon(t) \bigr>_{V} 
    \\
    & \qquad +\bigl( \partial_x p_\varepsilon(t), \partial_x \chi_\varepsilon(t) \bigr)_H +\bigl( \alpha''(\eta_\varepsilon^*(t)) f_\varepsilon(\partial_x \theta_\varepsilon^*(t)) p_\varepsilon(t), \chi_\varepsilon(t) \bigr)_H 
    \\
& \qquad +\bigl( g'(\eta_\varepsilon^*(t)) p_\varepsilon(t), \chi_\varepsilon(t) \bigr)_H +\bigl( \alpha'(\eta_\varepsilon^*(t))f_\varepsilon'(\partial_x \theta_\varepsilon^*(t)) \partial_x z_\varepsilon(t), \chi_\varepsilon(t) \bigr)_{H} \left. \rule{-1pt}{16pt} \right] \, dt
    \\
    & +\bigl( \alpha_0(T) z_\varepsilon(T), \gamma_\varepsilon(T) \bigr)_H -\bigl( \alpha_0(0) z_\varepsilon(0), \gamma_\varepsilon(0) \bigr)_H 
    \\
    & \qquad +\int_0^T \left[ \rule{-1pt}{16pt} \right. \bigl< -\partial_t \bigl(\alpha_0 z_\varepsilon)(t),  \gamma_\varepsilon(t) \bigr>_{V_0} +\bigl( \alpha'(\eta_\varepsilon^*(t)) f_\varepsilon'(\partial_x \theta_\varepsilon^*(t)) p_\varepsilon(t), \partial_x \gamma_\varepsilon(t) \bigr)_{H} 
    \\
    & \qquad +\bigl( \alpha(\eta_\varepsilon^*(t)) f_\varepsilon''(\partial_x \theta_\varepsilon^*(t)) \partial_x z_\varepsilon(t), \partial_x \gamma_\varepsilon(t) \bigr)_H +\nu^2 \bigl( \partial_x z_\varepsilon(t), \partial_x \gamma_\varepsilon(t) \bigr)_H \left. \rule{-1pt}{16pt} \right] \, dt
    \\
    = & ( u, \chi_\varepsilon )_{\mathscr{H}} +( v, \gamma_\varepsilon )_{\mathscr{H}} = \bigl( [u, v], \mathcal{P}_\varepsilon [h, k]  \bigr)_{[\mathscr{H}]^2}.
\end{align*}
\qed
}
\begin{remark}\label{Rem.GD02}
    Note that the operator $ \mathcal{P}_\varepsilon \in \mathscr{L}([\mathscr{H}]^2; \mathscr{Z}) $, as in Lemma \ref{Lem.GD02}, corresponds to the operator $ \bar{\mathcal{P}}_\varepsilon \in \mathscr{L}([\mathscr{H}]^2; \mathscr{Z}) $, as in the previous Lemma \ref{Lem.GD01}, under the special setting \eqref{set_GD02-01}. 
\end{remark}
\bigskip

Now, we are ready to prove the Main Theorem \ref{mainTh03}.

\paragraph{\textbf{Proof of (III-A) of Main Theorem \ref{mainTh03}.}}{
    Let $ [u_\varepsilon^*, v_\varepsilon^*] \in [\mathscr{H}]^2 $ be the optimal control of (OP)$_\varepsilon$, with the solution $ [\eta_\varepsilon^*, \theta_\varepsilon^*] \in [\mathscr{H}]^2 $ to the system (S)$_\varepsilon$ for the initial pair $[\eta_0, \theta_0]$, as in (A1), and forcing pair $ [u_\varepsilon^*, v_\varepsilon^*] $, and let $ \mathcal{P}_\varepsilon, \mathcal{P}_\varepsilon^* \in \mathscr{L}([\mathscr{H}]^2; \mathscr{Z}) $ be the two operators as in Lemma \ref{Lem.GD02}. Then, on the basis of the previous Lemmas \ref{Lem.GD01} and \ref{Lem.GD02}, Main Theorem \ref{mainTh03} (III-A) will be demonstrated as follows:
    \begin{align*}
        0 = &~  \left(\mathcal{J}_\varepsilon'(u_\varepsilon^*, v_\varepsilon^*), [h, k]\right)_{[\mathscr{H}]^2} = \lim_{\delta \to 0} \frac{1}{\delta} \bigl( \mathcal{J}_\varepsilon(u_\varepsilon^* +\delta h, v_\varepsilon^* +\delta k) -\mathcal{J}_\varepsilon(u_\varepsilon^*, v_\varepsilon^*) \bigr)
        \\
        = &~ \bigl( [M_\eta (\eta_\varepsilon^* -\eta_\mathrm{ad}), M_\theta (\theta_\varepsilon^* -\theta_\mathrm{ad})], \mathcal{P}_\varepsilon [M_u h, M_v k] \bigr)_{[\mathscr{H}]^2} +\bigl( [M_u u_\varepsilon^*, M_v v_\varepsilon^*], [h, k] \bigr)_{[\mathscr{H}]^2}
        \\
        = &~ \bigl( \mathcal{P}_\varepsilon^* [M_\eta (\eta_\varepsilon^* -\eta_\mathrm{ad}), M_\theta (\theta_\varepsilon^* -\theta_\mathrm{ad})], [M_u h, M_v k] \bigr)_{[\mathscr{H}]^2} +\bigl( [M_u u_\varepsilon^*, M_v v_\varepsilon^*], [h, k] \bigr)_{[\mathscr{H}]^2}
        \\
        = &~ \bigl([M_u p_\varepsilon^*, M_v z_\varepsilon^*] , [h, k] \bigr)_{[\mathscr{H}]^2} +\bigl( [M_u u_\varepsilon^*, M_v v_\varepsilon^*], [h, k] \bigr)_{[\mathscr{H}]^2}
        \\
        = &~ \bigl( [M_u( p_\varepsilon^* +u_\varepsilon^*), M_v (z_\varepsilon^* + v_\varepsilon^*)], [h, k] \bigr)_{[\mathscr{H}]^2}, \mbox{ for any $ [h, k] \in [\mathscr{H}]^2 $.}
    \end{align*}
    \qed
}

\paragraph{\textbf{Proof of (III-B) of Main Theorem \ref{mainTh03}.}}{
    Let $[\eta_{0}, \theta_{0}] \in V \times V_{0}$ be the fixed initial pair as in (A1). For any $ \varepsilon > 0 $, let $ [ u_\varepsilon^*, v_\varepsilon^* ] \in [\mathscr{H}]^2 $, $ [\eta_\varepsilon^*, \theta_\varepsilon^*] \in [\mathscr{H}]^2 $, and $ [p_\varepsilon^*, z_\varepsilon^*] \in  \mathscr{Z} $ be as in Main Theorem \ref{mainTh03} (III-A). Then, by Main Theorem \ref{mainTh02} (II-B), we find an optimal control $ [u^\circ, v^\circ] \in [\mathscr{H}]^2 $ of (OP)$_0$, with a zero-convergent sequence $ \{ \varepsilon_n \}_{n = 1}^\infty \subset (0, 1) $, such that:
\begin{subequations}\label{conv3B}
    \begin{equation}\label{conv3B-01}
        [u_n^*, v_n^*] := [u_{\varepsilon_n}^*, v_{\varepsilon_n}^*] \to [u^\circ, v^\circ] \mbox{ weakly in $ [\mathscr{H}]^2 $, as $ n \to \infty $.}
    \end{equation}
    Let $ [\eta^\circ, \theta^\circ] \in [\mathscr{H}]^2 $ be the solution to (S)$_0$, for the initial pair $[\eta_{0}, \theta_{0}]$ and forcing pair $ [u^\circ, v^\circ] $. Then, having in mind Main Theorem \ref{mainTh01} (I-B) and Remark \ref{Rem.mTh01Conv}, we can find a subsequence of $ \{ \varepsilon_n \}_{n=1}^{\infty} $ (not relabeled) and a function $ \nu^\circ \in L^\infty(Q) $, such that:
\begin{align}\label{conv3B-02}
    [\eta_n^*, \theta_n^*] &~ := [\eta_{\varepsilon_n}^*, \theta_{\varepsilon_n}^*] \to [\eta^\circ, \theta^\circ] \mbox{ in $ [C(\overline{Q})]^2 $, in $ \mathscr{Y} $,}
    \nonumber
    \\
     &~ \mbox{and weakly-$*$ in $ L^\infty(0, T; V) \times L^\infty(0, T; V_0) $, as $ n \to \infty $,}
\end{align}
\begin{align}\label{conv3B-10}
    [\partial_x \eta_n,&~  \partial_x \theta_n] \to [\partial_x \eta^\circ, \partial_x \theta^\circ] \mbox{ in $ [\mathscr{H}]^2 $,} 
    \nonumber
    \\
    &~ \mbox{and in the pointwise sense a.e. in $ Q $, as $ n \to \infty $,}
\end{align}
\begin{align}\label{conv3B-03}
    & 
    \begin{cases}
        \mu_n^* := \alpha''(\eta_n^*) f_{\varepsilon_n}(\partial_x \theta_n^*) \to \mu^\circ := \alpha''(\eta^\circ) |\partial_x \theta^\circ| 
        \\
        \quad~ \mbox{weakly-$ * $ in $ L^\infty(0, T; H) $,}
        \\
        \quad~ \mbox{and in the pointwise sense a.e. in $ Q $,}
        \\[1ex]
        \mu_n^*(t) \to \mu^\circ(t) \mbox{ in $ H $,}
        \\
        \quad~ \mbox{and in the pointwise sense for a.e. $ t \in (0, T) $,}
    \end{cases}
    \mbox{as $ n \to \infty $,}
\end{align}
\begin{align}\label{conv3B-04}
    \lambda_n^* &~ := g'(\eta_n^*) \to \lambda^\circ := g'(\eta^\circ) \mbox{ in $ C(\overline{Q}) $, as $ n \to \infty $,}
\end{align}
\begin{align}\label{conv3B-05}
    & 
    \begin{cases}
        f_{\varepsilon_n}'(\partial_x \theta_n^*) \to \nu^\circ \mbox{ weakly-$*$ in $ L^\infty(Q) $, as $ n \to \infty $,}
        \\
        |\nu^\circ| \leq 1 \mbox{ a.e. in $ Q $,}
    \end{cases}
\end{align}
and
\begin{align}\label{conv3B-06}
    \omega_n^* &~ := \alpha'(\eta_n^*) f_{\varepsilon_n}'(\partial_x \theta_n^*) \to \alpha'(\eta^\circ)\nu^\circ \mbox{ weakly-$*$ in $ L^\infty(Q) $, as $ n \to \infty $.} 
\end{align}
\end{subequations}
Besides, from \eqref{conv3B-10}, \eqref{conv3B-05}, Remark \ref{Rem.MG} (Fact\,1) and (Fact\,2), and \cite[Proposition 2.16]{MR0348562}, one can see that:
\begin{equation}\label{conv3B-07}
    \nu^\circ \in \partial f_0(\partial_x \theta^\circ) = \mathrm{Sgn}^1(\partial_x \theta^\circ) \mbox{ a.e. in $ Q $.}
\end{equation}

Next, let us put:
\begin{equation*}
    \begin{cases}
        [p_n^*, z_n^*] := [p_{\varepsilon_n}^*, z_{\varepsilon_n}^*] \mbox{ in $ [\mathscr{H}]^2 $,} 
        \\
        A_n^* := \alpha(\eta_n^*) f_{\varepsilon_n}''(\partial_x \theta_n^*) \mbox{ in $ L^\infty(Q) $,}
    \end{cases}
    \mbox{$ n = 1, 2, 3, \dots $.}
\end{equation*}
Then, from \eqref{Thm.5-00}--\eqref{Thm.5-03}, and \eqref{productrule}, it follows that:
\begin{subequations}\label{[p,z]}
\begin{equation}\label{[p,z]-00}
    [M_u (u_n^* + p_n^*), M_v (v_n^* + z_n^*)] = [0, 0] \mbox{ in $ [\mathscr{H}]^2 $, $ n = 1, 2, 3, \dots $,}
\end{equation}
\begin{align}\label{[p,z]-01}
    \bigl< -\partial_t & p_n^*, \varphi \bigr>_{\mathscr{V}} +\bigl( \partial_x p_n^*, \partial_x \varphi \bigr)_{\mathscr{H}} +\bigl( \mu_n^* p_n^*, \varphi \bigr)_{\mathscr{H}} +\bigl( \lambda_n^* p_n^* +\omega_n^* \partial_x z_n^*, \varphi \bigr)_{\mathscr{H}} 
    \nonumber
    \\
    & = \bigl( M_\eta  (\eta_n^* -\eta_\mathrm{ad}), \varphi \bigr)_{\mathscr{H}}, 
    \mbox{for any $ \varphi \in \mathscr{V} $, $ n = 1, 2, 3, \dots $,}
\end{align}
\begin{align}\label{[p,z]-02}
    \bigl< -\alpha_0 \partial_t & z_n^*, \psi \bigr>_{\mathscr{V}_0} +\bigl( ( -\partial_t \alpha_0) z_n^*, \psi \bigr)_{\mathscr{H}} +\bigl( A_n^* \partial_x z_n^*   +\nu^2 \partial_x z_n^* +\omega_n^* p_n^*, \partial_x \psi \bigr)_{\mathscr{H}}
    \nonumber
    \\
    & = \bigl( M_\theta (\theta_n^* -\theta_\mathrm{ad}), \psi \bigr)_{\mathscr{H}}, \mbox{ for any $ \psi \in \mathscr{V}_0 $, $ n = 1, 2, 3, \dots $,}
\end{align}
and
\begin{equation}\label{[p,z]-03}
    [p_n^*(T), z_n^*(T)] = [0, 0] \mbox{ in $ [H]^2 $, $ n = 1, 2, 3, \dots $.}
\end{equation}
\end{subequations}
Here, invoking the operators $ \mathcal{Q}_{\varepsilon}^* \in \mathscr{L}([\mathscr{H}]^2; \mathscr{Z}) $ and $ \mathcal{R}_T \in \mathscr{L}([\mathscr{H}]^2) $ as in Remark \ref{Rem.mTh03}, we apply Proposition \ref{Prop(I-B)} to the case when:
\begin{equation*}
    \begin{array}{c}
    \begin{cases}
        [a^1, b^1, \mu^1, \lambda^1, \omega^1, A^1] = [a^2, b^2, \mu^2, \lambda^2, \omega^2, A^2] 
        \\
        \hspace{7ex}= \mathcal{R}_T[\alpha_0, -\partial_t \alpha_0, \mu_n^*, \lambda_n^*, \omega_n^*, A_n^*],
        \\[1ex]
        [p_0^1, z_0^1] = [p_0^2, z_0^2] = [0, 0],
        \\[1ex]
        [h^1, k^1] = [\mathcal{R}_T \bigl( M_\eta (\eta_n^* -\eta_\mathrm{ad}) \bigr), \mathcal{R}_T \bigl( M_\theta (\theta_n^* -\theta_\mathrm{ad}) \bigr)], ~ [h^2, k^2] = [0, 0],
        \\[1ex]
        [p^1, z^1] = \mathcal{Q}_{\varepsilon_n}^* \bigl[ \mathcal{R}_T [M_\eta (\eta_n^* -\eta_\mathrm{ad}), M_\theta (\theta_n^* -\theta_\mathrm{ad})] \bigr], 
        \\[0.5ex]
        [p^2, z^2] = [0, 0] = \mathcal{Q}_{\varepsilon_n}^* \bigl[ \mathcal{R}_T [0, 0] \bigr],
    \end{cases}
    \mbox{for $ n \in \mathbb{N} $.}
    \end{array}
\end{equation*}
Then, with use of the constant $ \bar{C}_0^* $ as in \eqref{est3-01a}, we deduced that:
\begin{align}\label{Gronwall01}
    \frac{d}{dt} & \bigl( \bigl| (\mathcal{R}_T p_n^*)(t) \bigr|_H^2 +\bigl| \mathcal{R}_T \bigl( \sqrt{\alpha_0} z_n^* \bigr)(t) \bigr|_H^2 \bigr)
    \nonumber
    \\
    & \qquad +\bigl( \bigl|(\mathcal{R}_T p_n^*)(t)\bigr|_V^2 +\nu^2 \bigl|(\mathcal{R}_T z_n^*)(t)\bigr|_{V_0}^2 \bigr)
    \nonumber
    \\
    & \leq 3 \bar{C}_0^* \bigl( \bigl| (\mathcal{R}_T p_n^*)(t) \bigr|_H^2 +\bigl| \mathcal{R}_T \bigl( \sqrt{\alpha_0} z_n^* \bigr)(t) \bigr|_H^2 \bigr) 
    \nonumber
    \\
    & \qquad +2\bar{C}_0^* \big( \bigl| \mathcal{R}_T \bigl( M_\eta (\eta_n^* -\eta_\mathrm{ad}) \bigr)(t) \bigr|_{V^*}^2 +\bigl|\mathcal{R}_T \bigl( M_\theta (\theta_n^* -\theta_\mathrm{ad}) \bigr) (t) \bigr|_{V_0^*}^2 \bigr), 
    \\
    & \mbox{for a.e. $ t \in (0, T) $, $ n = 1, 2, 3, \dots $.}
    \nonumber
\end{align}
As a consequence of \eqref{est3-01a}, \eqref{conv3B-02}, \eqref{Gronwall01}, (A3), and Gronwall's lemma, it is observed that:
\begin{itemize}
    \item[\textmd{$(\star\,2)$}]the sequence $ \{ [p_n^*, z_n^*] \}_{n = 1}^\infty $ is bounded in $ [C([0, T]; H)]^2 \cap \mathscr{Y} $.
\end{itemize}

Furthermore, from \eqref{emb01}, \eqref{f_eps}, \eqref{conv3B-03}, \eqref{conv3B-04}, \eqref{conv3B-06}, \eqref{[p,z]-01}, \eqref{[p,z]-02}, and (A3), we can derive the following estimates:
\begin{align}\label{[p,z]-10}
    \bigl| \bigl< \partial_t p_n^*, \varphi \bigr>_{\mathscr{V}} \bigr| \leq &~ \bigl| \bigl( \mu_n^* p_n^*, \varphi \bigr)_{\mathscr{H}} \bigr| +\bigl| \bigl( \partial_x p_n^*, \partial_x \varphi \bigr)_{\mathscr{H}} \bigr| 
    \nonumber
    \\
    &~ +\bigl| \bigl( \lambda_n^* p_n^* +\omega_n^* \partial_x z_n^*, \varphi \bigr)_{\mathscr{H}} \bigr| +\bigl| \bigl( M_\eta  (\eta_n^* -\eta_\mathrm{ad}), \varphi \bigr)_{\mathscr{H}} \bigr|
    \\
    \leq &~ C_1^* |\varphi|_{\mathscr{V}}, \mbox{ for any $ \varphi \in \mathscr{V} $, $ n = 1, 2, 3, \dots $,}
    \nonumber
\end{align}
and 
\begin{align}\label{[p,z]-11}
    \bigl| \bigl< -\partial_x & (A_n^* \partial_x z_n^*), \psi \bigr>_{\mathscr{W}_{0}} \bigr| 
    = \bigl| \bigl( A_n^* \partial_x z_n^*, \partial_x \psi \bigr)_{\mathscr{H}} \bigr| 
    \nonumber
    \\
    \leq &~ \bigl| \bigl( \alpha_0 z_n^*, \partial_t \psi \bigr)_{\mathscr{H}} \bigr| +\bigl| \bigl( \nu^2 \partial_x z_n^* +\omega_n^* p_n^*, \partial_x \psi \bigr)_{\mathscr{H}} \bigr| +\bigl| \bigl( M_\theta (\theta_n^* -\theta_\mathrm{ad}), \psi \bigr)_{\mathscr{H}} \bigr|
    \\
    \leq &~ C_2^* |\psi|_{\mathscr{W}_0}, \mbox{ for any $ \psi \in C_\mathrm{c}^\infty(Q) $, $ n = 1, 2, 3, \dots $,}
    \nonumber
\end{align}
with $ n $-independent positive constants:
\begin{equation*}
    C_1^* := 2 \sup_{n \in \mathbb{N}} \left\{ \begin{array}{c}
            (1 +|\mu_n^*|_{L^\infty(0, T; H)} + |\lambda_n^*|_{L^\infty(Q)} + |\omega_n^*|_{L^\infty(Q)})  
        \\[1ex]
        \cdot \bigl( \bigl| [p_n^*, z_n^*] \bigr|_{\mathscr{Y}} +\bigl| M_\eta (\eta_n^* -\eta_\mathrm{ad}) \bigr|_{\mathscr{H}} \bigr) 
    \end{array} \right\} ~(< \infty),
\end{equation*}
and
\begin{equation*}
    C_2^* := 2 \sup_{n \in \mathbb{N}} \left\{ \begin{array}{c}
        (1 +\nu^2 +|\alpha_0|_{L^\infty(Q)} +|\omega_n^*|_{L^\infty(Q)})  
        \\[1ex]
        \cdot \bigl( \bigl| [p_n^*, z_n^*] \bigr|_{\mathscr{Y}} +\bigl| M_\theta (\theta_n^* -\theta_\mathrm{ad}) \bigr|_{\mathscr{H}} \bigr) 
    \end{array} \right\} ~(< \infty),
\end{equation*}
respectively. 
\medskip

Due to \eqref{conv3B-04}--\eqref{conv3B-06}, \eqref{[p,z]-10}, \eqref{[p,z]-11}, $ (\star\,2) $, and the compactness theory of Aubin's type (cf. \cite[Corollary 4]{MR0916688}), we can find subsequences of $ \{ [p_n^*, z_n^*] \}_{n = 1}^\infty \subset \mathscr{Y} $, $ \{ \omega_n^* \partial_x z_n^* \}_{n = 1}^\infty \subset \mathscr{H} $, and $ \{ -\partial_x (A_n^* \partial_x z_n^*) \}_{n = 1}^\infty \subset \mathscr{W}_0^* $ (not relabeled), together with the respective limits $ [p^\circ, z^\circ] \in \mathscr{Y} $, $ \xi^\circ \in \mathscr{H} $, and $ \zeta^\circ \in \mathscr{W}_0^* $, such that:
\begin{subequations}\label{conv3B_a}
\begin{align}\label{conv3B-10'}
    & \begin{cases}
        [p_n^*, z_n^*] \to [p^\circ, z^\circ] \mbox{ weakly in $ \mathscr{Y} $,}
        \\[1ex]
        p_n^* \to p^\circ \mbox{ in $ \mathscr{H} $, weakly in $ W^{1, 2}(0, T; V^*) $,}
        \\
        \quad \mbox{and in the pointwise sense a.e. in $ Q $,}
    \end{cases}
    \mbox{as $ n \to \infty $,}
\end{align}
\begin{align}\label{conv3B-11}
    & \begin{cases}
        \lambda_n^* p_n^* \to \lambda^\circ p^\circ \mbox{ in $ \mathscr{H} $,}
        \\[1ex]
        \omega_n^* p_n^* \to \alpha'(\eta^\circ) \nu^\circ p^\circ \mbox{ weakly in $ \mathscr{H} $,}
    \end{cases}
    \mbox{as $ n \to \infty $,}
\end{align}
\begin{align}\label{conv3B-12}
        \omega_n^* \partial_x z_n^* \to \xi^\circ \mbox{ weakly in $ \mathscr{H} $, as $ n \to \infty $,}
\end{align}
and
\begin{align}\label{conv3B-13}
    -\partial_x (A_n^* \partial_x z_n^*) \to \zeta^\circ \mbox{ weakly in $ \mathscr{W}_0^* $, as $ n \to \infty $.}
\end{align}
\end{subequations}

Now, the properties \eqref{Thm.5-10}--\eqref{Thm.5-12} will be verified through the limiting observations for \eqref{[p,z]-00}--\eqref{[p,z]-03}, as $ n \to \infty $, with use of \eqref{conv3B} and \eqref{conv3B_a}.
\medskip

Thus, we complete the proof. \qed
}

\section{{Appendix}}

The objective of the Appendix is to reorganize the  general theory of nonlinear evolution equation, which enables us to handle the state-systems (S)$_\varepsilon$, for all $ \varepsilon \geq 0 $ in a unified fashion. 
\medskip

In what follows, let $ X $ be an abstract Hilbert space. On this basis, the general theory will be stated by considering two Lemmas, and the proofs will be modified (mixed and reduced) versions of the existing theories, such as \cite{MR0348562,MR2582280,Kenmochi81}. 

\begin{lemma}\label{Lem.CP}
    Let $ \{ \mathcal{A}_0(t) \, | \, t \in [0, T] \} \subset \mathscr{L}(X) $ be a class of time-dependent bounded linear operators, let $ \mathcal{G}_0 : X \longrightarrow X $ be a given nonlinear operator, and let $ \Psi_0 : X \longrightarrow [0, \infty] $ be a non-negative, proper, l.s.c., and convex function, fulfilling the following conditions: 
    \begin{itemize}
        \item[\textmd{(cp.0)}]$ \mathcal{A}_0(t) \in \mathscr{L}(X) $ is positive and selfadjoint, for any $ t \in [0, T] $, and it holds that
            \begin{equation*}
                (\mathcal{A}_0(t) w, w)_X \geq \kappa_0 |w|_X^2 ,\ \mbox{for any}\ w \in X,
            \end{equation*}
            with some constant $\kappa_{0} \in (0, 1)$, independent of $t \in [0, T]$ and $w \in X$.
        \item[\textmd{(cp.1)}]$ \mathcal{A}_0 : [0, T] \longrightarrow \mathscr{L}(X) $ is Lipschitz continuous, so that $ \mathcal{A}_0 $ admits the (strong) time-derivative $ \mathcal{A}_0'(t) \in \mathscr{L}(X) $ a.e. in $ (0, T) $, and  
            \begin{equation*}
                A_T^* := \mathrm{ess} \sup_{\hspace{-3ex}t \in (0, T)} \left\{ \max \{ |\mathcal{A}_0(t)|_{\mathscr{L}(X)}, |\mathcal{A}_0'(t)|_{\mathscr{L}(X)} \} \right\} < \infty;
            \end{equation*}
        \item[\textmd{(cp.2)}]$ \mathcal{G}_0 : X \longrightarrow X $ is a Lipschitz continuous operator with a Lipschitz constant $ L_0 $, and $ \mathcal{G}_0 $ has a $ C^1 $-potential functional $ \widehat{\mathcal{G}}_0 : X \longrightarrow \mathbb{R} $, so that the G\^{a}teaux derivative $ \widehat{\mathcal{G}}_0'(w) \in X^* $ $ (= X) $ at any $ w \in X $ coincides with $ \mathcal{G}_0(w) \in X $;  
        \item[\textmd{(cp.3)}]$ \Psi_0 \geq 0 $ on $ X $, and the sublevel set $ \bigl\{ w \in X \, \bigl| \, \Psi_0(w) \leq r \bigr\} $ is compact in $ X $, for any $ r \geq 0 $.
    \end{itemize}
    Then, for any initial data $ w_0 \in D(\Psi_0) $ and a forcing term $ \mathfrak{f}_0 \in L^2(0, T; X) $, the following Cauchy problem of evolution equation:
    \begin{equation*}
        (\mathrm{CP})~~
        \begin{cases}
            \mathcal{A}_0(t) w'(t) +\partial \Psi_0(w(t)) +\mathcal{G}_0(w(t)) \ni \mathfrak{f}_0(t) \mbox{ in $ X $, ~ $ t \in (0, T) $,}
            \\
            w(0) = w_0 \mbox{ in $ X $;}
        \end{cases}
    \end{equation*}
    admits a unique solution $ w \in L^2(0, T; X) $, in the sense that:
    \begin{equation}\label{cp.s01}
        w \in W^{1, 2}(0, T; X), ~  \Psi_0(w) \in L^\infty(0, T),
    \end{equation}
    and
    \begin{equation}\label{cp.s02}
        \begin{array}{c}
            \displaystyle \bigl( \mathcal{A}_0(t)w'(t) +\mathcal{G}_0(w(t)) -\mathfrak{f}_0(t), w(t) -\varpi \bigr)_X +\Psi_0(w(t)) \leq \Psi_0(\varpi), 
            \\[1ex]
            \mbox{for any $ \varpi \in D(\Psi_0) $, a.e. $ t \in (0, T) $.} 
        \end{array}
    \end{equation}
    Moreover, both $ t \in [0, T] \mapsto \Psi_0(w(t)) \in [0, \infty) $ and $ t \in [0, T] \mapsto \widehat{\mathcal{G}}_0(w(t)) \in \mathbb{R} $ are absolutely continuous functions in time, and
    \begin{equation}\label{cp.s03}
        \begin{array}{c}
            \displaystyle |\mathcal{A}_0(t)^{\frac{1}{2}}w'(t)|_X^2 +\frac{d}{dt} \left( \Psi_0(w(t)) +\widehat{\mathcal{G}}_0(w(t)) \right) = (\mathfrak{f}_0(t), w'(t))_X,
            \\[1ex]
            \mbox{for a.e. $ t \in (0, T) $.}
        \end{array}
    \end{equation}
\end{lemma}
\begin{remark}\label{Rem.CP02}
    Under the assumptions (cp.0) and (cp.1), it is easily verified that:
    \begin{align*}
        \frac{d}{dt} \bigl( \mathcal{A}_0(t) & w(t), \varpi(t) \bigr)_X = \bigl( \mathcal{A}_0(t) w(t), \varpi'(t) \bigr)_X
        \\
        & + \bigl( \mathcal{A}_0'(t) w(t), \varpi(t) \bigr)_X +\bigl( \mathcal{A}_0(t) w'(t), \varpi(t) \bigr)_X, 
        \\
        & \qquad \qquad \mbox{ for a.e. $ t \in (0, T) $, and all $ w, \varpi \in W^{1, 2}(0, T; X) $.}
    \end{align*}
    Additionally, we can identify $ \mathcal{A}_0 \in \mathscr{L}(L^2(0, T; X)) $, and for arbitrary functions $ w, \varpi \in L^2(0, T; X) $ and arbitrary sequences $ \{w_n\}_{n = 1}^\infty, \{ \varpi_n \}_{n = 1}^\infty \subset L^2(0, T; X) $, we can compute that:
    \begin{align*}
        & \bigl( \mathcal{A}_0 w_n, \varpi_n \bigr)_{L^2(0, T; X)} = \bigl( w_n, \mathcal{A}_0 \varpi_n \bigr)_{L^2(0, T; X)}  
        \\
        & \qquad \to \bigl( w, \mathcal{A}_0 \varpi \bigr)_{L^2(0, T; X)} = \bigl( \mathcal{A}_0 w, \varpi \bigr)_{L^2(0, T; X)}, \mbox{ as $ n \to \infty $,}
        \\[1ex]
        & \parbox{14cm}{if $ \varpi_n \to \varpi $ in $ L^2(0, T; X) $, and $ w_n \to w $ weakly in $ L^2(0, T; X) $, 
as $ n \to \infty $.}
    \end{align*}
\end{remark}
\begin{remark}\label{Rem.CP}
    Note that the assumptions (cp.2) and (cp.3) imply that the potential $ \widehat{\mathcal{G}}_0 $ is the so-called $\lambda$-convex functional. More precisely, for every $ L > L_0 $, the functional:
    \begin{align}\label{F_L}
        \widehat{\mathcal{F}}_L : w & \in X \mapsto \widehat{\mathcal{F}}_L(w) := \widehat{\mathcal{G}}_0(w) +L |w|_X^2 +\widehat{C}_0 \in \mathbb{R},
        \nonumber
        \\
        & \mbox{with a constant } \widehat{C}_0 := |\widehat{\mathcal{G}}_0(0)| +\frac{|\mathcal{G}_{0}(0)|_X^2}{2L_0};
    \end{align}
    is nonnegative, strictly convex, and coercive on $ X $. Indeed, from the assumption (cp.2), we immediately see the strictly monotonicity property of the G\^{a}teaux differential $ \widehat{\mathcal{F}}_L' \in  \mathscr{L}(X) $, as follows:
    \begin{align*}
        (\widehat{\mathcal{F}}_L' & (w^1) -\widehat{\mathcal{F}}_L'(w^2), w^1 -w^2)_X = (\mathcal{G}_0(w^1) -\mathcal{G}_0(w^2), w^1 -w^2)_X +2L |w^1 -w^2|_X^2
        \\
        & \geq (2L -L_0) |w^1 -w^2|_X^2 >  0, \mbox{ if $ w^\ell \in X $, $ \ell = 1, 2 $, $w^{1} \ne w^{2}$, and $ L > L_0 $.}
    \end{align*}
    Hence, for every $ L  >  L_0 $, $ \widehat{\mathcal{F}}_L $ is strictly convex on $ X $ (cf. \cite[Theorem B in p. 99]{MR0442824}). 
    Moreover, with use of the mean-value theorem  (cf. \cite[Theorem 5 in p. 313]{lang1968analysisI}), one can verify the non-negativity and coercivity of $ \widehat{\mathcal{F}}_L $ as follows:
    \begin{align*}
        \widehat{\mathcal{F}}_L(w) = ~& \widehat{\mathcal{G}}_0(0) +\left( \int_0^1 \mathcal{G}_0(\varsigma w) \, d \varsigma, w \right)_{\hspace{-0.5ex}X} +\bigl( L|w|_X^2 +\widehat{C}_0 \bigr)
        \nonumber
        \\
        \geq &~ -|\widehat{\mathcal{G}}_0(0)| -L_0 |w|_X^2 \int_0^1 \varsigma \, d \varsigma +\bigl( \mathcal{G}_0(0), w \bigr)_X  +\bigl( L|w|_X^2 +\widehat{C}_0 \bigr)
        \nonumber
        \\
        \geq &~ (L -L_0)|w|_{X}^2 \geq 0, \mbox{ for all $ w \in X $.}
        \nonumber
    \end{align*}
\end{remark}
\paragraph{\textbf{Proof of Lemma \ref{Lem.CP}.}}{
    The existence result for the problem (CP) can be proved by means of standard time-discretization method, applied to the following iteration scheme:
    \begin{equation}\label{scheme01}
        \begin{array}{c}
            \displaystyle \frac{1}{\tau_n}\mathcal{A}_{0, i}(w_i -w_{i -1}) +2L(w_i -w_{i -1}) +\partial \Psi_0(w_i) +\mathcal{G}_0(w_i) \ni \mathfrak{f}_{0, i} \mbox{ in $ X $,}
            \\[1ex]
            \mbox{for $ i = 1, \dots, n $, starting from the initial data $ w_0 \in D(\Psi_0) $.}
        \end{array}
    \end{equation}
    In the context, $ n \in \mathbb{N} $ is a given (large) number, $ \tau_n := T/n $ is the time-step-size, $ \{ t_i \}_{i = 0}^n := \{ i \tau_n \}_{i = 0}^n $ is the partition of the time-interval $ [0, T] $, and
    \begin{equation}\label{F_0i}
        \begin{cases}
            \mathcal{A}_{0, i} := \mathcal{A}_{0}(t_i) \mbox{ in $ \mathscr{L}(X) $, $ i = 0, 1, \dots, n $,}
            \\[1ex]
            \displaystyle \mathfrak{f}_{0, i} := \frac{1}{\tau_n} \int_{t_{i -1}}^{t_i} \mathfrak{f}_{0}(\tau) \, d \tau \mbox{ in $ X $, $ i = 1, \dots, n $.} 
        \end{cases}
    \end{equation}
    Here, let us set:
    \begin{align*}
        & \left\{
            \begin{array}{l}
        \displaystyle [\widehat{w}]_n(t) := \chi_{(-\infty, 0]}(t) w_0 +\sum_{i = 1}^n \chi_{(t_{i -1}, t_i]}(t) \left( w_{i} +\frac{t -t_{i}}{\tau_n}(w_{i} -w_{i -1}) \right) \mbox{ in $X$,}
            \\
                \displaystyle [\overline{w}]_n(t) := \chi_{(-\infty, 0]}(t) w_0 +\sum_{i =1}^n \chi_{(t_{i -1}, t_i]}(t) w_i \mbox{ in $ X $,}
            \end{array}
        \right.
        \\
        & \hspace{20ex}\mbox{for all $ t \in [0, \infty) $, and $ n = 1, 2, 3, \dots $,}
    \end{align*}
    and
    \begin{align*}
        & \left\{
            \begin{array}{l}
                \displaystyle [\overline{\mathcal{A}_{0}}]_{n} := \chi_{(-\infty, 0]}(t) \mathcal{A}_{0, 0} +\sum_{i = 1}^n \chi_{(t_{i -1}, t_i]}(t) \mathcal{A}_{0, i} \mbox{ in $ \mathscr{L}(X) $,} 
                \\
                \displaystyle [\overline{\,\mathfrak{f}_0}]_n(t) := \sum_{i = 1}^n \chi_{(t_{i -1}, t_i]}(t) \, \mathfrak{f}_{0, i} \mbox{ in $ X $,}
            \end{array}
        \right.  
        \\
        & \hspace{14ex}\mbox{for all $ t \in [0, \infty) $, and $ n = 1, 2, 3, \dots $.}
    \end{align*}
    Then, it is easily checked from \eqref{F_0i}, (cp.1), and $ \mathfrak{f}_0 \in L^2(0, T; X) $ that
    \begin{equation}\label{conv02}
        \begin{cases}
            [\overline{\mathcal{A}_{0}}]_n \to \mathcal{A}_0 \mbox{ in $ C([0, T]; \mathscr{L}(X)) $,}
            \\[0.5ex]
            [\overline{\,\mathfrak{f}_0}]_n \to \mathfrak{f}_0 \mbox{ in $ L^2(0, T; X) $,}
        \end{cases}
        \mbox{as $ n \to \infty $.}
    \end{equation}

    Now, let us fix a constant $ L > L_0 $, and take $ n \in \mathbb{N} $ so large to satisfy $ (5 L + A_{T}^{*}) \tau_n < \kappa_0 $ $(< 1)$. Then, the existence and uniqueness of the scheme \eqref{scheme01} will be reduced to those of the minimization problems for the following proper, l.s.c., strictly convex, and coercive functions:
    \begin{align*}
        \varpi \in X \mapsto & \frac{1}{2 \tau_n} |\mathcal{A}_{0, i}^{\frac{1}{2}} (\varpi -w_{i -1})|_X^2 +\Psi_0(\varpi) +\widehat{\mathcal{F}}_{L}(\varpi) 
        \\
        & +L|\varpi -w_{i -1}|_X^2 -L|\varpi|_X^2 -\widehat{C}_0 - (\mathfrak{f}_{0, i}, \varpi)_{X} \in (-\infty, \infty], ~ i = 1, \dots, n.
    \end{align*}
    On this basis, let us multiply the both sides of the scheme \eqref{scheme01} by $ w_i -w_0 $. Then, as a consequence of (cp.0)--(cp.3), Remark \ref{Rem.CP}, and Young's inequality, we infer that:
    \begin{align}\label{1st.est00}
        \frac{1}{2 \tau_n} \bigl( \bigl| & \mathcal{A}_{0, i}^{\frac{1}{2}}(w_i -w_0) \bigr|_X^2 -\bigl| \mathcal{A}_{0, i -1}^{\frac{1}{2}}(w_{i -1} -w_0) \bigr|_X^2  \bigr)
        \nonumber
        \\
        \leq &~ \frac{5L +A_T^*}{\kappa_0} \left( \rule{-1pt}{16pt} \right. \frac{\bigl| \mathcal{A}_{0, i}^{\frac{1}{2}}(w_i -w_0) \bigr|_X^2 +\bigl| \mathcal{A}_{0, i -1}^{\frac{1}{2}}(w_{i -1} -w_0) \bigr|_X^2}{2} \left. \rule{-1pt}{16pt} \right)
        \\
        &~ +\frac{1 +2L^2}{2L} \bigl( |\mathfrak{f}_{0, i}|_X^2 +|w_0|_X^2 +\Psi_0(w_0) +\widehat{\mathcal{F}}_L(w_0) \bigr), \mbox{ for $ i = 1, \dots, n $;}
        \nonumber
    \end{align}
    via the following calculations:
    \begin{align*}
        \left( \rule{-1pt}{16pt} \right. \frac{1}{\tau_n} \mathcal{A}_{0, i} (w_i &\, -w_{i -1}), w_i -w_0 \left. \rule{-1pt}{16pt} \right)_{\hspace{-0.5ex}X} \geq \frac{1}{2 \tau_n} \bigl( \bigl| \mathcal{A}_{0, i}^{\frac{1}{2}}(w_i -w_0) \bigr|_X^2 -\bigl| \mathcal{A}_{0, i}^{\frac{1}{2}}(w_{i -1} -w_0) \bigr|_X^2 \bigr)
        \\
        = &~ \frac{1}{2 \tau_n} \bigl( \bigl| \mathcal{A}_{0, i}^{\frac{1}{2}}(w_i -w_0) \bigr|_X^2 -\bigl| \mathcal{A}_{0, i -1}^{\frac{1}{2}}(w_{i -1} -w_0) \bigr|_X^2 \bigr)
        \\
        &~ \qquad -\frac{1}{2} \left( \frac{1}{\tau_n} (\mathcal{A}_{0, i} -\mathcal{A}_{0, i -1})(w_{i -1} -w_0), w_{i -1} -w_0 \right)_{\hspace{-0.5ex}X}
        \\
        \geq &~ \frac{1}{2 \tau_n} \bigl( \bigl| \mathcal{A}_{0, i}^{\frac{1}{2}}(w_i -w_0) \bigr|_X^2 -\bigl| \mathcal{A}_{0, i -1}^{\frac{1}{2}}(w_{i -1} -w_0) \bigr|_X^2 \bigr) -\frac{A_T^*}{2} |w_{i -1} -w_0|_X^2,
    \end{align*}
    \begin{align}\label{w_i^*}
        & \bigl( w_i^*, w_i -w_0 \bigr)_X \geq \Psi_0(w_i) -\Psi_0(w_0),
        \nonumber
        \\
        \mbox{with } w_i^* := \mathfrak{f}_{0, i} -\frac{1}{\tau_n}&\, \mathcal{A}_{0, i}(w_i -w_{i -1}) -2L(w_i -w_{i -1}) -\mathcal{G}_0(w_i) \in \partial \Psi_0(w_i), 
    \end{align}
    \begin{align*}
        \bigl( 2L(w_i - &\,  w_{i -1}), w_i -w_0 \bigr)_X +\bigl( \mathcal{G}_0(w_i), w_i -w_0 \bigr)_X 
        \\
        = &~ 
        \bigl( \widehat{\mathcal{F}}_L'(w_i), w_i -w_0 \bigr)_X -2L(w_{i -1}, w_i -w_0)_{X}
        \\
        \geq &~ \widehat{\mathcal{F}}_L(w_i) -\widehat{\mathcal{F}}_L(w_0) -2L|w_i -w_0|_X |w_{i -1} -w_0|_X -2L|w_0|_X |w_i -w_0|_X
        \\
        \geq &~ \widehat{\mathcal{F}}_L(w_i) -\widehat{\mathcal{F}}_L(w_0) -2L|w_i -w_0|_X^2 -L |w_{i -1} -w_0|_X^2 -L |w_0|_X^2,
    \end{align*}
    \begin{align*}
        (\mathfrak{f}_{0, i}, w_i -w_0)_X \leq \frac{L}{2}|w_i -w_0|_X^2 +\frac{1}{2L} |\mathfrak{f}_{0, i}|_X^2, 
    \end{align*}
    and
    \begin{align*}
        |w_i -w_0|_X^2 \leq \frac{1}{\kappa_0} \bigl( \mathcal{A}_{0, i}(w_i -w_0), w_i -w_0 \bigr)_X = \frac{1}{\kappa_0} \bigl| \mathcal{A}_{0, i}^{\frac{1}{2}}(w_i -w_0) \bigr|_X^2, \mbox{ for $ i = 1, \dots, n $.}
    \end{align*}
    So, applying the discrete version of Gronwall's lemma (cf. \cite[Section 3.1]{emmrich1999discrete}) to \eqref{1st.est00}, and having in mind \eqref{conv02}, it is observed that:
    \begin{align*}
        \bigl| \mathcal{A}_{0, i}^{\frac{1}{2}} & (w_i -w_0) \bigr|_X^2 
        \\
        \leq &~ \frac{1 +2L^2}{L} e^{\frac{4T(A_T^* +5L)}{\kappa_0}} \left( \rule{-1pt}{14pt} \right. \sup_{n \in \mathbb{N}} \bigl| [\overline{\,\mathfrak{f}_0}]_n \bigr|_{L^2(0, T; X)}^2 +T(|w_0|_X^2 +\Psi_0(w_0) +\widehat{\mathcal{F}}_L(w_0)) \left. \rule{-1pt}{14pt} \right) 
        \\
        =: &~ r_0^* < \infty,  \mbox{ for $ i = 1, \dots, n $,}
    \end{align*}
    and
    \begin{align}\label{1st.est01}
        |w_i|_X^2 \leq &~ 2 \left( |w_0|_X^2 +\frac{1}{\kappa_0} |\mathcal{A}_{0, i}^{\frac{1}{2}}(w_i -w_0)|_X^2 \right)
        \nonumber
        \\
        \leq &~ 2 \left( |w_0|_X^2 +\frac{r_0^*}{\kappa_0} \right) =: r_1^* < \infty, \mbox{ for $ i = 1, \dots, n $.}
    \end{align}
    Additionally, multiplying the both sides of \eqref{scheme01} by $ w_i -w_{i -1} $, and using (cp.0)--(cp.3) and \eqref{1st.est01}, we infer that:
    \begin{align}\label{est01}
        \frac{\kappa_0}{2\tau_n} |w_i -w_{i -1} & |_X^2 + \bigl( \Psi_0(w_i) +\widehat{\mathcal{F}}_L(w_i) \bigr) -\bigl( \Psi_0(w_{i -1}) +\widehat{\mathcal{F}}_L(w_{i -1}) \bigr)
        \nonumber
        \\
        \leq &~ \frac{1 +4L^2}{\kappa_0} \cdot \tau_n \bigl( r_1^* +|\mathfrak{f}_{0, i}|_X^2 \bigr), \mbox{ for $ i = 1, \dots, n $,}
    \end{align}
    via the following calculations:
    \begin{align*}
        \bigl( w_i^*, &\, w_i -w_{i -1} \bigr)_X +2L|w_i -w_{i -1}|_X^2 +\bigl( \mathcal{G}_0(w_i), w_i -w_{i -1} \bigr)_X
        \\
        \geq &~ \Psi_0(w_i) -\Psi_0(w_{i -1}) +\bigl( \widehat{\mathcal{F}}_L'(w_i), w_i -w_{i -1} \bigr)_X -2L(w_{i -1}, w_i -w_{i -1})_X
        \\
        \geq &~ \bigl( \Psi_0(w_i) +\widehat{\mathcal{F}}_L(w_i) \bigr) -\bigl( \Psi_0(w_{i -1}) +\widehat{\mathcal{F}}_L(w_{i -1}) \bigr) -\frac{\kappa_0}{4 \tau_n} |w_i -w_{i -1}|_X^2 -\frac{4L^2}{\kappa_0} \cdot \tau_n r_1^*,
        \\
        & \mbox{with the element $ w_i^* \in \partial \Psi_0(w_i) $, as in \eqref{w_i^*},}
    \end{align*}
    and
    \begin{align*}
        (\mathfrak{f}_{0, i}, w_i -w_{i -1})_X \leq \frac{\kappa_0}{4 \tau_n} |w_i -w_{i -1}|_X^2 +\frac{1}{\kappa_0} \cdot \tau_n |\mathfrak{f}_{0, i}|_X^2, \mbox{ for $ i = 1, \dots, n $.}
    \end{align*}
    So, summing up \eqref{est01}, for $ i = 1, \dots, n $, and invoking \eqref{conv02}, we can derive the following estimate:
    \begin{align*}
        \frac{\kappa_0}{2} \int_0^t \bigl| [\widehat{w}]_n' & (\varsigma) \bigr|_X^2 \, d \varsigma +\Psi_0([\overline{w}]_n(t)) +\widehat{\mathcal{F}}_L([\overline{w}]_n(t)) 
        \nonumber
        \\
        & \leq \Psi_0(w_0) +\widehat{\mathcal{F}}_L(w_0)  +\frac{1+4L^{2}}{\kappa_0}\left(Tr_{1}^{*} + \sup_{n \in \mathbb{N}} \bigl| [\overline{\,\mathfrak{f}_0}]_n\bigr|_{L^2(0, T; X)}^2 \right)
        \\
        & =: r_2^* < \infty, \mbox{ for all $ t \in [0, T] $, and $ n = 1, 2, 3, \dots $.}
        \nonumber
    \end{align*}
    This estimate enable us to say that:
    \begin{enumerate}
        \item[$(\star\,3)$]$ \{ [\widehat{w}]_n \}_{n = 1}^\infty $ is bounded in $ W^{1, 2}(0, T; X) $, and $ \{ [\overline{w}]_n \}_{n = 1}^\infty $ is bounded in $ L^\infty(0, T; X) $;
        \item[$(\star\,4)$]$ \bigl\{ [\overline{w}]_n(t), [\widehat{w}]_n(t) \,\bigl|\, t \in [0, T], ~ n = 1, 2, 3, \dots \bigr\} $ is contained in a compact sublevel set $ \bigl\{ \varpi \in X \, \bigl| \, \Psi_0(\varpi) \leq r_2^* \bigr\} $. 
    \end{enumerate}
    By virtue of $(\star\,3)$ and $(\star\,4)$, we can apply the general theories of compactness, such as Ascoli's and Alaoglu's theorems (cf. \cite[Corollary 4]{MR0916688}, \cite[Section 1.2]{MR932730}, and so on), and we can find a limit function $ w \in W^{1, 2}(0, T; X) $ for some subsequences of $ \{ [\widehat{w}]_n \}_{n = 1}^\infty $ and $ \{ [\overline{w}]_n \}_{n = 1}^\infty $ (not relabeled), such that:
    \begin{subequations}\label{conv01}
    \begin{align}\label{conv01-01}
    [\widehat{w}]_{n} \to w \ & \mbox{in}\ C([0, T; X]), \nonumber \\
                              & \mbox{and weakly in}\ W^{1,2}(0, T; X),\ \mbox{as}\ n \to \infty.
    \end{align}
    Here, having in mind:
    \begin{align*}
    |[\widehat{w}]_{n} - [\overline{w}]_{n}|_{L^{\infty}(0, T; X)} \leq \tau_{n}^{\frac{1}{2}}|[\widehat{w}]_{n}'|_{L^{2}(0, T; X)} \to 0,\ \mbox{as}\ n \to \infty, 
    \end{align*}
    we can also see that
    \begin{align}\label{conv01-02}
    [\overline{w}]_{n} \to w\ \mbox{in}\ L^{\infty}(0, T; X), \ \mbox{as}\ n \to \infty.
    \end{align}
    \end{subequations}
    Taking into account \eqref{scheme01}, \eqref{conv02}, \eqref{conv01}, and (cp.0)--(cp.3), we deduce that:
    \begin{align*}
        \int_I \bigl( \mathcal{A}_0 & (t) w'(t), w(t) -\varpi \bigr)_X \, dt +\int_I \bigl( \mathcal{G}_0(w(t)) -\mathfrak{f}_0(t), w(t) -\varpi \bigr)_X \, dt 
        \\
        & +\int_I \Psi_0(w(t)) \, dt -\int_I \Psi_0(\varpi) \, dt  
        \\
        \leq & \lim_{n \to \infty} \int_I \bigl( [\widehat{w}]_n'(t), [\overline{\mathcal{A}}_0]_n(t) ([\overline{w}]_n(t) -\varpi) \bigr)_X \, dt 
        \\
        & +\lim_{n \to \infty} \tau_n \int_I \bigl( 2L [\widehat{w}]_n'(t), [\overline{w}]_n(t) -\varpi \bigr)_X \, dt 
        \\
        & +\lim_{n \to \infty} \int_I \bigl( \mathcal{G}_0([\overline{w}]_n(t)) -[\overline{\,\mathfrak{f}_0}]_n(t), [\overline{w}]_n(t) -\varpi \bigr)_X \, dt 
        \\
        & +\varliminf_{n \to \infty} \int_I \Psi_0([\overline{w}]_n(t)) \, dt -\int_I \Psi_0(\varpi) \, dt \leq 0, 
        \\
        & \mbox{ for any $ \varpi \in D(\Psi_0) $, and any open interval $ I \subset (0, T) $.}
    \end{align*}
This implies that $ w $ is a solution to the problem (CP). 
\medskip

    Next, for the proof of uniqueness, we suppose that the both $ w^\ell \in L^2(0, T; X) $, $ \ell = 1, 2 $, are solutions to (CP). Then, by virtue of  (cp.0)--(cp.3), it is immediately verified that:
    \begin{subequations}\label{comp}
    \begin{align}\label{comp03}
        \bigl( \mathfrak{f}_0 -\mathcal{A}_0 (w^\ell)' -\mathcal{G}_0(w^\ell) \bigr)(t) & \in \partial \Psi_0(w^\ell(t)) \mbox{ in $ X $,}
        \nonumber
        \\
        &  \mbox{for a.e. $ t \in (0, T) $, $ \ell = 1, 2 $,}
    \end{align}
    \begin{align}\label{comp01}
        \bigl( \mathcal{A}_0(t) & (w^1 -w^2)'(t), (w^1 -w^2)(t) \bigr)_X  
        \nonumber
        \\
        = & \frac{1}{2} \bigl( [\mathcal{A}_0(w^1 -w^2)]'(t), (w^1 -w^2)(t)\bigr)_X 
        \nonumber
        \\
        & \quad -\frac{1}{2}\bigl(\mathcal{A}_0'(t)(w^1 -w^2)(t), (w^1 -w^2)(t)\bigr)_X
        \nonumber
        \\
        & \quad +\frac{1}{2}\bigl( \mathcal{A}_0(t) (w^1 -w^2)(t), (w^1 -w^2)'(t) \bigr)_X
        \nonumber
        \\
        \geq & \frac{1}{2} \frac{d}{dt} |\mathcal{A}_0(t)^\frac{1}{2}(w^1 -w^2)(t)|_X^2 
        \nonumber
        \\
        & \quad -\frac{A_T^*}{2\kappa_0} |\mathcal{A}_0(t)^\frac{1}{2}(w^1 -w^2)(t)|_X^2, \mbox{ for a.e. $ t \in (0, T) $,}
    \end{align}
    and
    \begin{align}\label{comp02}
        \bigl( \mathcal{G}_0(w^1(t)) & -\mathcal{G}_0(w^2(t)), (w^1 -w^2)(t) \bigr)_X 
        \nonumber
        \\
        \geq & -\frac{L_0}{\kappa_0}|\mathcal{A}_0(t)^\frac{1}{2}(w^1 -w^2)(t)|_X^2, \mbox{ for a.e. $ t \in (0, T) $.}
    \end{align}
\end{subequations}
    Hence, the uniqueness for the problem (CP) will be verified via the following Gronwall type estimate:
    \begin{equation*}
        \begin{array}{c}
            \displaystyle \frac{d}{dt} |\mathcal{A}_0(t)^\frac{1}{2}(w^1 -w^2)(t)|_X^2 \leq \frac{A_T^* +2L_0}{\kappa_0} |\mathcal{A}_0(t)^\frac{1}{2}(w^1 -w^2)(t)|_X^2  
            \\[1ex]
            \mbox{for a.e. $ t \in (0, T) $,}
        \end{array}
    \end{equation*}
    that will be obtained by referring to the standard method, i.e.: by taking the difference between two equations, as in \eqref{comp03}; by multiplying the both sides by $ (w^1 -w^2)(t) $; and by applying \eqref{comp01} and \eqref{comp02}, the monotonicity of $ \partial \Psi_0 $ in  $ X \times X $, and the initial condition $ w^1(0) = w^2(0) = w_0 $ in $ X $. 
    \medskip

    Finally, we verify \eqref{cp.s03}. Owing to (cp.2) and \cite[Lemma 3.3]{MR0348562}, one can say that the both functions $ t \in [0, T] \mapsto \Psi_0(w(t)) \in [0, \infty) $ and $ t \in [0, T] \mapsto \widehat{\mathcal{G}}_0(w(t)) \in \mathbb{R} $ are absolutely continuous, and:
    \begin{equation}\label{comp04}
        \frac{d}{dt} \left( \Psi_0(w(t)) +\widehat{\mathcal{G}}_0(w(t)) \right) = \bigl( \mathfrak{f}_0(t) -\mathcal{A}_0(t) w'(t), w'(t) \bigr)_X,\ \mbox{for a.e. $ t \in (0, T) $.}
    \end{equation}
    The equality \eqref{cp.s03} will be obtained as a consequence of \eqref{comp04} and (cp.0).
    \qed
\begin{lemma}\label{Lem.CP02}
    Under the notations $ \mathcal{A}_0 $, $ \mathcal{G}_0 $, and $ \Psi_0 $, and assumptions (cp.0)--(cp.3) as in the previous Lemma \ref{Lem.CP}, let us fix $ w_0 \in D(\Psi_0) $ and $ \mathfrak{f}_0 \in L^2(0, T; X) $, and take the unique solution $ w \in L^2(0, T; X) $ to the Cauchy problem (CP). Let $ \{ \Psi_n \}_{n = 1}^\infty $, $ \{ w_{0, n} \}_{n = 1}^\infty \subset X $, and $ \{ \mathfrak{f}_n \}_{n = 1}^\infty $ be, respectively, a sequence of proper, l.s.c., and convex functions on $ X $, a sequence of initial data in $ X $, and a sequence of forcing terms in $ L^2(0, T; X) $, such that:
    \begin{itemize}
        \item[\textmd{(cp.4)}]$ \Psi_n \geq 0 $ on $ X $, for $ n = 1, 2, 3, \dots $, and the union $ \bigcup_{n = 1}^\infty \bigl\{  w \in X \, \bigl| \, \Psi_n(w) \leq r \bigr\} $ of sublevel sets is relatively compact in $ X $, for any $ r \geq 0 $;
    \item[\textmd{(cp.5)}]$ \Psi_n $ converges to $ \Psi_0 $ on $ X $, in the sense of Mosco, as $ n \to \infty $; 
    \item[\textmd{(cp.6)}]$ \sup_{n \in \mathbb{N}} \Psi_n(w_{0, n}) < \infty $, and $ w_{0, n} \to w_0 $ in $ X $, as $ n \to \infty $;
    \item[\textmd{(cp.7)}]$ \mathfrak{f}_n \to \mathfrak{f}_0 $ weakly in $ L^2(0, T; X) $, as $ n \to \infty $.
    \end{itemize}
    Let $ w_n \in W^{1, 2}(0, T; X)$ be the solution to the Cauchy problem (CP), for the initial data $ w_{0, n} \in D(\Psi_n) $ and forcing term $ \mathfrak{f}_n \in L^2(0, T; X) $. Then, 
    \begin{equation*}
        \begin{array}{c}
            \displaystyle w_n \to w\mbox{ in $ C([0, T]; X) $, weakly in $ W^{1, 2}(0, T; X) $,}
            \\[1ex]
            \displaystyle \int_0^T \Psi_n(w_n(t)) \, dt \to \int_0^T \Psi_0(w(t)) \, dt, 
            \mbox{ as $ n \to \infty $,}
        \end{array}
    \end{equation*}
    and
    \begin{equation*}
        \bigl| \Psi_0(w) \bigr|_{C([0, T])} \leq  \sup_{n \in \mathbb{N}} \, \bigl| \Psi_n(w_n) \bigr|_{C([0, T])} < \infty.
    \end{equation*}
\end{lemma}
\paragraph{\textbf{Proof.}}{
    This Lemma is proved by referring to the method of proof as in \cite[Theorem 2.7.1]{Kenmochi81} (also see \cite[Main Theorem 2]{MR3661429}).
\medskip

First, let us apply \eqref{cp.s03} to the solutions $ w_n $, for $ n = 1, 2, 3, \dots $. Then, we have:
\begin{equation}\label{cp.ineq}
    \begin{array}{c}
        \displaystyle |\mathcal{A}_0(t)^\frac{1}{2} w_n'(t)|_X^2 +\frac{d}{dt} \left( \Psi_n(w_n(t)) +\widehat{\mathcal{G}}_0(w_n(t)) \right) = \bigl( \mathfrak{f}_n(t), w_n'(t) \bigr)_X, 
        \\[1ex]
        \mbox{for a.e. $ t \in (0, T) $, $ n = 1, 2,3, \dots $.}
    \end{array}
\end{equation}
Besides, for simplicity of description, we define:
\begin{equation*}
    \begin{array}{c}
        \displaystyle \widehat{\Psi}_0(\varpi) := \int_0^T \Psi_0(\varpi(t)) \, dt \mbox{ and } \widehat{\Psi}_n(\varpi) := \int_0^T \Psi_n(\varpi(t)) \, dt, ~ n = 1, 2, 3, \dots,
        \\[2ex]
        \mbox{for any $ \varpi \in L^2(0, T; X) $.}
    \end{array}
\end{equation*}
By (cp.5), Remark \ref{Rem.MG} (Fact\,2), and \cite[Proposition 2.16]{MR0348562}, the above $ \widehat{\Psi}_0 $ and $ \widehat{\Psi}_n $, $ n = 1, 2, 3, \dots $, form proper, l.s.c., and convex functions on $ L^2(0, T; X) $, such that:
\begin{subequations}\label{cp.Psis}
\begin{equation}\label{cp.Psi_0}
    \begin{cases}
        \bigl[ w, \mathfrak{f}_0 -\mathcal{A}_0 w' -\mathcal{G}_0(w) \bigr] \in \partial \widehat{\Psi}_0 \mbox{ in $ L^2(0, T; X) \times L^2(0, T; X) $,}
        \\[1ex]
        \bigl[ w_n, \mathfrak{f}_{0, n} -\mathcal{A}_0 w_n' -\mathcal{G}_0(w_n) \bigr] \in \partial \widehat{\Psi}_n \mbox{ in $ L^2(0, T; X) \times L^2(0, T; X) $,}
        \\
        \qquad \mbox{for } n = 1, 2, 3, \dots,
    \end{cases}
\end{equation}
and
\begin{equation}\label{cp.Mosco}
    \widehat{\Psi}_n \to \widehat{\Psi}_0 \mbox{ on $ L^2(0, T; X) $, in the sense of Mosco, as $ n \to \infty $.}  
\end{equation}
\end{subequations}

Next, let us take arbitrary $ t \in [0, T] $, and integrate the both sides of \eqref{cp.ineq} over $ [0, t] $. Then, by using H\"{o}lder's and Young's inequalities, and by applying (cp.0), (cp.2), (cp.6), (cp.7), and the mean-value theorem (cf. \cite[Theorem 5 in p. 313]{lang1968analysisI}), we deduce that:
\begin{align}
    & \frac{\kappa_0}{2} \int_0^t |w_n'(\tau)|_X^2 \, d\tau  +\left( \Psi_n(w_n(t)) +\widehat{\mathcal{G}}_0(w_n(t)) \right) 
    \nonumber
    \\
    & \qquad \leq \left( \Psi_n(w_{0, n}) +\widehat{\mathcal{G}}_0(w_{0, n}) \right)  +\frac{1}{2\kappa_0} \int_0^T |\mathfrak{f}_n(t)|_{X}^2 \, dt
    \nonumber
    \\
    & \qquad \leq \sup_{n \in \mathbb{N}} \left( \rule{0pt}{18pt} \Psi_n(w_{0, n}) +\frac{1}{2 \kappa_0}|\mathfrak{f}_n|_{L^2(0, T; X)}^2 \right.
    \nonumber
    \\
    & \qquad \qquad +|\widehat{\mathcal{G}}_0(0)| +|w_{0, n}|_X \bigl( |\mathcal{G}_0(0)|_X +L_0 |w_{0, n}|_X \bigr) \left. \rule{0pt}{18pt} \right)
    \nonumber
    \\
    & \qquad =: r_3^* < \infty, \mbox{ for all $ t \in [0, T] $, and $ n = 1, 2, 3, \dots $.} 
    \label{cp.est01}
\end{align}
From the above estimate, one can say that: 
\begin{equation*}
    \begin{cases} \hspace{-3ex}
        \parbox{15.0cm}{
            \vspace{-2ex}
            \begin{itemize}
                \item $ \{ w_n \}_{n = 1}^\infty $ is bounded in $ W^{1, 2}(0, T; X) $, and is also bounded in $ C([0, T]; X) $,
            \vspace{-1ex}
        \item $ \bigl\{ w_n(t) \,\bigl|\, t \in [0, T], ~ n = 1 ,2, 3, \dots \bigr\} $ is contained in a relatively compact set $ \bigcup_{n = 1}^\infty \bigl\{ \varpi \in X \,\bigl|\, \Psi_n(\varpi) \leq r_3^* \bigr\} $.
                    \vspace{-2ex}
            \end{itemize}
        }
    \end{cases}
\end{equation*}
Therefore, applying (cp.1)--(cp.7), and the general theories of compactness, such as Ascoli's and Alaoglu's theorems (cf. \cite[Corollary 4]{MR0916688}, \cite[Section 1.2]{MR932730}, and so on), we find a limit function $ \bar{w} \in W^{1,2}(0, T; X) $, with a subsequence of $ \{ w_n \}_{n = 1}^\infty $ (not relabeled), such that:
\begin{subequations}\label{MT2}
\begin{align}\label{MT2.1}
    w_n \to &~ \bar{w} \mbox{ in $ C([0, T]; X) $, weakly in $ W^{1,2}(0, T; X) $,} 
    \nonumber
    \\
    & \mbox{and in particular, $ w_{0, n} = w_n(0) \to w_0 = \bar{w}(0) $, as $ n \to \infty $,} 
\end{align}
\begin{equation}\label{MT2.3}
    \begin{array}{c}
        \displaystyle \mathfrak{f}_n -\mathcal{A}_0 w_n' -\mathcal{G}_0(w_n) \to \mathfrak{f}_0 -\mathcal{A}_0 \bar{w}' -\mathcal{G}_0(\bar{w}) 
        \\[1ex]
        \mbox{weakly in $ L^2(0, T; X) $, as $ n \to \infty $,}
    \end{array}
\end{equation}
and
\begin{align}\label{MT2.2}
    0 \leq \Psi_0(\bar{w}(t)) \leq &~ \liminf_{n \to \infty} \Psi_n(w_n(t)) \leq \sup_{n \in \mathbb{N}} \bigl| \Psi_n(w_n) \bigr|_{C([0, T])} 
    \nonumber
    \\
    \leq &~ r_3^* < \infty, \mbox{ for any $ t \in [0, T] $.}
\end{align}
\end{subequations}

On account of \eqref{cp.Psis}, \eqref{MT2}, and Remark \ref{Rem.MG} (Fact\,1), we can observe that $ \bar{w}$ coincides with the unique solution $ w $ to the problem (CP), and we can conclude this Lemma. \qed 
}


\begin{thebibliography}{99}

\bibitem{MR2033382}
Andreu-Vaillo, F.; Caselles, V.; Maz\'on, J.~M.
\newblock {\em Parabolic quasilinear equations minimizing linear growth
  functionals}, Vol. 223 of {\em Progress in Mathematics}.
\newblock Birkh\"auser Verlag, Basel, 2004.

\bibitem{MR3888633}
Antil, H.; Shirakawa, K.; Yamazaki, N.
\newblock A class of parabolic systems associated with optimal controls of
  grain boundary motions.
\newblock {\em Adv. Math. Sci. Appl.}, {\bfseries 27}(2): 299--336, 2018.

\bibitem{MR0773850}
Attouch, H.
\newblock {\em Variational Convergence for Functions and Operators}.
\newblock Applicable Mathematics Series. Pitman (Advanced Publishing Program),
  Boston, MA, 1984.

\bibitem{MR2582280}
Barbu, V.
\newblock {\em Nonlinear Differential Equations of Monotone Types in {B}anach
  Spaces}.
\newblock Springer Monographs in Mathematics. Springer, New York, 2010.

\bibitem{MR0348562}
Br\'ezis, H.
\newblock {\em Op\'erateurs Maximaux Monotones et Semi-groupes de Contractions
  dans les Espaces de {H}ilbert}.
\newblock North-Holland Publishing Co., Amsterdam-London; American Elsevier
  Publishing Co., Inc., New York, 1973.
\newblock North-Holland Mathematics Studies, No. 5. Notas de Matem\'atica (50).

\bibitem{MR2436794}
Caselles, V.; Chambolle, A.; Moll, S.; Novaga, M.
\newblock A characterization of convex calibrable sets in {$\Bbb R^N$} with
  respect to anisotropic norms.
\newblock {\em Ann. Inst. H. Poincar\'{e} Anal. Non Lin\'{e}aire}, {\bfseries
  25}(4): 803--832, 2008.

\bibitem{MR3661429}
Colli, P.; Gilardi, G.; Nakayashiki, R.; Shirakawa, K.
\newblock A class of quasi-linear {A}llen--{C}ahn type equations with dynamic
  boundary conditions.
\newblock {\em Nonlinear Anal.}, {\bfseries 158}: 32--59, 2017.

\bibitem{emmrich1999discrete}
Emmrich, E.
\newblock Discrete versions of gronwall's lemma and their application to the
  numerical analysis of parabolic problems.
\newblock Technical Report 637, Institute of Mathematics, Technische
  Universit\"{a}t Berlin,
  ``\href{http://www3.math.tu-berlin.de/preprints/files/Preprint-637-1999.pdf}{http://www3.math.tu-berlin.de/preprints/files/Preprint-637-1999.pdf}'',
  1999.

\bibitem{MR2746654}
Giga, M.-H.; Giga, Y.
\newblock Very singular diffusion equations: second and fourth order problems.
\newblock {\em Jpn. J. Ind. Appl. Math.}, {\bfseries 27}(3): 323--345, 2010.

\bibitem{MR1865089}
Giga, M.-H.; Giga, Y.; Kobayashi, R.
\newblock Very singular diffusion equations.
\newblock In {\em Taniguchi {C}onference on {M}athematics {N}ara '98}, Vol.~31
  of {\em Adv. Stud. Pure Math.}, pp.  93--125. Math. Soc. Japan, Tokyo, 2001.

\bibitem{MR2096945}
Giga, Y.; Kashima, Y.; Yamazaki, N.
\newblock Local solvability of a constrained gradient system of total
  variation.
\newblock {\em Abstr. Appl. Anal.}, (8): 651--682, 2004.

\bibitem{MR3951294}
Hoppe, R.~H.~W.; Winkle, J.~J.
\newblock A splitting scheme for the numerical solution of the {KWC} system.
\newblock {\em Numer. Math. Theory Methods Appl.}, {\bfseries 12}(3): 661--680,
  2019.

\bibitem{MR2469586}
Ito, A.; Kenmochi, N.; Yamazaki, N.
\newblock A phase-field model of grain boundary motion.
\newblock {\em Appl. Math.}, {\bfseries 53}(5): 433--454, 2008.

\bibitem{MR2548486}
Ito, A.; Kenmochi, N.; Yamazaki, N.
\newblock Weak solutions of grain boundary motion model with singularity.
\newblock {\em Rend. Mat. Appl. (7)}, {\bfseries 29}(1): 51--63, 2009.

\bibitem{MR2836555}
Ito, A.; Kenmochi, N.; Yamazaki, N.
\newblock Global solvability of a model for grain boundary motion with
  constraint.
\newblock {\em Discrete Contin. Dyn. Syst. Ser. S}, {\bfseries 5}(1): 127--146,
  2012.

\bibitem{Kenmochi81}
Kenmochi, N.
\newblock Solvability of nonlinear evolution equations with time-dependent
  constraints and applications.
\newblock {\em Bull. Fac. Education, Chiba Univ.
  (\url{http://ci.nii.ac.jp/naid/110004715232})}, {\bfseries 30}: 1--87, 1981.

\bibitem{MR2668289}
Kenmochi, N.; Yamazaki, N.
\newblock Large-time behavior of solutions to a phase-field model of grain
  boundary motion with constraint.
\newblock In {\em Current advances in nonlinear analysis and related topics},
  Vol.~32 of {\em GAKUTO Internat. Ser. Math. Sci. Appl.}, pp.  389--403.
  Gakk\=otosho, Tokyo, 2010.

\bibitem{MR1712447}
Kobayashi, R.; Giga, Y.
\newblock Equations with singular diffusivity.
\newblock {\em J. Statist. Phys.}, {\bfseries 95}(5-6): 1187--1220, 1999.

\bibitem{MR1752970}
Kobayashi, R.; Warren, J.~A.; Carter, W.~C.
\newblock A continuum model of grain boundaries.
\newblock {\em Phys. D}, {\bfseries 140}(1-2): 141--150, 2000.

\bibitem{MR1794359}
Kobayashi, R.; Warren, J.~A.; Carter, W.~C.
\newblock Grain boundary model and singular diffusivity.
\newblock In {\em Free boundary problems: theory and applications, {II}
  ({C}hiba, 1999)}, Vol.~14 of {\em GAKUTO Internat. Ser. Math. Sci. Appl.},
  pp.  283--294. Gakk\=otosho, Tokyo, 2000.

\bibitem{lang1968analysisI}
Lang, S.
\newblock {\em Analysis I}.
\newblock Addison-Wesley Publishing Company, 1968.

\bibitem{MR3268865}
Moll, S.; Shirakawa, K.
\newblock Existence of solutions to the {K}obayashi--{W}arren--{C}arter system.
\newblock {\em Calc. Var. Partial Differential Equations}, {\bfseries 51}(3-4):
  621--656, 2014.

\bibitem{MR3670006}
Moll, S.; Shirakawa, K.; Watanabe, H.
\newblock Energy dissipative solutions to the {K}obayashi--{W}arren--{C}arter
  system.
\newblock {\em Nonlinearity}, {\bfseries 30}(7): 2752--2784, 2017.

\bibitem{MR0298508}
Mosco, U.
\newblock Convergence of convex sets and of solutions of variational
  inequalities.
\newblock {\em Advances in Math.}, {\bfseries 3}: 510--585, 1969.

\bibitem{MR3888636}
Nakayashiki, R.
\newblock Quasilinear type {K}obayaski-{W}arren-{C}arter system including
  dynamic boundary condition.
\newblock {\em Adv. Math. Sci. Appl.}, {\bfseries 27}(2): 403--437, 2018.

\bibitem{MR2836557}
Ohtsuka, T.; Shirakawa, K.; Yamazaki, N.
\newblock Optimal control problem for {A}llen-{C}ahn type equation associated
  with total variation energy.
\newblock {\em Discrete Contin. Dyn. Syst. Ser. S}, {\bfseries 5}(1): 159--181,
  2012.

\bibitem{MR0442824}
Roberts, A.~W.; Varberg, D.~E.
\newblock {\em Convex functions}.
\newblock Academic Press [A subsidiary of Harcourt Brace Jovanovich,
  Publishers], New York-London, 1973.
\newblock Pure and Applied Mathematics, Vol. 57.

\bibitem{MR2223383}
Shirakawa, K.
\newblock Stability for phase field systems involving indefinite surface
  tension coefficients.
\newblock In {\em Dissipative phase transitions}, Vol.~71 of {\em Ser. Adv.
  Math. Appl. Sci.}, pp.  269--288. World Sci. Publ., Hackensack, NJ, 2006.

\bibitem{MR2101878}
Shirakawa, K.; Kimura, M.
\newblock Stability analysis for {A}llen-{C}ahn type equation associated with
  the total variation energy.
\newblock {\em Nonlinear Anal.}, {\bfseries 60}(2): 257--282, 2005.

\bibitem{MR3082861}
Shirakawa, K.; Watanabe, H.
\newblock Energy-dissipative solution to a one-dimensional phase field model of
  grain boundary motion.
\newblock {\em Discrete Contin. Dyn. Syst. Ser. S}, {\bfseries 7}(1): 139--159,
  2014.

\bibitem{MR3462536}
Shirakawa, K.; Watanabe, H.
\newblock Large-time behavior for a {PDE} model of isothermal grain boundary
  motion with a constraint.
\newblock {\em Discrete Contin. Dyn. Syst.}, {\bfseries 1}(Dynamical systems,
  differential equations and applications. 10th AIMS Conference. Suppl.):
  1009--1018, 2015.
  
\bibitem{MR3038131}
Shirakawa, K.; Watanabe, H.; Yamazaki, N.
\newblock Solvability of one-dimensional phase field systems associated with
  grain boundary motion.
\newblock {\em Math. Ann.}, {\bfseries 356}(1): 301--330, 2013.

\bibitem{MR3362773}
Shirakawa, K.; Watanabe, H.; Yamazaki, N.
\newblock Phase-field systems for grain boundary motions under isothermal
  solidifications.
\newblock {\em Adv. Math. Sci. Appl.}, {\bfseries 24}(2): 353--400, 2014.

\bibitem{MR0916688}
Simon, J.
\newblock Compact sets in the space {$L^p(0,T;B)$}.
\newblock {\em Ann. Mat. Pura Appl. (4)}, {\bfseries 146}: 65--96, 1987.

\bibitem{MR932730}
Vrabie, I.~I.
\newblock {\em Compactness methods for nonlinear evolutions}, Vol.~32 of {\em
  Pitman Monographs and Surveys in Pure and Applied Mathematics}.
\newblock Longman Scientific \& Technical, Harlow; John Wiley \& Sons, Inc.,
  New York, 1987.
\newblock With a foreword by A. Pazy.

\bibitem{MR3203495}
Watanabe, H.; Shirakawa, K.
\newblock Qualitative properties of a one-dimensional phase-field system
  associated with grain boundary.
\newblock In {\em Nonlinear analysis in interdisciplinary
  sciences---modellings, theory and simulations}, Vol.~36 of {\em GAKUTO
  Internat. Ser. Math. Sci. Appl.}, pp.  301--328. Gakk\=otosho, Tokyo, 2013.

\bibitem{MR3238848}
Watanabe, H.; Shirakawa, K.
\newblock Stability for approximation methods of the one-dimensional
  {K}obayashi-{W}arren-{C}arter system.
\newblock {\em Math. Bohem.}, {\bfseries 139}(2): 381--389, 2014.

\end{thebibliography}

\end{document}